\documentclass[10pt,a4paper]{amsart}
\usepackage{amsfonts,layout}
\newcount\minute
\newcount\hour
\def\currenttime{%
	\minute\time
	\hour\minute
	\divide\hour60
	\the\hour:\multiply\hour60\advance\minute-\hour\the\minute}
\def\draftnote{{\it \today \quad  \currenttime \hfill  tex-file :   \jobname}}

\catcode`@=11
\@addtoreset{equation}{section}

\catcode`@=12
\newtheorem{Theorem}{Theorem}[section]
\newtheorem{Definition}{Definition}[section]
\newtheorem{Proposition}{Proposition}[section]
\newtheorem{Lemma}{Lemma}[section]

\newtheorem{Corollary}{Corollary}[section]

\newtheorem{Remark}{Remark}[section]

\newtheorem{Hyp.}{Hyp.}[section]

\begin{document}

\title[]{$L^p$ asymptotic stability of 1D damped wave equation with nonlinear distributed damping}

\author{Y. Chitour} 
\address{Laboratoire des Signaux et Syst\`emes, Universit\'e Paris Saclay, CNRS and Centrale Supelec, France}
\email{yacine.chitour@l2s.centralesupelec.fr}

\author{M. Kafnemer} 
\address{}
\email{meryem.kafnemer@centralesupelec.fr}

\author{P. Martinez} 
\address{Institut de Math\'ematiques de Toulouse; UMR 5219, Universit\'e de Toulouse; CNRS \\ 
UPS IMT F-31062 Toulouse Cedex 9, France} 
\email{patrick.martinez@math.univ-toulouse.fr}

\author{B. Mebkhout} 
\address{D\'epartement de Math\'ematiques, Universit\'e Abou Bekr Belkaid, Tlemcen, Algeria}
\email{b\_mebkhout@mail.univ-tlemcen.dz}


\thanks{This research was partially supported by the iCODE Institute, a research project of the IDEX Paris-Saclay, and by the Hadamard Mathematics LabEx (LMH) through the grant number ANR-11-LABX-0056-LMH in the ``Programme des Investissements d'Avenir''.}

\begin{abstract}

In this paper, we study the one-dimensional wave equation with localized nonlinear damping and Dirichlet boundary conditions, in the $L^p$ framework, with $p\in [1,\infty)$.

We start by addressing the well-posedness problem. We prove the existence and the uniqueness of weak and strong solutions for $p\in [1,\infty)$, under suitable assumptions on the damping function.

Then we study the asymptotic behaviour of the associated energy when $p \in (1,\infty)$, and we provide decay estimates that appear to be almost optimal as compared to a similar problem with boundary damping.

Our study is motivated by earlier works, in particular, \cite{Haraux2009, Chitour-Marx-Prieur-2020}. 
Our proofs combine arguments from \cite{KMJC2022} (wave equation in the $L^p$ framework with a linear damping) with 
a technique of weighted energy estimates (\cite{PM-COCV}) and new integral inequalities
when $p>2$, and with convex analysis tools when $p\in (1,2)$.
		
\end{abstract}
\maketitle


\section{Introduction}


In this paper, we consider the one-dimensional wave equation with nonlinear and localized damping, in a functional framework relying on $L^p$ spaces, with $p\in [1,+\infty)$. The problem is the following:
\begin{equation}
\label{eq-wave}
\begin{cases}
z_{tt} - z_{xx} + a(x) g(z_t) = 0 \quad  &\text{ for } (t,x) \in \Bbb R_+ \times (0,1), \\
z(t,0)=0=z(t,1) \quad  &\text{ for } t \in \Bbb R_+, \\
z(0,x)=z_0(x), z_t(0,x)=z_1 (x) \quad  &\text{ for } x \in  (0,1) ,
\end{cases}
\end{equation}
where the initial data $(z_0,z_1)$ belong to an $L^p$ functional space (defined later), with $p\in [1,\infty)$. The function $a$ is a continuous nonnegative function on $[0,1]$, bounded from below by a positive constant 
on some non-empty open interval $\omega$ of $(0,1)$, which represents the region of the domain where the damping term is active. The nonlinearity $g: \mathbb{R} \mapsto \mathbb{R}$ is a $C^1$ nondecreasing function such that $g(0)=0$ and
\begin{equation}
\label{cond-diss}
\forall x\in \mathbb{R}, \quad x \, g(x) \geq 0 ,
\end{equation}
the classical condition ensuring that the usual energy is nonincreasing in the $L^2$ framework.

For many years, the asymptotic behaviour of the solutions of \eqref{eq-wave} has been investigated in the classical functional space $H^1 _0 (0,1) \times L^2 (0,1)$. Our paper is devoted to the case of the following functional spaces
\begin{align}
\label{X_p}
X_p&:= W^{1,p}_0(0,1) \times L^p(0,1), \\
\label{Y_p}
Y_p&:= \left(W^{2,p}(0,1) \cap W^{1,p}_0(0,1)\right) \times W^{1,p}_0(0,1),
\end{align}
both equipped with their natural norm, and with $p\in [1,+\infty)$.

To the best of our knowledge, very few is known about the global asymptotic stability of \eqref{eq-wave} in those functional spaces. This is partly due to the fact that for linear wave equations on $\Bbb R^n$, $n\geq 2$, the ``semigroup'' associated with the d'Alembertian operator $\square z$ (equal to $z_{tt} - \Delta z$) is not, in general, a well-defined bounded operator for an $L^p$ space framework when $p\neq 2$, see Peral \cite{Peral1980}. However Haraux \cite{Haraux2009} succeeded to study the problem in the one dimensional case, and this is one of the motivations of our work.

Let us give our main results, and then we will compare them to the already existing results on this subject.


\section{Main results}

\subsection{Assumptions} \hfill

We will work under the following assumptions:

	$\mathbf{ (H_1)}$ $ a : [0,1] \rightarrow \mathbb{R}_+\ $ is continuous, and satisfies
	\begin{equation} \label{a0}
	\exists \ a_0 >0, \ a \geq a_0\ \ \hbox{on}\ \ \omega= (b,c) \subset [0,1] .
	\end{equation}
	
$\mathbf{ (H_2)}$
$g: \mathbb{R} \to \mathbb{R}$ is a $C^1$ nondecreasing function such that $g(0)=0$,  and
\begin{align}\label{gxx}
\forall x\neq 0, \quad x g(x) > 0.
\end{align}
(Note that hypothesis $\mathbf{ (H_2)}$ implies that $g'(0)\geq 0$.)


$\mathbf{ (H_3)}$
$g: \mathbb{R} \to \mathbb{R}$ is a globally Lipschitz function, \eqref{gxx} holds true and there exists an increasing and odd function $g_0: [-1,1] \to \Bbb R$ and $0<\alpha \leq \beta$ such that
\begin{equation}
\label{eq-comp-g}
\begin{cases}
\forall s \in [-1,1], \quad \vert g_0 ( s) \vert  \leq \vert g(s) \vert \leq \beta \vert s \vert , \\
\forall \vert s \vert \geq 1, \quad \alpha \vert s \vert \leq \vert g(s) \vert \leq \beta \vert s \vert .
\end{cases}
\end{equation}


\subsection{Well posedness: Weak/strong solutions of Problem \eqref{eq-wave}} \hfill

\subsubsection{Functional framework} \hfill

Given $p\in [1,+\infty)$, the space $X_p$ (defined in \eqref{X_p}) is equipped with the norm
\begin{align}
\Vert{(u,v)}\Vert_{X_p}:&= \left( \frac{1}{p}\int_0^1 \vert u' \vert ^p + \vert v \vert ^p \,  dx \right)^\frac{1}{p}.
\end{align}
which is equivalent to the natural norm $\Vert u' \Vert _{L^p} + \Vert v \Vert _{L^p}$.

Given $p\in [1,+\infty)$, the space $Y_p$ (defined in \eqref{Y_p}) is equipped with the norm
\begin{align}
\Vert(u,v)\Vert_{Y_p}:&= \left(\frac{1}{p}\int_0^1 \vert u'' \vert ^p + \vert v' \vert ^p  \, dx \right)^\frac{1}{p} .
\end{align}
which is equivalent to the natural norm $\Vert u'' \Vert _{L^p} + \Vert v' \Vert _{L^p}$.



\subsubsection{Weak/Strong solutions of \eqref{eq-wave}} \hfill 

We consider the following definitions:

\begin{Definition}(Weak solutions of Problem \eqref{eq-wave})

Given $p\in [1,+\infty)$, and $(z_0,z_1) \in X_p$, the function $$z\in L^\infty(\mathbb{R_+}, W^{1,p}_0(0,1)) \cap W^{1,\infty}(\mathbb{R_+},L^p(0,1)) $$ 
is said to be a weak solution of Problem \eqref{eq-wave} if it satisfies the problem in the dual sense (meaning that the equalities are taken in the weak topology of $X_p$).
\end{Definition}

\begin{Definition}(Strong solutions of Problem \eqref{eq-wave})		
		 
Given $p\in [1,+\infty)$, and  $(z_0,z_1) \in Y_p$, the function 
$$z\in L^\infty(\mathbb{R_+},W^{2,p}(0,1)\cap W^{1,p}_0(0,1))\cap W^{1,\infty}(\mathbb{R_+},W^{1,p}_0(0,1))$$ 
is said to be a strong solution of Problem \eqref{eq-wave} if it satisfies the problem in the classical sense.

\end{Definition}

Then we prove the following results, that complete the studies in \cite{Haraux2009, Chitour-Marx-Prieur-2020} in the $L^\infty$ framework:

\begin{Theorem}{Weak solutions for $p\in [2,+\infty)$.}
\label{wellposedness-w2}

Assume that Hypotheses $\mathbf{ (H_1)}$ and $\mathbf{ (H_2)}$ are satisfied, and take $2\leq p <\infty $. Then for every initial conditions $(z_0,z_1)\in X_p$, there exists a unique weak solution $z$ of Problem \eqref{eq-wave}, satisfying additionnally
	\begin{align}\label{eq:exist-weak}
z_t\in L^\infty(\mathbb{R_+},L^p(0,1)).
	\end{align}
\end{Theorem}

\begin{Theorem}{Strong solutions for $p\in [1,+\infty)$.}
\label{wellposedness-s1}

Assume that Hypotheses $\mathbf{ (H_1)}$ and $\mathbf{ (H_2)}$ are satisfied, and take $1\leq p <\infty $.
Then for every $(z_0,z_1)\in Y_p$, there exists a unique strong solution $z$ of Problem \eqref{eq-wave}, satisfying additionnally
	\begin{align}\label{eq:exist-strong}
z_t\in L^\infty(\mathbb{R_+},W^{1,p}_0(0,1)).
	\end{align}
	\end{Theorem}
	
\begin{Theorem}{Weak solutions for $p\in [1,+\infty)$.}	
\label{wellposedness-w1}

Assume that Hypotheses $\mathbf{ (H_1)}$ and $\mathbf{ (H_3)}$ are satisfied, and take $1\leq p <\infty $. Then, 
	for every initial conditions $(z_0,z_1)\in X_p$, there exists a unique weak solution $z$ verifying \eqref{eq:exist-weak}.
	\end{Theorem}
	
\noindent \begin{Remark} {\rm 
The argument based on D'Alembert formula and fixed point theory that is used in \cite{Chitour-Marx-Prieur-2020} cannot be used in the nonlinear case without imposing the extra assumption that $g$ is linearly bounded.
} 
\end{Remark}


\subsection{Stability results} \hfill

In the classical $L^2$ framework, the natural energy is defined by
$$ \mathcal E_2 (t) := \frac{1}{2} \int _0 ^1 \vert z_x \vert ^2 + \vert z_t \vert ^2 \, dx ,$$
and \eqref{cond-diss} implies that $\mathcal E_2$ is nonincreasing along solutions of \eqref{eq-wave}.
We will recall in section \ref{sec-litt} the main stability results in this context.
 
In the $L^p$ framework, the natural energy is defined by
\begin{equation}
\label{energyp-init}
\mathcal E_p (t) := \Vert (z,z_t) (t) \Vert _{X_p} ^p = \frac{1}{p} \int _0 ^1 \vert z_x \vert ^p + \vert z_t \vert ^p \, dx ,
\end{equation}
and is equivalent to the following one, which will reveal to be more useful:

\begin{equation}
\label{energyp}
E_p (t) := \frac{1}{p} \int _0 ^1 \vert z_x + z_t \vert ^p + \vert z_x - z_t \vert ^p \, dx .
\end{equation}
It is well known (\cite{Haraux2009}) that $E_p$ is nonincreasing along solutions of \eqref{eq-wave}. The main stability result of this paper is to estimate the decay of that energy for a large class of damping functions $g$. 

\subsubsection{Stability results when $p\in (2,+\infty)$} \hfill

When $p\in (2,+\infty)$, we prove the following.

\begin{Theorem}
 \label{thm-stab}
 
Fix $p\in (2,+\infty)$.  Assume that \eqref{a0} and \eqref{eq-comp-g} hold. Then we have the following:
 
 (i) If $g'(0) >0$, the energy $E_p$ decays exponentially to $0$, i.e., there exists $C>0$, $\omega  >0$ such that, for every 
 initial condition $(z_0,z_1) \in X_p$, we have
 \begin{equation}
\label{decay-exp}
\forall t\geq 0, \quad E_p (t) \leq C E_p (0) \, e^{-\omega t} .
\end{equation}
 
 (ii) If $g'(0)=0$, assume that there exists $\eta >0$ such that the function 
 $$H(s) := \frac{g_0(s)}{s}$$ is increasing on $[0,\eta]$. Then, for every
 initial condition $(z_0,z_1) \in X_p$ and $\delta >0$, there exists $C_\delta (E_p(0))$ such that
 \begin{equation}
\label{decay-gene}
\forall t\geq 1, \quad E_p (t) \leq C_\delta (E_p (0)) \, \Bigl( g_0 ^{-1} (\frac{1}{t}) \Bigr) ^{p-\delta} .
\end{equation} 
Moreover, fix $m\geq 1$, and denote 
$$t_m := \frac{\frac{2}{\eta ^{1/m}}}{H(\frac{\eta}{2^m})} .$$ 
Then there exists $C_{m,\delta} (E_p (0))$ such that
\begin{equation}
\label{decay-gene-m}
\forall t\geq t_m , \quad E_p (t) \leq \frac{C_{m,\delta} (E_p (0))}{s(t) ^{m(p-\delta)}} \quad \text{ where } \quad 
t = \frac{2s(t)}{H(\frac{1}{(2s(t))^m})} .
\end{equation} 

  \end{Theorem}
{\underline{Two typical examples of application:} 

\begin{itemize} 

\item Assume that there exists $\alpha >0$ and $q>1$ such that
\begin{equation}
\label{hyp-poly}
\forall s \in [0,1], \quad g_0 (s) = \alpha s^q .
\end{equation} 
The estimate \eqref{decay-gene} gives that
\begin{equation}
\label{decay-poly-1}
\forall t\geq 1 , \quad E_p (t) \leq \frac{C_\delta (E_p (0))}{t^{\frac{p}{q} - \delta}} ,
\end{equation}
and this can be improved using the estimate \eqref{decay-gene-m}, which is not ideal but gives almost optimal decay rates here: 
indeed, taking $m$ large enough, \eqref{decay-gene-m} gives that, for every $\delta >0$, there exists $C_\delta (E_p (0))$ such that
\begin{equation}
\label{decay-poly-lim}
\forall t\geq 1 , \quad E_p (t) \leq \frac{C_\delta (E_p (0))}{t^{\frac{p}{q-1} - \delta}} .
\end{equation}

\item Assume that there exists $\alpha >0$ and $q>0$ such that
\begin{equation}
\label{hyp-expo-0}
\forall s \in [0,1], \quad g_0 (s) = \alpha e^{-\alpha / s^q} .
\end{equation} 
Then \eqref{decay-gene} (and \eqref{decay-gene-m}) gives that, for every $\delta >0$, there exists $C_\delta (E_p (0))$ such that
\begin{equation}
\label{decay-expo-0-lim}
\forall t\geq 2 , \quad E_p (t) \leq \frac{C_\delta (E_p (0))}{(\ln t)^{\frac{p}{q} - \delta}} .
\end{equation} 
\end{itemize}

The estimates \eqref{decay-poly-lim} and \eqref{decay-expo-0-lim} are almost optimal in the following sense: 
consider the wave equation with boundary damping, acting at the boundary point:
\begin{equation}
\label{eq-wave-bd}
\begin{cases}
z_{tt} - z_{xx} = 0 \quad  &\text{ for } (t,x) \in \Bbb R_+ \times (0,1), \\
z(t,0)=0 \quad  &\text{ for } t \in \Bbb R_+, \\
z_x(t,1)  +  g(z_t (t,1)) = 0 \quad  &\text{ for } t \in \Bbb R_+, \\
z(0,x)=z_0(x), z_t(0,x)=z_1 (x) \quad  &\text{ for } x \in  (0,1) .
\end{cases}
\end{equation}
Then, using \cite{JV-PM-00} (see also \cite{Chitour-Marx-Mazanti2021}), there exist initial conditions $(z_0,z_1)$
such that 

\begin{itemize}
\item in the case \eqref{decay-poly-1}, the solution $z$ satisfies
$$ E_p (t) \sim \frac{C_\delta (E_p (0))}{t^{\frac{p}{q-1} }}  \quad \text{ as } t \to +\infty , $$
\item in the case \eqref{hyp-expo-0}, the solution $z$ satisfies
$$ E_p (t) \sim  \frac{C_\delta (E_p (0))}{(\ln t)^{\frac{p}{q} }}  \quad \text{ as } t \to +\infty .$$
\end{itemize}
(\cite{JV-PM-00} gives lower bounds of energy for several classes of damping function $g$ in the $L^2$ framework, but the proofs are immediately adaptable to the $L^p$ framework.)
This is why we say that \eqref{decay-poly-lim} and \eqref{decay-expo-0-lim} are almost optimal. The optimal ones would probably be the same ones with $\delta =0$, but we were not able to obtain such estimates. 


The proof of Theorem \ref{thm-stab} is based on the multiplier method with suitably chosen multipliers from \cite{KMJC2022} (done in the linear case), on the method developed in \cite{PM-COCV} (to study the nonlinear stabilization in the $L^2$ framework), and on new integral inequalities. We refer the reader to Lemmas \ref{lem-COCV} and \ref{lem-COCV-m}.

\medskip



\subsubsection{Stability results when $p\in (1,2)$} \hfill

When $p\in (1,2)$, we need to combine previous arguments with convex analysis tools, in particular, the convex conjugate of some suitably chosen convex function $\tilde F$ (that will replace the function $F(y)=\vert y \vert ^p$), and the Fenchel's inequality, see section \ref{sec-Fenchel}. We obtain the following result.

 \begin{Theorem}
 \label{thm-stab2}

 Fix $p\in (1,2)$.  Assume that \eqref{a0} and \eqref{eq-comp-g} hold. Then we have the following:
 
 (i) If $g$ is linear, then there exists $C, \omega >0$ such that, for every initial condition $(z_0,z_1) \in X_p$, we have
 \begin{equation}
\label{decay-exp*}
\forall t\geq 0, \quad E_p (t) \leq C E_p (0) \, e^{-\omega t} .
\end{equation}
 
 (ii) If $g'(0) >0$, the energy $E_p$ decays exponentially to $0$, i.e., for every
 initial condition $(z_0,z_1) \in X_p$,  there exists $C^{(0)} >0$, $\omega ^{(0)}  >0$
 depending on $E_p (0)$ such that, 
 \begin{equation}
\label{decay-exp*2}
\forall t\geq 0, \quad E_p (t) \leq C^{(0)} E_p (0) \, e^{-\omega ^{(0)} t} .
\end{equation}
 
 (iii) If $g'(0)=0$, assume that there exists $\eta >0$ such that the function 
 $$H(s) := \frac{g_0(s)}{s}$$ is increasing on $[0,\eta]$. Then, for every
 initial condition $(z_0,z_1) \in X_p$, there exists $C (E_p(0))$ such that
 \begin{equation}
\label{decay-gene*}
\forall t\geq 1, \quad E_p (t) \leq C (E_p (0)) \, \Bigl( g_0 ^{-1} (\frac{1}{t}) \Bigr) ^{p} .
\end{equation} 
Moreover, we also have the following: then, given $m\geq 1$, we have
\begin{equation}
\label{decay-gene-m*}
\forall t\geq t_m , \quad E_p (t) \leq \frac{C_\delta (E_p (0))}{s(t) ^{pm}} \quad \text{ where } \quad 
t = \frac{2s(t)}{H(\frac{1}{(2s(t))^m})} 
\end{equation} 
and $t_m$ is given in Theorem \ref{thm-stab}.
 \end{Theorem}
 
 {\underline{Two typical examples of application:} 
 
 \begin{itemize}
\item if $g_0$ satisfies \eqref{hyp-poly}, then the estimate \eqref{decay-gene-m*} (with $m$ large enough) gives that, for every $\delta >0$, there exists $C_\delta (E_p (0))$ such that
\eqref{decay-poly-lim} holds;

\item if \eqref{hyp-expo-0} holds, then
\eqref{decay-gene*} gives that there exists $C (E_p (0))$ such that
\begin{equation}
\label{decay-expo-0-lim*}
\forall t\geq 2 , \quad E_p (t) \leq \frac{C (E_p (0))}{(\ln t)^{\frac{p}{q}}} .
\end{equation} 
\end{itemize}
 

\subsection{Organization of the paper} \hfill

  The paper is organized as follows: 
  \begin{itemize}
  \item In Section \ref{sec-litt}, we compare our results with the already existing ones,
  \item In Section \ref{sec-wellposed}, we prove Theorems \ref{wellposedness-w2}-\ref{wellposedness-w1}, starting by rewriting \eqref{eq-wave} using the Riemann invariants.
  \item In section \ref{sect4}, we prove Theorem \ref{thm-stab}.
   \item In section \ref{sect6}, we prove Theorem \ref{thm-stab2}.
\end{itemize}


\section{Comparison with earlier results} 
 \label{sec-litt}
 
 \subsection{The $L^2$ framework} \hfill
 
 The nonlinear problem \eqref{eq-wave} has been already hugely studied in the Hilbertian framework, i.e., with $p=2$. 
 The well-posedness is a classical result and is a consequence of the theory of maximal monotone operators. 
 Decay estimates have been known for a long time, and we refer in particular to Nakao \cite{Nakao1977}, Haraux \cite{Haraux1978},
 Lagnese \cite{Lagnese1983}, Zuazua \cite{Zuazua1990}, Komornik \cite{Komornik1993}
 for estimates when the damping function has a polynomial behaviour near $0$, to Lasiecka and Tataru \cite{Las-Tat1993}, 
 Martinez \cite{PM-COCV}, Alabau-Boussouira \cite{Alabau2005, Alabau2010} under more general conditions on the damping function, and to Vancostenoble and Martinez \cite{JV-PM-00} for lower bound estimates, essentially proving the optimality of the previous upper estimates. These upper estimates are mainly obtained thanks to the multiplier method that allows one to obtain integral estimates on the energy, and then to some integral inequalities of the Gronwall type, or comparison with ODE, for which the asymptotic behaviour of the solutions can be more easily estimated. In particular, Alabau-Boussouira \cite{Alabau2005, Alabau2010} succeeded into proving optimal estimates under very few assumptions, and for a large class of problems, combining cleverly convexity arguments with weighted integral inequalities.


\subsection{The $L^p$ framework} \hfill
  
The non-Hilbertian framework has only been considered recently, even if, from an applied point of view, the $L^\infty$ setting is particularly interesting (and challenging).  
The usual well-posedness Hilbertian techniques (maximal monotone operators for instance) are not easy to handle in this framework. 

The main results that exist and are relevant to our context come primarily from Haraux \cite{Haraux2009},
and were extended in several directions in Chitour, Marx and Prieur \cite{Chitour-Marx-Prieur-2020}, and Kafnemer, Mebkhout, Jean and Chitour \cite{KMJC2022}. 

Haraux \cite{Haraux2009} proves that the semigroup associated with the d'Alembertian operator in one space dimension
 and with Dirichlet boundary conditions is well-posed for every $X_p$, $p\in [2,+\infty]$. Haraux also provides 
some decay rate estimates of an energy equivalent to the natural one for solutions associated to smooth initial data, and when $a \equiv 1$ and $g$ behaves polynomially near $0$ and at infinity: if there exists $k_0, k_1 >0$, and $R\geq r >0$ such that
\begin{equation}
\label{hyp-Haraux2009}
\forall s \in \Bbb R, \quad k_0 \vert s \vert ^r \leq g'(s) \leq k_1 (\vert s \vert ^r  + \vert s \vert ^R),
\end{equation}
and $(z_0,z_1) \in Y_{\infty}$, then 
\begin{equation}
\label{estim-Haraux2009}
\forall t \geq 0, \quad \Vert z \Vert ^p _{W^{1,p}_0(0,1)} + \Vert z_t \Vert ^p _{L^p (0,1)} \leq \frac{C(z_0,z_1)}{(1+t)^{\frac{2(p+1)}{3r}}} .
\end{equation}

In \cite{Amadori2019}, Amadori, Aqel and Dal Santo provide an asymptotic analysis when $a$ is positive and $g$ has a linear behaviour:
$$ \begin{cases} 0 < k_1 \leq a(x) \leq k_2, \\ \forall x \in (0,1) , \end{cases} \quad \text{ and } \quad \quad 
\begin{cases} 0 < g_1 \leq g'(s) \leq g_2 , \\ \forall s \in \Bbb R  . \end{cases} $$
The techniques used in this paper are based on the theory of scalar conservation laws, which makes the analysis in those $L^p$ spaces really natural.

In \cite{Chitour-Marx-Prieur-2020}, Chitour, Marx and Prieur study \eqref{eq-wave} when $p \in [2,+\infty]$, and provide well-posedness results in the associated $L^p$ framework (relying on the d'Alembert formula); they also provide semi-global exponential stability results when the damping function satisfies mainly $g'(0) >0$: given $R >0$, there exists $K(R), \beta (R) >0$ such that
$$ \Vert (z_0,z_1) \Vert _{X_p} \leq R \quad \implies \quad \Vert (z(t),z_t (t) ) \Vert _{X_p} \leq K(R) e^{-\beta (R) t} ,
\quad \forall t \geq 0 .$$

Finally, in \cite{KMJC2022}, Kafnemer, Mebkhout, Jean and Chitour study \eqref{eq-wave} when $g(s) \equiv k_0 s$ with some $k_0 >0$, extend the well-posedness results of \cite{Chitour-Marx-Prieur-2020} to the case $p>1$, obtain exponential decay
for all $1<p<+\infty$ using a generalized multiplier method and exponential stability when $p=\infty$ 
in some cases of a constant in space global damping.

Here we use as a basis \cite{K2022, KMJC2022} alongside with some techniques from \cite{PM-COCV} and \cite{Chitour-Marx-Prieur-2020} to prove our stability results. However, combining the existing arguments was not sufficient, and we needed new integral inequalities necessary to conclude when $p>2$, see Lemmas \ref{lem-COCV} and \ref{lem-COCV-m}, and to extend or precise some estimates when $p \in (1,2)$.


\subsection{Open questions} \hfill

There are some natural open questions associated to our work:

\begin{itemize}
\item Optimality of the decay estimates of Theorem \ref{thm-stab}: when $p>2$, we believe that the optimal decay estimates would be the ones of Theorem \ref{thm-stab}, with $\delta =0$. It would be interesting to know if such decay estimates hold true.

\item Optimality of the decay estimates of Theorem \ref{thm-stab2}: in case (ii), we do not know if the decay rate can be uniform with respect to the norm of the initial conditions. 
It would also be interesting to know if the estimate \eqref{decay-poly-lim} holds with $\delta =0$.

\item Our proofs are not valid in the case $p\in\{1,\infty\}$. However, there are some results in some particular cases in \cite{K2022}, when the function $a(x)$ is constant on $(0,1)$, and with some constraints on the value of that constant. It would be interesting to go further in that direction, and to obtain estimates in $X_p$, $p\in\{1,\infty\}$, with less restrictive assumptions on the feedback term.
\end{itemize}


\section{Well-posedness}
\label{sec-wellposed}

In this section, we study the well-posedness of Problem \eqref{eq-wave}. Before, we introduce a crucial tool, namely, the Riemann invariants.


\subsection{Riemann invariants} \hfill

We define the Riemann invariants for all $(t,x) \in \mathbb{R_+}\times (0,1)$ by
\begin{align}\label{Riemann}
\rho(t,x)= {z_x(t,x)+z_t(t,x)}, \\
\xi(t,x)= {z_x(t,x)-z_t(t,x)}.
\end{align}
Along strong solutions of \eqref{eq-wave}, we deduce that
\begin{equation}
\label{eq-invariants}
\begin{cases}
\begin{cases}
\rho _t - \rho _x = -a(x) g(\frac{\rho - \xi}{2}) \quad  &\text{ for } (t,x) \in \Bbb R_+ \times (0,1), \\
\xi _t + \xi _x = a(x) g(\frac{\rho - \xi}{2}) \quad  &\text{ for } (t,x) \in \Bbb R_+ \times (0,1), 
\end{cases}
\\
\smallskip
\\
\begin{cases}
\rho (t,0) - \xi (t,0) =0 \quad  &\text{ for } t \in \Bbb R_+, \\
\rho (t,1) - \xi (t,1) =0 \quad  &\text{ for } t \in \Bbb R_+, 
\end{cases}
\\
\smallskip
\\
\begin{cases}
\rho(0,x)=z_0 ' (x) + z_1 (x) =: \rho _0 (x)\quad  &\text{ for } x \in  (0,1), \\
\xi(0,x)= z_0 ' (x) - z_1 (x) =: \xi _0 (x) \quad  &\text{ for } x \in  (0,1) .
\end{cases}
\end{cases}
\end{equation}
Note that $\left(\rho_0, \xi_0 \right) \in W^{1,p}(0,1)\times W^{1,p}(0,1)$.


Before proving Theorems \ref{wellposedness-w2}-\ref{wellposedness-w1}, we need to study the monotonicity of the associated energy.


\subsection{The energy $E_p$ of strong solutions is nonincreasing} \hfill

	For $r\geq 0$, we introduce the following notation
	\begin{align}
\lfloor x\rceil^{r}:=\textrm{sgn}(x)|x|^r ,\ \   \forall x \in \mathbb{R},
	\end{align}
	where $\textrm{sgn}(x)=\frac{x}{\vert x\vert}$ for nonzero $x\in\mathbb{R}$ and {$\textrm{sgn}(0)=[-1,1]$}.
	We have the following obvious formulas, which will be repeatedly used later on:
	\begin{align}
	\frac{d}{dx}(\lfloor x\rceil^{r})=r|x|^{r-1}, \ \ \forall r\geq 1,\ x\in\mathbb{R},\\
	\frac{d}{dx}(|x|^r)=r\lfloor x\rceil^{r-1}, \ \ \forall r> 1,\ x\in\mathbb{R}.
	\end{align}

First, we recall that the $p$th-energy $E_p$ of a solution $z$ is defined in \eqref{energyp}. The Riemann invariants allow us to have the following expression of $E_p$:
\begin{align}
E_p(t)= \frac{1}{p} \int_0^1 (|{\rho}|^p+|{\xi}|^p) dx.
\end{align}

Here is the natural extension of Proposition 2.1 of \cite{KMJC2022}.

\begin{Proposition}\label{prop1}
	Let $p \in [1, \infty)$ and suppose that a strong solution $z$ of \eqref{eq-wave} exists and is defined on
	a non trivial interval $I \subset \mathbb{R_+}$ containing $0$, for some initial conditions $(z_0,z_1)\in
	Y_p$. Let ${\mathcal{F}} $ be a $C^1$ convex function. For $t\in I$, define
\begin{equation}\label{phidef}
	\Phi (t) := \int_0^1 [{\mathcal{F}}(\rho)+{\mathcal{F}}(\xi)] \, dx,
\end{equation}
	where $\rho$, $\xi$ are the Riemann invariants, defined in \eqref{Riemann}. Then $\Phi$ is well defined
	for $t\in I$ and satisfies
	\begin{equation}\label{phidec}
	 \Phi '(t) =-\int_0^1 a(x) g\left(\frac{\rho-\xi}{2}\right)({\mathcal{F}}'(\rho)- {\mathcal{F}}'(\xi)) dx.
	\end{equation}
As a consequence, the function $\Phi$ is nonincreasing on $I$.
\end{Proposition}

As an immediate consequence of Proposition \ref{prop1}, we will have the following

	\begin{Corollary}\label{cor:decrease}
Let $p \in (1, \infty)$. Assume that $z$ is a strong solution of \eqref{eq-wave}. Denote $\rho$ and $\xi$ the Riemann invariants, defined in \eqref{Riemann}. Then
		\begin{align}\label{E_p'}
		E_p'(t) = -\int_0^1 a(x)g\left(\frac{\rho - \xi}{2}\right)\left( \lfloor \rho\rceil^{p-1}-\lfloor \xi\rceil^{p-1} \right) dx .
		\end{align}
As a consequence, the energy $E_p$ is nonincreasing on its existence interval.  
	\end{Corollary}

{\it Proof of Corollary \ref{cor:decrease} assuming Proposition \ref{prop1}.}
If $p>1$, it is sufficient to apply Proposition \ref{prop1} with $\mathcal F (x) = \vert x \vert ^p$, and
it is clear that $\mathcal F$ is $C^1$ on $\Bbb R$ and convex. \qed

{\it Proof of Proposition \ref{prop1}.} The proof of the formula of $\Phi'(t)$ is similar to the proof of \cite[Proposition 2.1]{KMJC2022}, but for the sake of completeness, we provide it. For $t\geq 0$, it holds
$$ \Phi '(t) =  \int_0^1 \mathcal{F}'(\rho) \rho _t +\mathcal{F}'(\xi) \xi _t \, dx .$$
Using \eqref{eq-invariants}, we have
$\rho _t = \rho _x  -a(x) g(\frac{\rho - \xi}{2})$, and 
$\xi _t =- \xi _x + a(x) g(\frac{\rho - \xi}{2})$. Therefore
\begin{multline*}
\Phi '(t) =  \int_0^1 \mathcal{F}'(\rho) \rho _x - \mathcal{F}'(\xi)\xi _x \, dx - \int _0 ^1 a(x) g(\frac{\rho - \xi}{2}) \Bigl( \mathcal{F}'(\rho) - \mathcal{F}'(\xi) \Bigr) \, dx
\\
= \Bigl[ \mathcal{F}(\rho) - \mathcal{F}(\xi) \Bigr] _0 ^1 - \int _0 ^1 a(x) g(\frac{\rho - \xi}{2}) \Bigl( \mathcal{F}'(\rho) - \mathcal{F}'(\xi) \Bigr) \, dx .
\end{multline*}
Thanks to the Dirichlet boundary conditions, we have $z_t=0$ at the boundary, and then
$$ \mathcal{F}(\rho) _{\vert x = 1} - \mathcal{F}(\xi) _{\vert x = 1} = \mathcal{F}(z_x + z_t ) _{\vert x = 1}
- \mathcal{F}(z_x - z_t ) _{\vert x = 1} = \mathcal{F}(z_x) _{\vert x = 1} - \mathcal{F}(z_x) _{\vert x = 1} = 0 ,$$
and the same thing at $x=0$. This implies that
$$
\Phi '(t) =  - \int _0 ^1 a(x) g(\frac{\rho - \xi}{2}) \Bigl( \mathcal{F}'(\rho) - \mathcal{F}'(\xi) \Bigr) \, dx .
$$
To conclude, since ${\mathcal{F}}'$ is nondecreasing, the sign of ${\mathcal{F}}'(\rho)- {\mathcal{F}}'(\xi)$ is the same than the sign of $\rho - \xi$; and the sign of $g\left(\frac{\rho-\xi}{2}\right)$ is also the same than the sign of $\rho - \xi$. This implies that $g\left(\frac{\rho-\xi}{2}\right)({\mathcal{F}}'(\rho)- {\mathcal{F}}'(\xi)) $ is nonnegative, and then that 
$\Phi '(t) \leq 0$. \qed


\begin{Remark}
{\rm 
	The previous proposition was first introduced in \cite{Haraux1978}  and reused in \cite{Chitour-Marx-Prieur-2020} to prove that the energy functional is nonincreasing. Then it was improved in \cite{KMJC2022} by omitting the hypothesis that function $\mathcal{F}$ should be even on top of being convex. 
}
\end{Remark}


\subsection{Proof of Theorem \ref{wellposedness-w2}} \hfill

Fix $2\leq p< +\infty$ and let $(z_0,z_1)\in X_p$. 
Since $Y_p$ is dense in $X_p$ for all $2\leq p<\infty$ and $Z_0=(z_0,z_1)\in X_p$, there exists a sequence $( Z_{0} ^{(n)} =(z_{0} ^{(n)}, z_{1} ^{(n)})_n$
such that  $Z_{0} ^{(n)} \in Y_p$ and $Z_{0} ^{(n)} \to Z_0$ in $X_p$.

Now, since $Y_p\subset X_{\infty}$, we have $Z_{0} ^{(n)} \in X_\infty$, and  \cite[Theorem 1]{Chitour-Marx-Prieur-2020} gives us 
that there exists a unique solution
$$ z^{(n)} \in L^\infty(\mathbb{R_+};W^{1,\infty}_0(0,1))\cap W^{1,\infty}(\mathbb{R_+};L^\infty(0,1)) $$ 
of \eqref{eq-wave} associated to the initial condition $(z_{0} ^{(n)}, z_{1} ^{(n)})$, 
and we have, for all $t\geq 0$,
\begin{align}
\Vert (z^{(n)},z^{(n)} _t)\Vert _{X_\infty} \leq 2 \max\left(\Vert {z_{0} ^{(n)}}' \Vert _{L^\infty(0,1)}, 
\Vert z_{1} ^{(n)} \Vert _{L^\infty(0,1)}  \right).
\end{align}
Denote $Z^{(n)} := (z^{(n)}, z^{(n)} _t )$. We prove now that, for every $t_0\geq 0$, the sequence $( Z^{(n)}(t_0,\cdot)$ is a Cauchy sequence in $X_p$. This will follow from the following observation:
\begin{multline*} 
\Vert Z^{(n)}(t_0) - Z^{(m)}(t_0) \Vert _{X_p} ^p
= \Bigl\Vert \Bigl( z^{(n)}(t_0) - z^{(m)}(t_0), z^{(n)}_t (t_0) - z^{(m)} _t (t_0) \Bigr) \Bigr\Vert _{X_p} ^p
\\
= \frac{1}{p} \int _0 ^1 \Bigl\vert z_x ^{(n)}(t_0) - z_x ^{(m)}(t_0) \Bigr \vert ^p 
+ \Bigl\vert z^{(n)}_t (t_0) - z^{(m)} _t (t_0) \Bigr\vert ^p \, dx .
\end{multline*}
Define for all $(t,x)\in \mathbb{R_+}\times (0,1)$ the quantity $e^{(n,m)}$ as 
\begin{align}
e^{(n,m)}=z^{(n)}-z^{(m)} .
\end{align}
Then
$$ 
\Vert Z^{(n)}(t_0) - Z^{(m)}(t_0) \Vert _{X_p} ^p
= \frac{1}{p} \int _0 ^1 \Bigl\vert e^{(n,m)} _x (t_0) \Bigr\vert ^p
+ \Bigl\vert  e^{(n,m)}_t (t_0) \Bigr\vert ^p \, dx = \mathcal E _p (e^{(n,m)}) (t_0) ,
$$
where $\mathcal E _p$ is defined in \eqref{energyp-init}. But $\mathcal E _p$ is equivalent to the energy defined in \eqref{energyp}: there exists $C^* > C_* '$ such that
$$ C_* E_p \leq \mathcal E_p \leq C^* E_p ,$$
and therefore
we obviously have
\begin{equation}
\label{rel-evid}
\Vert Z^{(n)}(t_0) - Z^{(m)}(t_0) \Vert _{X_p} ^p \leq C^* E_p ( e^{(n,m)}) (t_0) .
\end{equation}
Then the property that the sequence $( Z^{(n)}(t_0,\cdot))_n$ is a Cauchy sequence in $X_p$ will derive from a suitable estimate on 
$E_p ( e^{(n,m)}) (t_0)$, and this will be a consequence of the problem satisfied by the function $e^{(n,m)}$, which is the following one: 
\begin{equation}\label{probe}
\begin{cases}
e^{(n,m)}_{tt} -e^{(n,m)}_{xx} + a(x)\left(g(z^{(n)} _t )-g(z^{(m)} _t)\right)=0, \quad (t,x) \in  \mathbb{R_+} \times (0,1), \\
e^{(n,m)}(t,0)= e^{(n,m)}(t,1)=0 , \quad    t\geq 0 ,\\
e^{(n,m)}(0,\cdot)= z_0^{(n)}-z_0^{(m)}\ ,\ e^{(n,m)}_t(0,\cdot)=z_1^{(n)}-z_1^{(m)}.
\end{cases}
\end{equation}
Since 
$$ g(y)-g(x) = \int _x ^y g'(s) \, ds = \Bigl( \int _0 ^1 g'( \tau y + (1-\tau) x) \, d\tau \Bigr) (y-x) ,$$
we have
\begin{multline*}
 g(z^{(n)} _t )-g(z^{(m)} _t) = \Bigl( \int _0 ^1 g'( \tau z^{(n)} _t + (1-\tau) z^{(m)} _t) \, d\tau \Bigr) (z^{(n)} _t-z^{(m)} _t) 
 \\
 = \Bigl( \int _0 ^1 g'( \tau z^{(n)} _t + (1-\tau) z^{(m)} _t) \, d\tau \Bigr) e^{(n,m)} _t .
 \end{multline*}
Denoting
 $$ A^{(n,m)} (t,x) := a(x) \Bigl( \int _0 ^1 g'( \tau z^{(n)} _t + (1-\tau) z^{(m)} _t) \, d\tau \Bigr) ,$$
we see that  $e^{(n,m)}$ satisfies
\begin{equation}\label{probe2}
\begin{cases}
e^{(n,m)}_{tt} -e^{(n,m)}_{xx} + A^{(n,m)} (t,x) e^{(n,m)} _t =0, \quad (t,x) \in  \mathbb{R_+} \times (0,1), \\
e^{(n,m)}(t,0)= e^{(n,m)}(t,1)=0 , \quad    t\geq 0 ,\\
e^{(n,m)}(0,\cdot)= z_0^{(n)}-z_0^{(m)}\ ,\ e^{(n,m)}_t(0,\cdot)=z_1^{(n)}-z_1^{(m)}.
\end{cases}
\end{equation}
Proceeding as in Proposition \ref{prop1}, one has in the same way
$$ \frac{d}{dt} E_p (e^{(n,m)}) (t) = 
-\frac{1}{2} \int_0^1 A^{(n,m)} (t,x) \Bigl( \rho ^{(n,m)} - \xi ^{(n,m)} \Bigr) \Bigl( \mathcal F'(\rho ^{(n,m)}) -\mathcal F'(\xi ^{(n,m)}) \Bigr) dx , $$
where $\mathcal F(x)=\vert x \vert ^p$, and $\rho ^{(n,m)}$ and $\xi ^{(n,m)}$ are the Riemann invariants associated to $e^{(n,m)}$: 
$\rho ^{(n,m)} = e^{(n,m)} _x + e^{(n,m)} _t$, and $\xi ^{(n,m)} = e^{(n,m)} _x - e^{(n,m)} _t$.
Since $g$ is nonincreasing, $g'$ is nonnegative, and thus $A^{(n,m)} \geq 0$. This implies that $E_p (e^{(n,m)})$ is nonincreasing on $[0,+\infty)$.
It follows that 
$$
\forall t\geq 0, \quad E_p(e^{(n,m)})(t) \leq E_p(e^{(n,m)})(0),
$$
and using \eqref{rel-evid}, we have
$$
 \forall t_0\geq 0, \quad \Vert Z^{(n)}(t_0) - Z^{(m)}(t_0) \Vert _{X_p} ^p
\leq C^* E_p(e^{(n,m)})(0)
\leq \frac{C^*}{C_*} \Vert Z^{(n)}_0  - Z^{(m)} _0 \Vert _{X_p} ^p .
$$
Since $(Z^{(n)}_0)_n$ is a Cauchy sequence in $X_p$, we obtain that for all $t_0 \geq 0$, the sequence $(Z^{(n)} (t_0))_n$ is also a Cauchy sequence in $X_p$. It follows then that $( Z^n(t_0,\cdot) )_n $ converges to a limit in $X_p$ that we denote by $Z(t_0,\cdot)=(z(t_0,\cdot),z_t(t_0,\cdot))$. In particular, for $t_0=0$, we have that $Z(0,\cdot)=Z_0(\cdot)$. Note also that the convergence is uniform with respect to $t_0\geq 0$. We define the function $(t,x)\mapsto Z(t,x)$, where $Z(t_0,\cdot)$ is the limit of $Z^n(t_0,\cdot)$ in $X_p$ for every $t_0\geq 0$. 

It remains to prove now that the limit $Z$ is a weak solution of \eqref{eq-wave}. Fix $T>0$ and denote $\psi$ a test function that belongs to $C^1([0,T]\times[0,1])$
 also verifying $\psi(T,\cdot)=\psi(0,\cdot)\equiv 0$ and {$\psi(\cdot,0)=\psi(\cdot,1)\equiv 0$}. Define
\begin{align}
B_T^{(n)}(\psi)=	\int_0^T \int_0^1 (z_{tt}^{(n)}-z_{xx}^{(n)})\, \psi \, dx dt.
\end{align}
We have that
$$
B_T^{(n)} (\psi) - B_T^{(m)} (\psi)= \int_0^T \int_0^1 (z_{tt}^{(n)} -z_{xx}^{(n)})\psi \, dxdt -\int_0^T \int_0^1 (z_{tt}^{(m)} -z_{xx}^{(m)} )\psi \, dxdt .
$$
By integrating by parts, it follows that
$$
B_T^{(n)} (\psi)-B_T^{(m)} (\psi)= -\int_0^T \int_0^1 (z_{t}^{(n)} -z_{t}^{(m)})\psi_t \, dxdt +\int_0^T \int_0^1 (z_{x}^{(m)}-z_{x}^{(n)})\psi_x \, dxdt.
$$
Using H\"older's inequality,
\begin{align*}
\Bigl\vert B_T^{(n)} (\psi)-B_T^{(m)} (\psi) \Bigr\vert &\leq  \left(\int_0^T \int_0^1|z_t^{(n)}-z_t^{(m)}|^p dxdt\right)^{\frac{1}{p}}   \left(\int_0^T \int_0^1|\psi_t|^q dxdt\right)^{\frac{1}{q}}
\notag 
\\&+\left(\int_0^T \int_0^1 |z_x^{(n)}-z_x^{(m)}|^p dxdt\right)^{\frac{1}{p}}   \left(\int_0^T \int_0^1 |\psi_x|^q dxdt\right)^{\frac{1}{q}},
\end{align*}
which means that 
$$
\Bigl\vert B_T^{(n)} (\psi)-B_T^{(m)} (\psi) \Bigr\vert \leq T^{\frac{1}{p}} E_p(e^{(n,m)})^{\frac{1}{p}}(0) \left(||\psi_t||_{L^q((0,T)\times (0,1))}+||\psi_x||_{L^q((0,T)\times (0,1))}\right).
$$
By a density argument, we obtain that for all $\psi$ in the space
\begin{multline*}
\mathcal X _q ^T  = \{ \psi: [0,T] \times [0,1] \mapsto \Bbb R_+, 
(\psi,\psi_t)\in W^{1,q}((0,T)\times(0,1)) \times L^q((0,T)\times(0,1)), 
\\
\psi(T,\cdot)=\psi(0,\cdot)\equiv 0 \text{ and } \psi(\cdot,0)=\psi(\cdot,1)\equiv 0 \} ,
\end{multline*}
where $q$ is the conjugate exponent of $p$ and is equal to $\frac{p}{p-1}$,
we have that
$$
\Bigl\vert B_T^{(n)} (\psi)-B_T^{(m)} (\psi) \Bigr\vert \leq T^{\frac{1}{p}} E_p(e^{(n,m)})^{\frac{1}{p}}(0) \left(||\psi_t||_{L^q((0,T)\times (0,1))}+||\psi_x||_{L^q((0,T)\times (0,1))}\right) .
$$
Then
$$
\Bigl\vert B_T^{(n)} (\psi)-B_T^{(m)} (\psi) \Bigr\vert \leq T^{\frac{1}{p}} E_p(e^{(n,m)})^{\frac{1}{p}}(0) ||\psi||_{\mathcal{X}_q^T},
$$
which gives that $( B_T^{(n)} )_n $ is a Cauchy sequence in $(\mathcal{X}_q^T)'$, the dual of $\mathcal{X}_q^T$, since $( Z_0^{(n)})_n$ is in $X_p$. We conclude then that $(B_T^{(n)} )_n $ converges in $(\mathcal{X}_q^T)'$, i.e.,  $(  z^{(n)}_{tt} -z^{(n)} _{xx})_n$ converges to ${z}_{tt} -{{z}}_{xx}$ as a linear functional on $\mathcal{X}_q^T$.
We also know that for all $n\in \mathbb{N}$, 
$$
z^{(n)} _{tt} -z_{xx}^{(n)} = -a(x) g(z^{(n)} _t), \text{ on } \mathbb{R_+}\times (0,1),
$$
which means that, the sequence of linear functionals $( -a(x) g(z^{(n)} _t) )_n$  defined on $\mathcal{X}_q^T$ also converges ${z}_{tt}(t,\cdot)-{{z}}_{xx}(t,\cdot)$ in $(\mathcal{X}_q^T)'$. 
 We use now the existence of weak solutions in $L^2$ framework (see \cite{PM-JV2000}, \cite{KMC2021}). For  $(z_0,z_1)\in X_p$ with $p\geq 2$, we have the existence of a unique weak solution ${z}\in W^{1,\infty}(\mathbb{R_+},L^2(0,1))\cap L^\infty(\mathbb{R_+},H^1_0(0,1))$. The solution $z$ satisfies that for almost every $t\in \mathbb{R_+}$
\begin{align}
{z}_{tt}(t,\cdot)-{z}_{xx}(t,\cdot)=-a(x) g({z}_t(t,\cdot) \ \ \ \text{ in } H^{-1}(0,1).
\end{align}
In particular, for every $T>0$, ${z}_{tt}(t,\cdot)-{z}_{xx}(t,\cdot)=-a(x)g({z}_t)$ belongs to $(\mathcal{X}_2^T)'$. Since $q\leq 2$, one has that  $(\mathcal{X}_q^T)'\subset (\mathcal{X}_2^T)'$, yielding in particular that  $z$ is a weak solution of System \eqref{eq-wave} in $X_p$. \qed



\subsection{Proof of Theorem \ref{wellposedness-s1}} \hfill

Take $Z_0=(z_0,z_1)\in Y_p$ for all $1\leq p<\infty$. Using Sobolev embeddings classical results, we have that $Z_0\in X_\infty$.
Then we are in position to use once again \cite[Theorem 1]{Chitour-Marx-Prieur-2020}, and we have the existence of a unique solution $z$ of \eqref{eq-wave} such that $$(z,z_t)\in L^\infty(\mathbb{R_+}; W^{1,\infty}_0(0,1) )\times W^{1,\infty}(\mathbb{R_+};L^{\infty}(0,1)),$$ and,  for all $t\geq 0$,
\begin{align}
\Vert (z,z_t) \Vert _{X_\infty} \leq 2 \max\left(\Vert z_0' \Vert _{L^\infty(0,1)}, \Vert z_1 \Vert _{L^\infty(0,1)}  \right).
\end{align}
To have more information on $z_t$, we note that $w:=z_t$ satisfies
\begin{equation}\label{weq2}
\left\lbrace
\begin{array}{ll}
w_{tt} - w_{xx}  = - a(x)\,  g'(w) \, w_t  &\text{in }   \mathbb{R_+} \times (0,1), \\
w(t,0)=w(t,1)=0   &  \forall t \in \mathbb{R_+} ,\\
w(0,\cdot)=  z_1,\  w_t(0,\cdot)=z_0''-g(z_1).
\end{array}
\right.
\end{equation}
%
%
We define for all $(t,x)\in \mathbb{R_+}\times (0,1)$ the associated Riemann invariants for \eqref{weq2}:
\begin{align}
u=w_x+w_t,\\
v=w_x-w_t.
\end{align}
Along strong solutions of \eqref{weq2}, we have 
\begin{equation}\label{weqr}
\left\lbrace
\begin{array}{ll}
u_{t} - u_{x}  = - \frac{a(x)}{2}\,  g'(w) \, (u-v) &\text{in }   \mathbb{R_+} \times (0,1), \\
v_{t} + v_{x}  = \frac{a(x)}{2}\,  g'(w) \, (u-v) &\text{in }   \mathbb{R_+} \times (0,1), \\
u(t,0)-v(t,0)=u(t,1)-v(t,1)=0   &  \forall t \in \mathbb{R_+}.
\end{array}
\right.
\end{equation}
Introducing
$$ \tilde a (t,x) := a(x) g'(w) ,$$
Problem \eqref{weqr} can be written
\begin{equation}\label{weqr2}
\left\lbrace
\begin{array}{ll}
u_{t} - u_{x}  = - \tilde a(t,x) \frac{u-v}{2} &\text{in }   \mathbb{R_+} \times (0,1), \\
v_{t} + v_{x}  = \tilde a(t,x) \frac{u-v}{2}  &\text{in }   \mathbb{R_+} \times (0,1), \\
u(t,0)-v(t,0)=u(t,1)-v(t,1)=0   &  \forall t \in \mathbb{R_+}.
\end{array}
\right .
\end{equation}
Then, we consider the $p$-th energy associated with $w$:
\begin{align}
E_p(w)(t)=\frac{1}{p}\int_0^1(|u|^p+|v|^p)\, dx ,
\end{align}
and proceeding as in Corollary \ref{cor:decrease}, we see that $E_p(w)$ is nonincreasing, the fundamental remark to obtain that property being that the function $\tilde a$ is nonnegative, since $g' \geq 0$. Hence
$$
E_p(w)(t) \leq E_p(w) (0),\quad \textrm{ for a.e. } t\geq 0.
$$
Then, using the fact that $E_p(w)$ is an equivalent energy to $\mathcal E_p (w)$, 
it follows that 
$$
\int_0^1|w_x|^p\, dx \leq p C_ p E_p(w) (0) .
$$
Therefore
$$
\Vert z_t \Vert _{W^{1,p}(0,1)} \leq  (p C^* E_p(w) (0))^{\frac{1}{p}} ,
$$
and this implies that $(z,z_t)\in Y_p$ for all $t\geq 0$, and yields the required regularity for a strong solution. \qed

 \begin{Remark}
{\rm 	For strong solutions in the case $p\geq 2$, we can easily use the results that have been proved in \cite{Haraux2009} for $p\geq 2$ to prove the well-posedness. Indeed, let $(z_0,z_1)\in Y_p$, with $ p\geq 2$ this implies that $(z_0,z_1)\in \left(H^2(0,1)\cap H^1_0(0,1)\right)\times H^1_0(0,1)$. We have then the existence of a unique strong solution $$z\in C(\mathbb{R_+},H^1_0(0,1))\cap C^1(\mathbb{R_+},L^2(0,1)),$$ such that $$(z(t,\cdot),z_t(t,\cdot) )\in \left(H^2(0,1)\cap H^1_0(0,1)\right)\times H^1_0(0,1), \ \ \forall t\in \mathbb{R_+},$$ which means that $z_t(t,\cdot)\in L^\infty(0,1)$ we can then use \cite[Corollary 2.3, item (ii) ]{Haraux2009} that implies in our context that if $(z_0,z_1)\in Y_p$ then the solution $z\in L^\infty(\mathbb{R_+},W^{2,p}(0,1)\cap W^{1,p}_0(0,1))$ and $z_t\in L^\infty(\mathbb{R_+},W^{1,p}_0(0,1))$ which guarantees the well-posedness in $Y_p$.
} \end{Remark}


\subsection{Proof of Theorem \ref{wellposedness-w1}} \hfill

Now assume Hypotheses $\mathbf{ (H_1)}$ and $\mathbf{ (H_3)}$.
In that case, the argument follows closely \cite[Corollary 2.3]{Haraux1978} or \cite[Theorem 3.3]{KMJC2022} and
we refer the reader to these references. \qed


\section{Proof of Theorem  \ref{thm-stab}}
\label{sect4}

\subsection{Asymptotic stability: the key integral inequality on the energy} \hfill

For the rest of the paper, we will manipulate (essentially) only strong solutions of \eqref{eq-wave}, and a standard density will complete the arguments in order to obtain conclusions relative to weak solutions of \eqref{eq-wave}.

In the following we assume that $\phi: \Bbb R_+ \to \Bbb R_+$ is increasing and concave. The goal of this section is to prove 
the next lemma, which will be useful in determining the asymptotic behaviour of the energy. The function $\phi$ will be chosen later.

\begin{Lemma}
\label{lem-int-estim}
There exists $C_u >0$ (independent of the initial conditions) such that, for every pair of times $0\leq S\leq T$ and, it holds
\begin{multline}
\label{eq-int-estim}
\int _S ^T E_p (t) ^2 \phi ' (t) \, dt
\leq C_u E_p (S) ^2  \phi ' (S) 
\\ + C_u \int _S ^T E_p (t) \phi ' (t) \int _0 ^1 a (x) \, \frac{\rho - \xi}{2} \Bigl( f(\rho) - f(\xi)\Bigr) \, dx \, dt .
\end{multline}
\end{Lemma}
The proof of Lemma \ref{lem-int-estim} is based on the estimates \eqref{mult-1st-mult-cons}, \eqref{mult-2nd-mult-cons} and \eqref{mult-3rd-mult-cons}, that are consequences of the
identities \eqref{mult-1st-mult}, \eqref{mult-2nd-mult} and \eqref{mult-3rd-mult}.
These identities will be obtained using the multiplier method, and we will use some cut-off functions in order to localize some information in the damping region.

In order to make the paper easier to read, first we give the main identities (obtained by the multiplier method) and their consequences (see Lemmas
\ref{lem-1st-estim}-\ref{lem-3rd-estim}), then we prove that Lemmas
\ref{lem-1st-estim}-\ref{lem-3rd-estim} imply Lemma \ref{lem-int-estim}; then we show how Lemma \ref{lem-int-estim} implies Theorem \ref{thm-stab}; and at the end we prove Lemmas \ref{lem-1st-estim}-\ref{lem-3rd-estim} with the multiplier method and suitable estimates.


\subsection{Multiplier method: useful estimates} \hfill

\subsubsection{Some notations} \hfill

Consider
\begin{equation}
\label{def-F-f}
F(s) := \frac{\vert s \vert ^p}{p} , \quad f(s) := F'(s) = \lfloor s\rceil^{p-1}.
\end{equation}

Choose $a_3 < a_2 < a_1 < a_0 < b_0 < b_1 < b_2 < b_3$ such that $(a_3, b_3) \subset \subset \omega$, and 
$0 < \varepsilon _0 < \varepsilon _1 < \varepsilon _2 < \varepsilon _3$ and 
introduce the following cut-off functions:

\begin{itemize}
\item function $\psi _1: \Bbb R \to [0,1]$, smooth, equal to $1$ on $(-\infty, a_1] \cup [b_1,+\infty)$
and equal to $0$ on $[a_0, b_0]$;
\item function $ \psi _2: \Bbb R \to [0,1]$, smooth, equal to $0$ on $(-\infty, a_2] \cup [b_2,+\infty)$
and equal to $1$ on $[a_1, b_1]$;
\item function $ \psi _3: \Bbb R \to [0,1]$, smooth, equal to $0$ on $(-\infty, a_3] \cup [b_3,+\infty)$ and equal to $1$ on $[a_2, b_2]$.
\end{itemize}
Of course, we have that there exists $C >0$ such that $ \psi _3 \leq C a$.


\subsubsection{First multiplier, first identity, first estimate} \hfill

Multiplying the first equation of \eqref{eq-invariants} by $E_p \phi ' x \psi_1 (x) f(\rho)$ and the second equation of \eqref{eq-invariants} by $E_p \phi ' \,  x \psi _1 (x) f(\xi)$, we obtain the following first identity: along strong solutions of \eqref{eq-wave}, it holds
\begin{multline}
\label{mult-1st-mult}
\int _S ^T E_p (t) ^2 \phi ' (t) \, dt
= \Bigl[  E_p (t) \phi ' (t) \int _0 ^1 x \psi _1 (x) \Bigl( F(\xi) - F(\rho) \Bigr) \, dx \Bigr] _S ^T 
\\
+ \int _S ^T E_p (t) \phi ' (t) \int _0 ^1 \Bigl( 1 - (x\psi _1 (x))' \Bigr) \Bigl( F(\rho) + F(\xi) \Bigr) \, dx \, dt
\\
+ \int _S ^T (E_p (t) \phi '(t) )' \int _0 ^1 x \psi _1 (x) \Bigl( F(\rho) - F(\xi) \Bigr) \, dx \, dt
\\
- \int _S ^T E_p (t) \phi ' (t) \int _0 ^1 x a(x) \psi _1 (x) g(\frac{\rho - \xi}{2}) \Bigl( f(\rho) + f(\xi) \Bigr) \, dx \, dt .
\end{multline}

This first identity \eqref{mult-1st-mult} allows us to obtain the following first estimate:

\begin{Lemma}
\label{lem-1st-estim}
There exists $C_1 >0$ independent of the initial conditions such that
\begin{multline}
\label{mult-1st-mult-cons}
\int _S ^T E_p (t) ^2 \phi ' (t) \, dt
\leq \frac{7}{2} E_p (S) ^2  \phi ' (S)
\\
+ C_1 \int _S ^T E_p(t)  \phi ' (t) \int _{0} ^1  \psi _2 (x) \Bigl( \rho f(\rho) + \xi f (\xi) \Bigr) \, dx \, dt
\\
- \int _S ^T E_p (t) \phi ' (t) \int _0 ^1 x a(x) \psi _1 (x) g(\frac{\rho - \xi}{2}) \Bigl( f(\rho) + f(\xi) \Bigr) \, dx \, dt.
\end{multline}
\end{Lemma}


\subsubsection{Second multiplier, second identity, second estimate} \hfill
\label{sec-sec-mult}

Multiplying the first equation of \eqref{eq-invariants} by $E_p \phi ' \psi _2 (x) f'(\rho) z $ and the second equation of \eqref{eq-invariants} by $E_p \phi ' \,  \psi _2 (x) f'(\xi) z$, we obtain the following second identity:

\begin{multline}
\label{mult-2nd-mult}
\int _S ^T E_p (t) \phi ' (t) \int _0 ^1 \psi _2 (x) \Bigl( \rho f(\rho) + \xi f(\xi) \Bigr) \, dx \, dt
\\
= \Bigl[ E_p (t) \phi ' (t) \int _0 ^1 \psi _2 (x) z \Bigl( f(\xi) - f(\rho)\Bigr) \, dx \Bigr] _S ^T 
\\
+ \int _S ^T (E_p (t) \phi '(t) )' \int _0 ^1 \psi _2 (x) z \Bigl( f(\rho) - f(\xi)\Bigr) \, dx \, dt
\\
- \int _S ^T E_p (t) \phi ' (t) \int _0 ^1 \psi _2 ' (x)  \, z \Bigl( f(\rho) + f(\xi)\Bigr) \, dx \, dt
\\
- \int _S ^T E_p (t) \phi ' (t) \int _0 ^1 \psi _2 (x)   a(x)  \, z  g(\frac{\rho - \xi}{2}) \Bigl( f'(\rho) + f'(\xi)\Bigr) \, dx \, dt
\\
+ 2 \int _S ^T E_p (t) \phi ' (t) \int _0 ^1 \psi _2 (x) \, \frac{\rho - \xi}{2} \Bigl( f(\rho) - f(\xi)\Bigr) \, dx \, dt .
\end{multline}

This second identity \eqref{mult-2nd-mult} allows us to obtain the following second estimate:
 
\begin{Lemma}
\label{lem-2nd-estim}
There exists $C_2 ' >0$ independent of the initial conditions such that
\begin{multline}
\label{mult-2nd-mult-cons}
\int _S ^T E_p (t) ^2 \phi ' (t) \, dt
\leq C_2 ' E_p (S) ^2  \phi ' (S) 
+ C_2 ' \int _S ^T E_p (t) \phi '(t)  \Bigl( \int _{a_2} ^{b_2} \vert z \vert ^p \, dx \Bigr) \, dt
\\
+ 2 C_1 \int _S ^T E_p (t) \phi ' (t) \int _0 ^1 \psi _2 (x) \, \frac{\rho - \xi}{2} \Bigl( f(\rho) - f(\xi)\Bigr) \, dx \, dt
\\
- 2\int _S ^T E_p (t) \phi ' (t) \int _0 ^1 x a(x) \psi _1 (x) g(\frac{\rho - \xi}{2}) \Bigl( f(\rho) + f(\xi) \Bigr) \, dx \, dt .
\end{multline}

\end{Lemma}


\subsubsection{Third multiplier, third identity, third estimate} \hfill

Given $t >0$, consider the solution $v^{(t)}$ of the elliptic problem

\begin{equation}
\label{def-v}
\begin{cases}
v_{xx} = \psi _3 (x) f(z(t,x)) \quad & x\in (0,1) ,\\
v(0)=0=v(1) .
\end{cases}
\end{equation}

Multiplying the first equation and the second equation of \eqref{eq-invariants} by $E_p \phi ' v $, we obtain the following third identity:

\begin{multline}
\label{mult-3rd-mult}
2 \int _S ^T E_p \phi ' \int _0 ^1 \psi _3 (x) \vert z \vert ^p 
= \Bigl[ E_p \phi ' \int _0 ^1 v (\rho - \xi) \Bigr] _S ^T 
\\
- \int _S ^T (E_p \phi ')' \int _0 ^1  v (\rho - \xi)
-  \int _S ^T E_p \phi ' \int _0 ^1  v_t  (\rho - \xi)
\\
+ 2 \int _S ^T E_p \phi ' \int _0 ^1 a (x) \, v g(\frac{\rho - \xi}{2}) .
\end{multline}

This third identity \eqref{mult-3rd-mult} allows us to obtain the following third estimate:

\begin{Lemma}
\label{lem-3rd-estim}
There exists $C_4>0$ independent of the initial conditions such that
\begin{multline}
\label{mult-3rd-mult-cons}
\int _S ^T E_p ^2 \phi ' 
\leq C_4 E_p (S) ^2  \phi ' (S) + C_4 \int _S ^T E_p \phi ' \int _0 ^1 \psi _3 (x) \, \vert \rho - \xi \vert  ^p \, dx
\\
+ 4C_3 \int _S ^T E_p \phi ' \int _0 ^1 a (x) \, v g(\frac{\rho - \xi}{2}) 
+ C_4 \int _S ^T E_p \phi ' \int _0 ^1 \psi _2 (x) \, \frac{\rho - \xi}{2} \Bigl( f(\rho) - f(\xi)\Bigr)
\\
- 4\int _S ^T E_p \phi ' \int _0 ^1 x a(x) \psi _1 (x) g(\frac{\rho - \xi}{2}) \Bigl( f(\rho) + f(\xi) \Bigr) .
\end{multline}

\end{Lemma}


\subsection{Proof of Lemma \ref{lem-int-estim} assuming Lemmas \ref{lem-1st-estim}-\ref{lem-3rd-estim} } \hfill

The starting point is \eqref{mult-3rd-mult-cons}, and the goal is to estimate the second, third, and fifth terms.

\underline{The third term of the right hand side of \eqref{mult-3rd-mult-cons}:} using \eqref{estim-v-Lq} and the fact that $g$ is globally Lipschitz, we have
\begin{multline*}
\Bigl\vert \int _0 ^1 a (x) \, v (t,x) g(\frac{\rho - \xi}{2}) \, dx \Bigr \vert
\leq C \Bigl( \int _0 ^1 \vert v (t,x) \vert ^q \, dx \Bigr) ^{1/q} \Bigl( \int _0 ^1 a(x) ^p \vert g(\frac{\rho - \xi}{2}) \vert ^p \, dx \Bigr) ^{1/p}
\\
\leq C' \Bigl( C C_p E_p (t) \Bigr) ^{1/q} \Bigl( \int _0 ^1 a(x)  \vert \rho - \xi \vert ^p \, dx \Bigr) ^{1/p}.
\end{multline*}
Then, as we did before, choosing $\sigma >0$, we have
$$
\Bigl\vert \int _0 ^1 a (x) \, v (t,x) g(\frac{\rho - \xi}{2}) \, dx \Bigr \vert
\leq C' \Bigl[ \frac{\sigma C C_p}{q} E_p (t) + \frac{1}{p \sigma ^{p/q}} \int _0 ^1 a(x)  \vert \rho - \xi \vert ^p \, dx \Bigr] .
$$
Therefore
\begin{multline*}
C_4 \int _S ^T E_p \phi ' \int _0 ^1 a (x) \, v g(\frac{\rho - \xi}{2}) 
\leq C_4 C' \frac{\sigma C C_p}{q} \int _S ^T E_p ^2  \phi ' \, dt 
\\
+ \frac{C_4 C'}{p \sigma ^{p/q}} \int _S ^T E_p \phi ' \int _0 ^1 a(x)  \vert \rho - \xi \vert ^p \, dx \, dt ,
\end{multline*}
and choosing $\sigma$ such that $C_4 C' \frac{\sigma C C_p}{q} = \frac{1}{2}$, and remembering that
$\psi _3 \leq C a$, we obtain that
\begin{multline}
\label{mult-1st-mult-cons4}
\int _S ^T E_p ^2 \phi ' 
\leq 2 C_4 E_p (S) ^2  \phi ' (S) + C_4 ' \int _S ^T E_p \phi ' \int _0 ^1 a (x) \, \vert \rho - \xi \vert  ^p \, dx
\\
+ C_4 \int _S ^T E_p \phi ' \int _0 ^1 \psi _2 (x) \, \frac{\rho - \xi}{2} \Bigl( f(\rho) - f(\xi)\Bigr)
\\
- 4\int _S ^T E_p \phi ' \int _0 ^1 x a(x) \psi _1 (x) g(\frac{\rho - \xi}{2}) \Bigl( f(\rho) + f(\xi) \Bigr) .
\end{multline}

\underline{The fifth term of the right-hand side of \eqref{mult-3rd-mult-cons}:}
We note that
\begin{multline*}
- 4\int _S ^T E_p \phi ' \int _0 ^1 x a(x) \psi _1 (x) g(\frac{\rho - \xi}{2}) \Bigl( f(\rho) + f(\xi) \Bigr)
\\
= 4\int _S ^T E_p \phi ' \int _0 ^1 x a(x) \psi _1 (x) g(\frac{\rho - \xi}{2}) \Bigl( f(\rho) - f(\xi) -2f(\rho) \Bigr)
\\
= 4\int _S ^T E_p \phi ' \int _0 ^1 x a(x) \psi _1 (x) g(\frac{\rho - \xi}{2}) \Bigl( f(\rho) - f(\xi) \Bigr)
\\
+ 4\int _S ^T E_p \phi ' \int _0 ^1 x a(x) \psi _1 (x) g(\frac{\rho - \xi}{2}) \Bigl( -2f(\rho) \Bigr)
\\
\leq 4 \int _S ^T E_p \phi ' \int _0 ^1  a(x)  g(\frac{\rho - \xi}{2}) \Bigl( f(\rho) - f(\xi) \Bigr)
\\
+ 8 \int _S ^T E_p \phi ' \int _0 ^1 x a(x) \psi _1 (x) \vert g(\frac{\rho - \xi}{2}) \vert \, \vert \rho \vert ^{p-1} .
\end{multline*}
On the other hand, one has
\begin{multline*}
4 \int _S ^T E_p \phi ' \int _0 ^1  a(x)  g(\frac{\rho - \xi}{2}) \Bigl( f(\rho) - f(\xi) \Bigr)
= 4 \int _S ^T E_p \phi ' (-E_p ' ) \, dt
\\
\leq 4 \phi '(S) \int _S ^T  -E_p 'E_p  \, dt
= 4 \phi '(S) \Bigl[ \frac{-1}{2} E_p (t) ^2 \Bigr] _S ^T 
\leq 2 E_p (S) ^2 \phi '(S) ,
\end{multline*}
and also
\begin{multline*}
8 \int _S ^T E_p \phi ' \int _0 ^1 x a(x) \psi _1 (x) \vert g(\frac{\rho - \xi}{2}) \vert \, \vert \rho \vert ^{p-1} 
\leq 8 C_g \int _S ^T E_p \phi '\int _0 ^1 a(x) \vert \frac{\rho - \xi}{2} \vert \, \vert \rho \vert ^{p-1} 
\\
\leq 4 C_g \int _S ^T E_p \phi ' \Bigl( \int _0 ^1 a(x)^p \vert \rho - \xi \vert ^p \Bigr) ^{1/p} 
\Bigl( \int _0 ^1 \vert \rho \vert ^{p} \Bigr) ^{1/q} \, dt
\\
\leq 4 C_g \int _S ^T E_p \phi ' \Bigl( \int _0 ^1 a(x)^p \vert \rho - \xi \vert ^p \Bigr) ^{1/p} 
\Bigl( E_p(t) \Bigr) ^{1/q} \, dt
\\
\leq \frac{1}{4} \int _S ^T E_p ^2 \phi ' \, dt + C \int _S ^T E_p \phi ' \int _0 ^1 a(x) \vert \rho - \xi \vert ^p \, dx \, dt .
\end{multline*}
Therefore
\begin{multline*}
- 4\int _S ^T E_p \phi ' \int _0 ^1 x a(x) \psi _1 (x) g(\frac{\rho - \xi}{2}) \Bigl( f(\rho) + f(\xi) \Bigr)
\\
\leq 2 E_p (S) ^2 \phi '(S) + \frac{1}{4} \int _S ^T E_p ^2 \phi ' \, dt + C \int _S ^T E_p \phi ' \int _0 ^1 a(x) \vert \rho - \xi \vert ^p \, dx \, dt ,
\end{multline*}
and with \eqref{mult-1st-mult-cons4}, we obtain that
\begin{multline}
\label{mult-1st-mult-cons5}
\int _S ^T E_p ^2 \phi ' 
\leq 3 C_4 E_p (S) ^2  \phi ' (S) + C_4 '' \int _S ^T E_p \phi ' \int _0 ^1 a (x) \, \vert \rho - \xi \vert  ^p \, dx
\\
+ C_4 \int _S ^T E_p \phi ' \int _0 ^1 \psi _2 (x) \, \frac{\rho - \xi}{2} \Bigl( f(\rho) - f(\xi)\Bigr) .
\end{multline}

\underline{The second term of the right-hand side of \eqref{mult-3rd-mult-cons}:}
we claim that there exists $C>0$ such that
\begin{equation}
\label{estim-f-pasf}
\forall a, b \in \Bbb R, \quad \vert a - b \vert ^p \leq C (a - b) (f(a)-f(b)) .
\end{equation}
Indeed, \eqref{estim-f-pasf} is obvious if $ab=0$: assume that $b=0$, then $\vert a - b \vert ^p = \vert a \vert ^p = af(a) = (a - b) (f(a)-f(b))$. 

Next, \eqref{estim-f-pasf} also holds if $ab <0$: assume that $b<0<a$, and denote $d:=-b$; then
$$ \vert a - b \vert ^p = (a+d)^p \leq C (a^p + d^p),$$
and
$$ (a - b) (f(a)-f(b)) = (a +d) (f(a)+f(d)) = (a+d) (a^{p-1} + d^{p-1}) \geq a^p + d^p ,$$
which implies that \eqref{estim-f-pasf} holds if $ab<0$.

Last, \eqref{estim-f-pasf} also holds if $ab>0$. One can assume that $a$ and $b$ are positive and that $b< a$. Then $b=\theta a$ with some $\theta \in [0,1]$; and
$$  \vert a - b \vert ^p = a^p (1-\theta ) ^p , \quad (a - b) (f(a)-f(b))
= a^p (1-\theta) (1 -\theta ^{p-1}) .$$
Denote $s:= 1-\theta$.  The function $s\in (0,1] \mapsto \frac{s^{p-1}}{1-(1-s)^{p-1}}$ is continuous on $(0,1]$
and goes to $0$ as $s\to 0^+$, therefore this function is bounded on $(0,1]$. Hence there exists $C>0$ such that
$$ \forall s \in (0,1], \quad \frac{s^{p-1}}{1-(1-s)^{p-1}} \leq C ,$$
hence
$$ \forall s \in (0,1], \quad s^{p-1} \leq C \Bigl( 1-(1-s)^{p-1} \Bigr) .$$
Coming back to $\theta = 1-s$:
$$ \forall \theta \in [0,1), \quad (1-\theta) ^{p-1} \leq C \Bigl( 1-\theta ^{p-1} \Bigr) .$$
Therefore
$$ \forall \theta \in [0,1), \quad (1-\theta) ^{p} \leq C (1-\theta) \Bigl( 1-\theta ^{p-1} \Bigr) ,$$
which gives that
$$ \vert a - b \vert ^p = a^p (1-\theta ) ^p \leq C a^p (1-\theta) \Bigl( 1-\theta ^{p-1} \Bigr) = (a - b) (f(a)-f(b)) .$$
Therefore \eqref{estim-f-pasf} holds.

And then, we deduce from \eqref{estim-f-pasf} that
\begin{multline}
\label{estim-puissp}
 C_4 '' \int _S ^T E_p \phi ' \int _0 ^1 a (x) \, \vert \rho - \xi \vert  ^p \, dx
\\
\leq C C_4 '' \int _S ^T E_p \phi ' \int _0 ^1 a (x) \,  ( \rho - \xi ) (f(\rho) - f(\xi))\, dx .
\end{multline}
Since $\psi _3 \leq C a$ on $(0,1)$, \eqref{mult-1st-mult-cons5} and \eqref{estim-puissp} give \eqref{eq-int-estim}.
This concludes the proof of Lemma \ref{lem-int-estim}. \qed 


\subsection{Asymptotic stability: Proof of Theorem \ref{thm-stab}, part (i)} \hfill

We recall that we know from \eqref{E_p'} that
$$
E_p'(t) 
=  -\int_0^1 a(x)g\left(\frac{\rho - \xi}{2}\right)\left( f(\rho)-f(\xi) \right) dx  .
$$
Now assume that there exists $0 < \alpha < \beta$ such that
$$ \forall s \in \Bbb R, \quad \alpha s \leq g(s) \leq \beta s .$$
Then choose $\phi: \Bbb R_+ \to \Bbb R_+$, $\phi (t)=t$: Lemma \ref{lem-int-estim}
gives that
\begin{equation}
\label{eq-int-estim-lin}
\int _S ^T E_p (t) ^2 \, dt
\leq C_u E_p (S) ^2   + C_u \int _S ^T E_p (t) \int _0 ^1 a (x) \, \frac{\rho - \xi}{2} \Bigl( f(\rho) - f(\xi)\Bigr) \, dx \, dt .
\end{equation}
If $\rho > \xi$, we have
$$ 0 < \frac{\rho - \xi}{2} \leq \frac{1}{\alpha} g(\frac{\rho - \xi}{2}) 
\quad \text{ and } \quad \Bigl( f(\rho) - f(\xi)\Bigr) >0 ,$$
hence
$$ \frac{\rho - \xi}{2} \Bigl( f(\rho) - f(\xi)\Bigr) \leq \frac{1}{\alpha} g(\frac{\rho - \xi}{2}) \Bigl( f(\rho) - f(\xi)\Bigr) ;$$
if $\rho \leq \xi$, we have

$$ \frac{\rho - \xi}{2} \geq \frac{1}{\beta} g(\frac{\rho - \xi}{2}) 
\quad \text{ and } \quad \Bigl( f(\rho) - f(\xi)\Bigr) < 0 ,$$
hence 
$$ \frac{\rho - \xi}{2} \Bigl( f(\rho) - f(\xi)\Bigr) \leq \frac{1}{\beta} g(\frac{\rho - \xi}{2}) \Bigl( f(\rho) - f(\xi)\Bigr) \leq \frac{1}{\alpha} g(\frac{\rho - \xi}{2}) \Bigl( f(\rho) - f(\xi)\Bigr)  .$$
Therefore \eqref{eq-int-estim-lin} gives that

\begin{multline*}
\int _S ^T E_p  ^2
\leq C_u E_p (S) ^2   + \frac{C_u}{\alpha} \int _S ^T E_p \int _0 ^1 a (x) \, g(\frac{\rho - \xi}{2}) \, \Bigl( f(\rho) - f(\xi)\Bigr) 
\\
= C_u E_p (S) ^2   + \frac{C_u}{\alpha} \int _S ^T E_p (-E_p ') 
\\
=  C_u E_p (S) ^2   + \frac{C_u}{\alpha} \Bigl[ \frac{-1}{2} E_p (t) ^2 \Bigr] _S ^T 
\leq \Bigl( C_u + \frac{C_u}{2\alpha} \Bigr) E_p (S) ^2 ,
\end{multline*}
and then a classical Gronwall type inequality (\cite{Komornik1994}) gives that $E_p ^2$ decays exponentially to $0$, as proved in \cite{KMJC2022} in the linear case. \qed


\subsection{Asymptotic stability: Proof of Theorem \ref{thm-stab}, part (ii)} \hfill
\label{sec-thm-ii}

Consider the assumptions of Theorem \ref{thm-stab}, part (ii): $g'(0)=0$, and there is an increasing and odd function $g_0: [-1,1] \to \Bbb R$ and $0<\alpha \leq \beta$ such that
\eqref{eq-comp-g} holds.
Moreover, the function
\begin{equation}
\label{def-H}
H(s) := \frac{g_0(s)}{s} ,
\end{equation}
 is $C^1$, satisfies $H(0)=0$ and is increasing on $[0,\eta]$ for some $\eta >0$.

Now, assume that $\varepsilon: \Bbb R_+ \to \Bbb R_+$ is $C^1$, positive, $\varepsilon (0) \leq \eta$ and decreasing to $0$ as $t\to +\infty$.
The function $\varepsilon$ will be chosen later. And we choose
\begin{equation}
\label{def-phi'}
\phi '(t) = H (\varepsilon (t)) ,
\end{equation}
which is well defined and decreasing on $\Bbb R_+$.

\subsubsection{Decomposition in subdomains and application} \hfill

Given $t\geq 0$, consider the following partition of $(0,1)$:
\begin{equation}
\label{def-partition}
\begin{cases}
\Omega _1 ^{(t)} := \{x \in (0,1), \vert \rho (t,x) - \xi (t,x) \vert \leq 2 \varepsilon (t) \}, \\
\Omega _2 ^{(t)} := \{x \in (0,1), 2 \varepsilon (t) < \vert \rho (t,x) - \xi (t,x) \vert \leq 2\eta \}, \\
\Omega _3 ^{(t)} := \{x \in (0,1), \vert \rho (t,x) - \xi (t,x) \vert > 2\eta \} .
\end{cases}
\end{equation}
In the following we decompose the integral appearing in the right-hand side of \eqref{eq-int-estim} according to the above domain decomposition, 
\begin{multline*}
\int _S ^T E_p \phi ' \int _0 ^1 a (x) \, \frac{\rho - \xi}{2} \Bigl( f(\rho) - f(\xi)\Bigr) \, dx \, dt
\\
= \int _S ^T E_p \phi ' \int _{\Omega _1 ^{(t)}} a (x) \, \frac{\rho - \xi}{2} \Bigl( f(\rho) - f(\xi)\Bigr)  \, dx \, dt
\\
+ \int _S ^T E_p \phi ' \int _{\Omega _2 ^{(t)}} a (x) \, \frac{\rho - \xi}{2} \Bigl( f(\rho) - f(\xi)\Bigr)  \, dx \, dt
\\
+ \int _S ^T E_p \phi ' \int _{\Omega _3 ^{(t)}} a (x) \, \frac{\rho - \xi}{2} \Bigl( f(\rho) - f(\xi)\Bigr)  \, dx \, dt .
\end{multline*}

\underline{Estimate on $\Omega _3 ^{(t)}$.}
Here we have
$$
\frac{\rho - \xi}{2} \Bigl( f(\rho) - f(\xi)\Bigr) 
=  \Bigl\vert \frac{\rho - \xi}{2} \Bigr\vert  \, \Bigl\vert f(\rho) - f(\xi) \Bigr\vert .
$$
There exists $\alpha ' \in (0, \alpha ]$ such that
$$ \forall \vert s \vert  \geq \eta , \quad \vert g(s) \vert \geq \alpha ' \vert s \vert .$$
Therefore
$$
\Bigl\vert \frac{\rho - \xi}{2} \Bigr\vert  \, \Bigl\vert f(\rho) - f(\xi) \Bigr\vert
\leq  \frac{1}{\alpha '} \Bigl\vert g(\frac{\rho - \xi}{2}) \Bigr\vert  \, \Bigl\vert f(\rho) - f(\xi) \Bigr\vert
= \frac{1}{\alpha '} g(\frac{\rho - \xi}{2}) \, \Bigl( f(\rho) - f(\xi)\Bigr).
$$
Therefore
\begin{multline*}
\int _S ^T E_p \phi ' \int _{\Omega _3 ^{(t)}} a (x) \, \frac{\rho - \xi}{2} \Bigl( f(\rho) - f(\xi)\Bigr)  \, dx \, dt
\\
\leq \int _S ^T E_p \phi ' \int _{\Omega _3 ^{(t)}} a (x) \frac{1}{\alpha '} g(\frac{\rho - \xi}{2}) \, \Bigl( f(\rho) - f(\xi)\Bigr)  \, dx \, dt
\\
\leq \frac{1}{\alpha '} \int _S ^T E_p \phi ' \int _0 ^1 a (x) g(\frac{\rho - \xi}{2}) \, \Bigl( f(\rho) - f(\xi)\Bigr)  \, dx \, dt
\\
= \frac{1}{\alpha '} \int _S ^T E_p \phi ' (-E_p ') \, dt 
\leq \frac{1}{\alpha '} \phi '(S) \Bigl[\frac{-1}{2} E_p (t) ^2 \Bigr] _S ^T ,
\end{multline*}
hence
\begin{equation}
\label{omega-1-t}
\int _S ^T E_p \phi ' \int _{\Omega _3 ^{(t)}} a (x) \, \frac{\rho - \xi}{2} \Bigl( f(\rho) - f(\xi)\Bigr)  \, dx \, dt
\leq \frac{1}{2 \alpha '} E_p (S) ^2 \phi '(S) .
\end{equation}

\underline{Estimate on $\Omega _2 ^{(t)}$.} If $x\in \Omega _2 ^{(t)}$, then $ \Bigl\vert \frac{\rho - \xi}{2} \Bigr\vert \in 
[\varepsilon (t), \eta]$, hence
$$ H\Bigl( \Bigl\vert \frac{\rho - \xi}{2} \Bigr\vert \Bigr) \geq H(\varepsilon (t)) = \phi '(t) ,$$
which brings
$$ g_0 \Bigl( \Bigl\vert \frac{\rho - \xi}{2} \Bigr\vert \Bigr) \geq \phi '(t) \, \Bigl\vert \frac{\rho - \xi}{2} \Bigr\vert .$$
We deduce that
$$ \forall x\in \Omega _2 ^{(t)}, \quad \phi '(t) \, \Bigl\vert \frac{\rho - \xi}{2} \Bigr\vert
\leq g_0 \Bigl( \Bigl\vert \frac{\rho - \xi}{2} \Bigr\vert \Bigr)
= \Bigl\vert g_0 \Bigl( \frac{\rho - \xi}{2} \Bigr) \Bigr\vert
\leq \Bigl\vert g \Bigl( \frac{\rho - \xi}{2} \Bigr) \Bigr\vert .$$
Then
\begin{multline*}
\int _S ^T E_p \phi ' \int _{\Omega _2 ^{(t)}} a (x) \, \frac{\rho - \xi}{2} \Bigl( f(\rho) - f(\xi)\Bigr)  \, dx \, dt
\\
= \int _S ^T E_p \int _{\Omega _2 ^{(t)}} a (x) \, \phi '(t) \, \Bigl\vert \frac{\rho - \xi}{2} \Bigr\vert \, \Bigl\vert f(\rho) - f(\xi) \Bigr\vert  \, dx \, dt
\\
\leq \int _S ^T E_p \int _{\Omega _2 ^{(t)}} a (x) \, \Bigl\vert g \Bigl( \frac{\rho - \xi}{2} \Bigr) \, \Bigl\vert \Bigl\vert f(\rho) - f(\xi) \Bigr\vert  \, dx \, dt
\\
= \int _S ^T E_p \int _{\Omega _2 ^{(t)}} a (x) \,  g \Bigl( \frac{\rho - \xi}{2} \Bigr) \Bigl( f(\rho) - f(\xi) \Bigr)  \, dx \, dt
\\
\leq \int _S ^T E_p \int _0 ^1 a (x) \,  g \Bigl( \frac{\rho - \xi}{2} \Bigr) \Bigl( f(\rho) - f(\xi) \Bigr)  \, dx \, dt
= \int _S ^T E_p (-E_p') \, dt .
\end{multline*}
Therefore
\begin{equation}
\label{omega-2-t}
\int _S ^T E_p \phi ' \int _{\Omega _1 ^{(t)}} a (x) \, \frac{\rho - \xi}{2} \Bigl( f(\rho) - f(\xi)\Bigr)  \, dx \, dt
\leq \frac{1}{2} E_p (S)^2 .
\end{equation}

\underline{Estimate on $\Omega _1 ^{(t)}$.} Assume that $\vert \rho (t,x) - \xi (t,x) \vert \leq 2 \varepsilon (t)$.
If $\rho \geq \xi$, then the mean value theorem gives us that
$$  \Bigl\vert f(\rho) - f(\xi) \Bigr\vert
\leq \Bigl( \sup _{(\xi, \xi + 2 \varepsilon (t))} \vert f' \vert \Bigr) \, \vert \rho  - \xi  \vert
\leq (p-1) \Bigl( \vert \xi \vert + 2 \varepsilon (t) \Bigr) ^{p-2} \, 2 \varepsilon (t) .$$
And if $\rho \leq \xi$, then
$$  \Bigl\vert f(\rho) - f(\xi) \Bigr\vert \leq (p-1) \Bigl( \vert \rho \vert + 2 \varepsilon (t) \Bigr) ^{p-2} \, 2 \varepsilon (t) .$$
Then, in any case, we have
$$ \Bigl\vert f(\rho) - f(\xi) \Bigr\vert \leq (p-1) \Bigl( \vert \rho \vert + \vert \xi \vert + 2 \varepsilon (t) \Bigr) ^{p-2} \, 2 \varepsilon (t) .$$
Then we obtain that
\begin{multline*}
\int _S ^T E_p \phi ' \int _{\Omega _1 ^{(t)}} a (x) \, \frac{\rho - \xi}{2} \Bigl( f(\rho) - f(\xi)\Bigr)  \, dx \, dt
\\
\leq \int _S ^T E_p \phi ' \int _{\Omega _1 ^{(t)}} 2 (p-1) \, a (x)\,  \varepsilon (t)^2 \Bigl( \vert \rho \vert + \vert \xi \vert + 2 \varepsilon (t) \Bigr) ^{p-2}  \, dx \, dt
\\
\leq C \int _S ^T E_p (t) \phi '(t)  \varepsilon (t)^2 \Bigl( \int _{\Omega _1 ^{(t)}} \Bigl( \vert \rho \vert + \vert \xi \vert + 2 \varepsilon (t) \Bigr) ^{p}  \, dx \Bigr) ^{(p-2)/p}\, dt
\\
\leq C' \int _S ^T E_p (t) \phi '(t)  \varepsilon (t)^2 \Bigl( \int _{\Omega _1 ^{(t)}}  \vert \rho \vert ^p + \vert \xi \vert ^p + \varepsilon (t) ^p  \, dx \Bigr) ^{(p-2)/p}\, dt
\\
\leq C' \int _S ^T E_p (t) \phi '(t)  \varepsilon (t)^2 \Bigl( p E_p (t) + \varepsilon (t) ^p  \Bigr) ^{(p-2)/p}\, dt
\\
\leq C'' \int _S ^T E_p (t) \phi '(t)  \varepsilon (t)^2 \Bigl( E_p (t)^{(p-2)/p} + \varepsilon (t) ^{p-2}  \Bigr)\, dt .
\end{multline*}
Therefore
\begin{multline}
\label{omega-3-t}
\int _S ^T E_p \phi ' \int _{\Omega _1 ^{(t)}} a (x) \, \frac{\rho - \xi}{2} \Bigl( f(\rho) - f(\xi)\Bigr)  \, dx \, dt
\\
\leq C'' \int _S ^T E_p (t) ^{1+ (p-2)/p} \phi '(t)  \varepsilon (t)^2 \, dt
+   C'' \int _S ^T E_p (t) \phi '(t)  \varepsilon (t)^p \, dt .
\end{multline}

\underline{Consequence on Lemma \ref{lem-int-estim}.} Using \eqref{omega-1-t}-\eqref{omega-3-t}, we deduce from \eqref{eq-int-estim} that
\begin{multline}
\label{eq-int-estim-cons*}
\int _S ^T E_p ^2 \phi ' 
\leq C' _u E_p (S) ^2 
\\
+ C' _u  \int _S ^T E_p (t) ^{1+ (p-2)/p} \phi '(t)  \varepsilon (t)^2 \, dt
+   C' _u  \int _S ^T E_p (t) \phi '(t)  \varepsilon (t)^p \, dt .
\end{multline}
In the following we choose the function $\varepsilon$ and we deduce some decay of the energy.

\subsubsection{Some possible choices for the functions $\phi$ and  $\varepsilon$ and their properties} \hfill

A way to obtain an estimate on the decay of the energy is to note that since $E_p$ is nonincreasing, we can derive from \eqref{eq-int-estim-cons*} that
\begin{multline}
\label{eq-int-estim-cons*2}
\int _S ^T E_p ^2 \phi ' 
\leq C' _u E_p (S) ^2 
\\
+ C' _u  E_p (S) ^{1+ (p-2)/p} \int _S ^T  \phi '(t)  \varepsilon (t)^2 \, dt
+   C' _u  E_p (S) \int _S ^T  \phi '(t)  \varepsilon (t)^p \, dt  .
\end{multline}
The trick is now to choose the function $\varepsilon$ such that we can compute the integrals appearing on the right-hand side of 
\eqref{eq-int-estim-cons*2}:
using \eqref{def-phi'}, we have 
$$\varepsilon (t) = H^{-1} (\phi '(t)),$$
which gives
\begin{multline*}
\int _S ^T  \phi '(t)  \varepsilon (t)^2 \, dt
= \int _S ^T  \Bigl[ H^{-1} (\phi '(t)) \Bigr] ^2 \phi '(t) \, dt
\\
= \int _{\phi (S)} ^{\phi (T)} \Bigl[ H^{-1} (\phi '(\phi ^{-1} (s))) \Bigr] ^2 \, ds
= \int _{\phi (S)} ^{\phi (T)} \Bigl[ H^{-1}  \Bigl( \frac{1}{(\phi ^{-1})' (s)} \Bigr) \Bigr]^2 \, ds .
\end{multline*}
Then, computations would be easy if, for some $m\geq 1$, we had
$$ H ^{-1} \Bigl( \frac{1}{(\phi ^{-1})' (s)} \Bigr) = \frac{1}{(s + A)^m}, $$
where $A$ is a constant. With that idea, 

\begin{itemize}
\item first, given $m\geq 1$, let us consider 
\begin{equation}
\label{choice-psi}
\psi _m (s) := \int _0 ^s \frac{1}{H(\frac{1}{(\sigma + A)^m})} \, d\sigma :
\end{equation}
the function $\psi _m$ is well defined, of class $C^2$, and 
$$ \psi _m '(s) = \frac{1}{H(\frac{1}{(s + A)^m})} >0 .$$
Hence $\psi _m$ and $\psi _m '$ are increasing, hence $\psi _m$ is convex; and $\psi _m ' (s) \to +\infty$ as $s\to +\infty$.

\item Next, let us consider 
\begin{equation}
\label{choice-phi}
\phi _m := \psi _m ^{-1} .
\end{equation}
Of course $\phi _m $ is $C^2$, increasing and concave. Moreover, fix $A = \frac{1}{\eta ^{1/m}}$, take $t \geq 0$ and denote $s:= \phi _m (t)$; then $t= \phi _m ^{-1} (s)= \psi _m (s)$, and
$$ \phi _m '(t) = \phi _m '(\psi _m (s)) = \frac{1}{\psi _m '(s)} = H(\frac{1}{(s + A)^m}).$$
In particular, $\phi _m '(0)= H(\frac{1}{ A^m }) = H(\eta)$, and $\phi _m '(t) \to 0$ as $t\to +\infty$, hence $\phi _m ' ([0,+\infty))
= (0, H(\eta)]$.

\item Finally, let us choose
\begin{equation}
\label{choice-eps}
\varepsilon _m (t) := H^{-1} (\phi _m '(t)):
\end{equation}
from the previous remarks, $\phi _m '$ takes its values in the region where $H^{-1}$ is well-defined, hence $\varepsilon _m$ is well-defined on $[0,+\infty)$, and $C^1$ by composition, and decreasing, with $\varepsilon _m (0)= \eta$, and $\varepsilon _m (t) \to 0$ as $t\to +\infty$.
\end{itemize}
Therefore, these choices have all the requirements we assumed to derive \eqref{eq-int-estim-cons*2}, and are going to allow us to estimate the decay of the energy $E_p$.

To conclude, we will need to have an estimate of the asymptotic behaviour of $\phi _m$ as $t\to +\infty$: it is related to the behaviour of $H$ (and thus of $g_0$) as $x\to 0^+$:

$$\forall s\geq A, \quad \psi _m (s) \leq \frac{s}{H(\frac{1}{(s + A)^m})} \leq \frac{s}{H(\frac{1}{(2s)^m})} \leq \frac{2s}{H(\frac{1}{(2s)^m})} ,$$
hence
\begin{equation}
\label{comp-infty-phi_m}
 \forall s \geq A, \quad s \leq \phi _m (t), \quad \text{ where } \quad t = \frac{2s}{H(\frac{1}{(2s)^m})} .
 \end{equation}
Then, 
\begin{itemize} 
\item With $m=1$, this brings
$$
t = \frac{2s}{\frac{g_0 (\frac{1}{2s})}{\frac{1}{2s}}} = \frac{1}{g_0 (\frac{1}{2s})} ,
$$
hence 
\begin{equation}
\label{growth-phi_1}
\phi _1 (t) \geq s = \frac{1}{2 \, g_0 ^{-1} (\frac{1}{t})} .
\end{equation}

\item If there is some $\alpha >0$ and some $q>1$ such that 
\begin{equation}
\label{assump-poly-0}
\forall x \in [0,\eta), \quad g_0 (x) \geq \alpha  x^q ,
\end{equation}
then 
$$ t = \frac{2s}{H(\frac{1}{(2s)^m})} \leq \frac{1}{\alpha} (2s) ^{m(q-1)+1} ,$$
therefore
\begin{equation}
\label{growth-phi_n}
\phi _m (t) \geq s \geq \frac{1}{2}  (\alpha t) ^{1/[m(q-1)+1]} .
\end{equation}
\end{itemize}
We end this section computing the two integrals that were the motivation for our choice of $\phi _m$ and $\varepsilon _m$: given $m\geq 1$, we have
\begin{multline}
\label{comp-int-2}
\int _S ^ {+\infty}  \phi _m '(t)  \varepsilon _m (t)^2 \, dt
= \int _{\phi _m (S)} ^{+\infty} \Bigl[ H^{-1}  \Bigl( \frac{1}{(\phi _m ^{-1})' (s)} \Bigr) \Bigr]^2 \, ds 
\\
= \int _{\phi _m (S)} ^{+\infty} \Bigl( \frac{1}{(s+A)^m} \Bigr) ^2 \, ds 
= \Bigl[ \frac{(s+A)^{-2m+1}}{-2m+1} \Bigr] _{\phi _m (S)} ^{+\infty}
\\
= \frac{1}{2m-1} \frac{1}{(\phi _m (S)+A) ^{2m-1}} .
\end{multline}
In the same way,
\begin{equation}
\label{comp-int-p}
\int _S ^ {+\infty}  \phi _m '(t)  \varepsilon _m (t)^p \, dt
= \frac{1}{pm-1} \frac{1}{(\phi _m (S)+A) ^{pm-1}} . 
\end{equation}


\subsubsection{A first estimate on the decay of the energy, using $\phi _1$} \hfill

Let us choose $m=1$, $\phi := \phi _1$ and $\varepsilon := \varepsilon _1$. 
Then we derive from \eqref{eq-int-estim-cons*2}, \eqref{comp-int-2} (with m=1) and \eqref{comp-int-p} (with m=1) that
\begin{equation}
\label{eq-int-estim-cons*3}
\int _S ^T E_p ^2 \phi ' 
\leq C' _u E_p (S) ^2 
+ C' _u \frac{E_p (S) ^{1+ (p-2)/p} }{\phi (S)+A}
+   C' _u \frac{E_p (S) }{(\phi (S)+A)^{p-1}} .
\end{equation}
This implies the following (weakened) estimate
\begin{multline}
\label{eq-int-estim-cons*4}
\int _S ^T E_p (t) ^2 \phi ' (t) \, dt
\\
\leq C' _u E_p (S) ^2 
+  C' _u E_p (0) ^{(p-2)/p} \frac{ E(S) }{\phi (S)+A}
+   C' _u \frac{E_p (S) }{\phi (S)+A} 
\\
= C' _u E_p (S) ^2 
+ \Bigl(C' _u  E_p (0) ^{(p-2)/p} + C_u '\Bigr) \frac{ E(S) }{\phi (S)+A} ,
\end{multline}
and we can now apply an integral estimate of the Gronwall type (see \cite{PM-COCV}):

\begin{Lemma} (Lemma 2 in \cite{PM-COCV})
\label{lem2-COCV}
Let $f: \Bbb R_+ \to\Bbb R_+$ be a nonincreasing continuous function. Assume that there exists $\sigma >0$, $\sigma ' >0$ and $c>0$ such that
\begin{equation}
\label{hyp-COCV}
\forall t \geq 0, \quad \int _t ^{+\infty} f(\tau) ^{1+\sigma} \, d\tau \leq c f(t)^{1+\sigma} + \frac{c}{(1+t)^{\sigma '}} f(0)^\sigma f(t) .
\end{equation}
Then there exists $C>0$ such that
\begin{equation}
\label{conc-COCV}
\forall t \geq 0, \quad f(t) \leq  f(0) \frac{C}{(1+t)^{(1+\sigma ')/\sigma}} .
\end{equation}
\end{Lemma}

\noindent As an immediate corollary, we have the following.

\begin{Lemma} (Lemma 3 in \cite{PM-COCV})
\label{lem3-COCV}
Let $E: \Bbb R_+ \to\Bbb R_+$ be a nonincreasing continuous function and $\phi: \Bbb R_+ \to\Bbb R_+$ an increasing function of class $C^1$ such that
$$ \phi (0)=0 \quad \text{ and } \quad \phi (t) \to +\infty \text{ as } t \to +\infty . $$
Assume that there exists $\sigma >0$, $\sigma ' >0$ and $c>0$ such that
\begin{equation}
\label{hyp-COCV2}
\forall S \geq 0, \quad \int _S ^{+\infty} E(t) ^{1+\sigma} \phi '(t) \, dt \leq c E(S)^{1+\sigma} + \frac{c}{(1+\phi (S))^{\sigma '}} E(0)^\sigma E(S) .
\end{equation}
Then there exists $C>0$ such that
\begin{equation}
\label{conc-COCV2}
\forall t \geq 0, \quad E (t) \leq  E(0) \frac{C}{(1+\phi (t))^{(1+\sigma ')/\sigma}} .
\end{equation}
\end{Lemma}

Then Lemma \ref{lem3-COCV} gives that there exists $C(E_p(0))$ that depends on $E_p(0)$ such that
\begin{equation}
\label{estim-decay1}
\forall t \geq 0, \quad E_p (t) \leq  \frac{C(E_p(0)) }{(\phi (S)+A)^2} .
\end{equation}
Using \eqref{growth-phi_1}, we obtain that
\begin{equation}
\label{estim-decay1*}
\forall t \geq 0, \quad E_p (t) \leq  C(E_p(0)) \Bigl( g_0 ^{-1} (\frac{1}{t}) \Bigr) ^2 .
\end{equation}
This is the first estimate of the decay of the energy.


\subsubsection{A new integral inequality and its application} \hfill

We claim that we have the following more adapted result:

\begin{Lemma}
\label{lem-COCV}
Let $f: \Bbb R_+ \to\Bbb R_+$ be a nonincreasing continuous function. Assume that there exists $p>2$, and $c>0$ such that
\begin{equation}
\label{eq-int-estim-cons*3*}
\int _t ^{+\infty} f(\tau) ^2 \, d\tau
\leq c f(t) ^2
+ c \frac{f(t) ^{1+ (p-2)/p} }{1+t}
+   c \frac{f(t) }{(1+t)^{p-1}} .
\end{equation}
Then, for all $\delta >0$, there exists $C_\delta >0$ such that
\begin{equation}
\label{conc-COCV3}
\forall t \geq 0, \quad f(t) \leq  f(0) \frac{C_\delta }{(1+t)^{p-\delta}} .
\end{equation}
\end{Lemma}

\begin{Remark} {\rm We were not able to prove that \eqref{conc-COCV3} is true for $\delta =0$. This would be interesting in itself, and would of course give a better estimate on the decay of the energy $E_p$.}
\end{Remark}

\noindent {\it Proof of Lemma \ref{lem-COCV}.}
Denote $\theta := \frac{p-2}{p}$. Assumption \eqref{eq-int-estim-cons*3*} implies that
\begin{equation}
\label{conc-COCV-step1.1}
\int _t ^{+\infty} f(\tau) ^2 \, d\tau
\leq c f(t) ^2
+ c f(0) ^\theta \frac{f(t)}{1+t} + 
+   c \frac{f(t) }{1+t} ,
\end{equation}
therefore we are in position to apply Lemma \ref{lem2-COCV} with $\sigma =1=\sigma '$, and we obtain that there exists $C_1$ such that
\begin{equation}
\label{conc-COCV-step1.2}
f(t) \leq \frac{C_1}{(1+t)^2} .
\end{equation}
This is the first estimate, that we can plug into \eqref{eq-int-estim-cons*3*} to obtain that
\begin{equation}
\label{conc-COCV-step2.1}
\int _t ^{+\infty} f(\tau) ^2 \, d\tau
\leq c f(t) ^2
+ c C_1 ^{\theta} \frac{f(t) }{(1+t)^{1+2\theta}}
+   c \frac{f(t) }{(1+t)^{p-1}} .
\end{equation}
Before using Lemma \ref{lem2-COCV}, we need to compare $1+2\theta$ and $p-1$; we see that
$$ 1+2\theta < p-1:$$
indeed,
$$ (p-1)- (1+2\theta) = p- 2 - 2\frac{p-2}{p} = \frac{p^2 - 4p + 4}{p} 
= \frac{(p-2)^2}{p} >0 .$$
Therefore we deduce from \eqref{conc-COCV-step2.1} that
$$ \int _t ^{+\infty} f(\tau) ^2 \, d\tau
\leq c f(t) ^2
+ (c C_1 ^{\theta} + c) \frac{f(t) }{(1+t)^{1+2\theta}} ,$$
therefore we are in position to apply Lemma \ref{lem2-COCV} with $\sigma =1$, and $\sigma '=1+2\theta$, and we obtain that there exists $C_2$ such that
\begin{equation}
\label{conc-COCV-step2.2}
f(t) \leq \frac{C_2}{(1+t)^{2+2\theta}} ,
\end{equation}
which improves the decay estimate \eqref{conc-COCV-step1.2}.

Let us prove by induction that, for all $n \in \Bbb N$, there exists $C_n$ such that
\begin{equation}
\label{conc-COCV-stepn.n}
(H_n): \quad f(t) \leq \frac{C_n}{(1+t)^{2 \sum _{k=0} ^n \theta ^k}} .
\end{equation}
$(H_0)$ is true, thanks to \eqref{conc-COCV-step2.1}, and $(H_1)$ is true, thanks to \eqref{conc-COCV-step2.2}.
If $(H_n)$ is true, then we derive from 
\eqref{eq-int-estim-cons*3} to obtain that
\begin{multline}
\label{conc-COCV-stepn.1}
\int _t ^{+\infty} f(\tau) ^2 \, d\tau
\leq c f(t) ^2
+ c C_n ^{\theta} \frac{f(t) }{(1+t)\, (1+t)^{2\theta \sum _{k=0} ^n \theta ^k}}
+   c \frac{f(t) }{(1+t)^{p-1}} 
\\
= c f(t) ^2
+ c C_n ^{\theta} \frac{f(t) }{(1+t)^{1+2 \sum _{k=1} ^{n+1} \theta ^k}}
+   c \frac{f(t) }{(1+t)^{p-1}} .
\end{multline}
In the same way, we need to compare the exponents, and we remark that
\begin{multline*} 
(p-1) - \Bigl( 1+2 \sum _{k=1} ^{n+1} \theta ^k \Bigr) 
= p - 2 \sum _{k=0} ^{n+1} \theta ^k 
= p - 2 \frac{1 - \theta ^{n+2}}{1-\theta} 
\\
= p - 2 \frac{1 - \theta ^{n+2}}{1-\frac{p-2}{p}}  
=  p - 2 \frac{1 - \theta ^{n+2}}{\frac{2}{p}}
= p - p (1 - \theta ^{n+2}) = p \theta ^{n+2} >0 .
\end{multline*}
Therefore we deduce from \eqref{conc-COCV-step2.1} that
$$ \int _t ^{+\infty} f(\tau) ^2 \, d\tau
\leq c f(t) ^2
+ (c C_n ^{\theta}+c) \frac{f(t)}{(1+t)^{1+2 \sum _{k=1} ^{n+1} \theta ^k}} ,$$
therefore we are in position to apply Lemma \ref{lem2-COCV} with $\sigma =1$, and $\sigma '=1+2 \sum _{k=1} ^{n+1} \theta ^k$, and we obtain that there exists $C_{n+1}$ such that
\begin{equation}
\label{conc-COCV-stepn.2}
f(t) \leq \frac{C_{n+1}}{(1+t)^{2+2 \sum _{k=1} ^{n+1} \theta ^k}} 
= \frac{C_{n+1}}{(1+t)^{2 \sum _{k=0} ^{n+1} \theta ^k}}
\end{equation}
which is \eqref{conc-COCV-stepn.n} at the rank $n+1$. This concludes the end of the induction and proves the validity of \eqref{conc-COCV-stepn.n} for all $n\in \Bbb N$.
Since we already noted that
$$ 2 \sum _{k=0} ^n \theta ^k = 2 \frac{1 - \theta ^{n+1}}{1-\theta}
= 2 \frac{1 - \theta ^{n+1}}{1-\frac{p-2}{p}}
= 2 \frac{1 - \theta ^{n+2}}{\frac{2}{p}}
= p (1 - \theta ^{n+2}) \to p \quad \text{ as } n\to \infty ,$$
we obtain \eqref{conc-COCV3}. \qed

Then Lemma \ref{lem-COCV} can be applied to \eqref{eq-int-estim-cons*3}, and we obtain that

\begin{equation*}
\forall \delta >0, \exists C_\delta, \forall t \geq 0, \quad E_p (t) \leq  \frac{C_\delta (E_p(0)) }{(\phi _1 (S)+A)^{p-\delta}} 
\leq  C_\delta (E_p(0)) \Bigl( g_0 ^{-1} (\frac{1}{t}) \Bigr) ^{p-\delta} ,
\end{equation*}
hence \eqref{decay-gene}. \qed


\subsubsection{An extension of the new integral inequality and its application} \hfill

Now, let us choose $m\geq 1$,  $\phi := \phi _m$ (defined in \eqref{choice-phi}), and $\varepsilon := \varepsilon _m$ (defined in \eqref{choice-eps}). 
We derive from \eqref{eq-int-estim-cons*2}, \eqref{comp-int-2} and \eqref{comp-int-p} that
\begin{equation}
\label{eq-int-estim-cons*3m}
\int _S ^T E_p ^2 \phi ' 
\leq C' _u E_p (S) ^2 
+ C' _u \frac{E_p (S) ^{1+ (p-2)/p} }{(\phi (S)+A)^{2m-1}}
+   C' _u \frac{E_p (S) }{(\phi (S)+A)^{mp-1}} .
\end{equation}
We claim that we have the following extension of Lemma \ref{lem-COCV}:

\begin{Lemma}
\label{lem-COCV-m}
Let $f: \Bbb R_+ \to\Bbb R_+$ be a nonincreasing continuous function. Assume that there exist $m\geq 1$, $p>2$, and $c>0$ such that
\begin{equation}
\label{eq-int-estim-cons*3*-m}
\int _t ^{+\infty} f(\tau) ^2 \, d\tau
\leq c f(t) ^2
+ c \frac{f(t) ^{1+ (p-2)/p} }{(1+t)^{2m-1}}
+   c \frac{f(t) }{(1+t)^{mp-1}} .
\end{equation}
Then, for all $\delta >0$, there exists $C_{m,\delta} >0$ such that
\begin{equation}
\label{conc-COCV3-m}
\forall t \geq 0, \quad f(t) \leq  f(0) \frac{C_{m,\delta} }{(1+t)^{m(p-\delta)}} .
\end{equation}
\end{Lemma}

\noindent {\it Proof of Lemma \ref{lem-COCV-m}.} The proof is similar to the proof of Lemma \ref{lem-COCV}, the only difference being the induction formula: one has to replace \eqref{conc-COCV-stepn.n} by
\begin{equation}
\label{conc-COCV-stepn.n-m}
(H_n): \quad f(t) \leq \frac{C_n}{(1+t)^{2 m \sum _{k=0} ^n \theta ^k}} ,
\end{equation}
and \eqref{conc-COCV3-m} follows in the same way. \qed

Then Lemma \ref{lem-COCV-m} can be applied to \eqref{eq-int-estim-cons*3m}, and we obtain that

\begin{equation}
\label{estim-decay1**-m}
\forall \delta >0, \exists C_\delta, \forall t \geq 0, \quad E_p (t) \leq  \frac{C_{m,\delta} (E_p(0)) }{(\phi _m (t)+A)^{m(p-\delta)}} .
\end{equation}
Then \eqref{decay-gene-m} follows from the estimate \eqref{comp-infty-phi_m} on the growth of $\phi _m$. \qed

Finally, assume that \eqref{assump-poly-0} holds true.
Then we deduce from the growth estimate of $\phi _m$ given in \eqref{growth-phi_n} that
\begin{equation}
\label{growth-phi_n-appl}
E_p (t) \leq  \frac{C_{m,\delta} (E_p(0)) }{(1+t) ^{\frac{m(p-\delta)}{m(q-1)+1}}}
\end{equation}
Taking $m$ large enough, we find \eqref{decay-poly-lim}. \qed


\subsection{Proof of Lemmas \ref{lem-1st-estim}-\ref{lem-3rd-estim}} \hfill

\subsubsection{Proof of Lemma \ref{lem-1st-estim}} \hfill

The starting point is \eqref{mult-1st-mult}, obtained multiplying the first equation of \eqref{eq-invariants} by $E_p \phi ' x \psi_1 (x) f(\rho)$ and the second equation of \eqref{eq-invariants} by $E_p \phi ' \,  x \psi _1 (x) f(\xi)$,
as in \cite{KMJC2022} (the only difference being here the presence of the multiplier $\phi'$.)

Then, first, we estimate the integrated term and we get, for $t \geq 0$ that
\begin{multline*}
 \Bigl\vert E_p (t) \phi ' (t) \int _0 ^1 x \psi _1 (x) \Bigl( F(\xi (t,x)) - F(\rho (t,x)) \Bigr) \, dx \Bigr\vert
\\
\leq E_p (t) \phi ' (t) \int _0 ^1  F(\xi (t,x)) + F(\rho (t,x)) \, dx
\\
= \frac{1}{p} E_p (t) \phi ' (t) \int _0 ^1  \vert \xi (t,x) \vert ^p  + \vert \rho (t,x) \vert ^p \, dx
= E_p (t) ^2 \phi ' (t) .
\end{multline*}
Therefore
$$ \Bigl\vert \Bigl[  E_p \phi ' \int _0 ^1 x \psi _1 (x) \Bigl( F(\xi) - F(\rho) \Bigr) \Bigr] _S ^T \Bigr\vert 
\leq E_p (T) ^2 \phi ' (T) + E_p (S) ^2 \phi ' (S),$$
and recalling that $\phi'$ and $E_p$ are nonincreasing, we obtain

\begin{equation}
\label{estim1-1}
\Bigl\vert \Bigl[  E_p \phi ' \int _0 ^1 x \psi _1 (x) \Bigl( F(\xi) - F(\rho) \Bigr) \Bigr] _S ^T \Bigr\vert 
\leq 2 E_p (S) ^2 \phi ' (S) .
\end{equation}

Next, in the same way, we have
\begin{multline*}
\Bigl\vert \int _S ^T (E_p (t) \phi ' (t) )' \int _0 ^1 x \psi _1 (x) \Bigl( F(\rho) - F(\xi) \Bigr) \, dx \, dt \Bigr\vert
\\
\leq \int _S ^T \Bigl\vert (E_p (t) \phi ' (t) )' \Bigr\vert \int _0 ^1 \Bigl( F(\rho) + F(\xi) \Bigr) \, dx \, dt
\\
= \int _S ^T \Bigl\vert E_p ' (t) \phi ' (t) + E_p (t) \phi ''  (t)\Bigr\vert E_p (t) \, dt
\\
= \int _S ^T \Bigl( - E_p ' (t) \phi ' (t) - E_p (t) \phi ''(t) \Bigr) E_p (t) \, dt
\\
= \int _S ^T \Bigl( - E_p ' (t)  E_p (t) \Bigr) \phi ' (t)  \, dt 
+ \int _S ^T E_p (t) ^2  \Bigl( -  \phi '' (t) \Bigr) \, dt
\\
\leq \phi ' (S) \int _S ^T \Bigl( - E_p ' (t)  E_p (t) \Bigr) \, dt 
+ E_p (S) ^2 \int _S ^T  \Bigl( -  \phi '' (t) \Bigr) \, dt
\\
= \phi ' (S) \Bigl[ \frac{-1}{2} E_p (t) ^2 \Bigr] _S ^T + E_p (S) ^2 \Bigl[ - \phi '(t) \Bigr] _S ^T
\\
= \phi ' (S) \Bigl( \frac{1}{2} E_p (S) ^2 - \frac{1}{2} E_p (T) ^2 \Bigr)
+ E_p (S) ^2  \Bigl( \phi '(S) - \phi '(T) \Bigr) .
\end{multline*}
Therefore
\begin{equation}
\label{estim1-2}
\Bigl\vert \int _S ^T (E_p (t) \phi ' (t) )' \int _0 ^1 x \psi _1 (x) \Bigl( F(\rho) - F(\xi) \Bigr) \, dx \, dt \Bigr\vert 
\leq \frac{3}{2} E_p (S) ^2  \phi ' (S) .
\end{equation}

Finally, we note that, for all $x \in (0, a_1) \cup (b_1,1)$ we have $\psi _1 (x)=1$, hence
$$ \forall x \in (0, a_1) \cup (b_1,1) , \quad 1 - (x\psi _1 (x))' = 1-1 = 0 .$$
Therefore, denoting $C_1 := \sup _{(0,1)} \vert 1 - (x\psi _1 (x))' \vert$, we have
\begin{multline*}
\int _S ^T E_p \phi ' \int _0 ^1 \Bigl( 1 - (x\psi _1 (x))' \Bigr) \Bigl( F(\rho) + F(\xi) \Bigr)
= \int _S ^T E_p \phi ' \int _{a_1} ^{b_1} \Bigl( 1 - (x\psi _1 (x))' \Bigr) \Bigl( F(\rho) + F(\xi) \Bigr)
\\
\leq C_1 \int _S ^T E_p \phi ' \int _{a_1} ^{b_1}  \Bigl( F(\rho) + F(\xi) \Bigr)
= \frac{C_1}{p} \int _S ^T E_p \phi ' \int _{a_0} ^{b_0}  \vert \rho \vert ^p + \vert \xi \vert ^p .
\end{multline*}
Since $s f(s) = \vert s \vert ^p$, and $\psi _2$ is nonnegative and equal to $1$ on $(a_1,b_1)$, we have
\begin{multline}
\label{estim1-3}
 \int _S ^T E_p \phi ' \int _0 ^1 \Bigl( 1 - (x\psi _1 (x))' \Bigr) \Bigl( F(\rho) + F(\xi) \Bigr)
\\
\leq \frac{C_1}{p} \int _S ^T E_p \phi ' \int _{0} ^1  \psi _2 (x) \Bigl( \rho f(\rho) + \xi f (\xi) \Bigr) .
\end{multline}

Then we deduce \eqref{mult-1st-mult-cons} from \eqref{mult-1st-mult} and \eqref{estim1-1}-\eqref{estim1-3}.\qed



\subsubsection{Proof of Lemma \ref{lem-2nd-estim}} \hfill

The starting point is \eqref{mult-2nd-mult}, obtained multiplying the first equation of \eqref{eq-invariants} by $E_p \phi ' \psi _2 (x) f'(\rho) z $ and the second equation of \eqref{eq-invariants} by $E_p \phi ' \,  \psi _2 (x) f'(\xi) z$,
as in \cite{KMJC2022} (the only difference being here the presence of the multiplier $\phi'$.)

Then we start by estimating the first four terms of the right-hand side of \eqref{mult-2nd-mult}.

{\underline{The first term of the right hand side of \eqref{mult-2nd-mult}:}
given $t\geq 0$, and using H\"older's inequality with 
$$ q = \frac{p}{p-1}, \quad \text{ such that } \quad \frac{1}{p}+\frac{1}{q} =1,$$
we have

\begin{multline*}
\Bigl\vert \int _0 ^1 \psi _2 (x) z \Bigl( f(\xi) - f(\rho)\Bigr) \, dx \Bigr\vert 
\leq \int _0 ^1 \vert z \vert \vert f(\xi) \vert + \vert z \vert \vert f(\rho) \vert \, dx
\\
\leq \Bigl( \int _0 ^1 \vert z \vert ^p \, dx \Bigr) ^{1/p} \Bigl( \int _0 ^1 (\vert \xi \vert ^{p-1} ) ^q \, dx \Bigr) ^{1/q}
+ \Bigl( \int _0 ^1 \vert z \vert ^p \, dx \Bigr) ^{1/p} \Bigl( \int _0 ^1 (\vert \rho \vert ^{p-1} ) ^q \, dx \Bigr) ^{1/q}
\\
\leq \frac{2}{p} \int _0 ^1 \vert z \vert ^p \, dx + \frac{1}{q} \int _0 ^1 (\vert \xi \vert ^{p-1} ) ^q + (\vert \rho \vert ^{p-1} ) ^q \, dx
\\
= \frac{2}{p} \int _0 ^1 \vert z \vert ^p \, dx + \frac{p-1}{p} \int _0 ^1 \vert \xi \vert ^p + \vert \rho \vert ^p \, dx
= \frac{2}{p} \int _0 ^1 \vert z \vert ^p \, dx  + (p-1) E_p (t) .
\end{multline*}
Moreover, using Dirichlet boundary condition and Poincar\'e's inequality, we have
$$ \int _0 ^1 \vert z \vert ^p \, dx \leq \int _0 ^1 \vert z_x \vert ^p \, dx
= \frac{1}{2^p} \int _0 ^1 \vert \rho + \xi \vert ^p \, dx \leq C \int _0 ^1 \vert \rho \vert ^p + \vert \xi \vert ^p \, dx = p C E_p (t),$$
hence there exists some (explicit) $C_p '$ such that
$$ \Bigl\vert \int _0 ^1 \psi _2 (x) z \Bigl( f(\xi) - f(\rho)\Bigr) \, dx \Bigr\vert 
\leq C_p ' E_p (t) .$$
Therefore, and once again using that $\phi '$ is nonincreasing, we have
\begin{equation}
\label{estim-2.1}
\Bigl\vert \Bigl[ E_p \phi ' \int _0 ^1 \psi _2 (x) z \Bigl( f(\xi) - f(\rho)\Bigr) \Bigr] _S ^T \Bigr\vert 
\leq  2 C_p ' E_p (S) ^2 \phi '(S) .
\end{equation}

{\underline{The second term of the right hand side of \eqref{mult-2nd-mult}:}

Using the previous estimates, we have
$$
\Bigl\vert \int _S ^T (E_p \phi ')' \int _0 ^1 \psi _2 (x) z \Bigl( f(\rho) - f(\xi)\Bigr) \Bigl\vert
\leq \int _S ^T \Bigl\vert E_p ' \phi ' + E_p \phi '' \Bigl\vert C_p ' E_p (t) \, dt
$$
Therefore, as in the previous subsection, we have
\begin{equation}
\label{estim-2.2}
\Bigl\vert \int _S ^T (E_p \phi ')' \int _0 ^1 \psi _2 (x) z \Bigl( f(\rho) - f(\xi)\Bigr) \Bigr\vert \leq 
\frac{3}{2} C_p ' E_p (S) ^2 \phi '(S) .
\end{equation}

{\underline{The third term of the right hand side of \eqref{mult-2nd-mult}:} once again we use H\"older estimates:
\begin{multline*}
\Bigl\vert - \int _S ^T E_p \phi ' \int _0 ^1 \psi _2 ' (x)  \, z \Bigl( f(\rho) + f(\xi)\Bigr) \Bigr\vert 
\leq C \int _S ^T E_p \phi ' \int _{a_2} ^{b_2} \vert z \vert \Bigl( \vert f(\rho)\vert + \vert f(\xi)\vert \Bigr) 
\\
\leq C \int _S ^T E_p \phi ' \Bigl[ \Bigl( \int _{a_2} ^{b_2} \vert z \vert ^p \, dx \Bigr) ^{1/p} \Bigl( \int _{a _2} ^{b_2} (\vert \xi \vert ^{p-1} ) ^q \, dx \Bigr) ^{1/q}
\\
+ \Bigl( \int _{a_2} ^{b_2} \vert z \vert ^p \, dx \Bigr) ^{1/p} \Bigl( \int _{a_2} ^{b_2} (\vert \rho \vert ^{p-1} ) ^q \, dx \Bigr) ^{1/q} \Bigr] .
\end{multline*}
For $\eta >0$, we have
$$ ab = (\frac{a}{\eta}) (\eta \, b) \leq \frac{a ^p}{p \, \eta ^p } + \frac{\eta ^q \, b^q}{q}.$$
This gives
\begin{multline*}
\Bigl[ \Bigl( \int _{a_2} ^{b_2} \vert z \vert ^p \, dx \Bigr) ^{1/p} \Bigl( \int _{a_2} ^{b_2} (\vert \xi \vert ^{p-1} ) ^q \, dx \Bigr) ^{1/q}
+ \Bigl( \int _{a_2} ^{b_2} \vert z \vert ^p \, dx \Bigr) ^{1/p} \Bigl( \int _{a_2} ^{b_2} (\vert \rho \vert ^{p-1} ) ^q \, dx \Bigr) ^{1/q} \Bigr]
\\
\leq \frac{2}{p \, \eta ^p} \Bigl( \int _{a_2} ^{b_2} \vert z \vert ^p \, dx \Bigr)
+ \frac{\eta ^q}{q} \Bigl( \int _{a_2} ^{b_2} \vert \xi \vert ^p + \vert \rho \vert ^p \, dx \Bigr) 
\\ 
\leq \frac{2}{p \, \eta ^p} \Bigl( \int _{a_2} ^{b_2} \vert z \vert ^p \, dx \Bigr)
+ \eta ^q (p-1)  E_p (t) .
\end{multline*}
Therefore
\begin{multline}
\label{estim-2.3}
\Bigl\vert - \int _S ^T E_p \phi ' \int _0 ^1 \psi _2 ' (x)  \, z \Bigl( f(\rho) + f(\xi)\Bigr) \Bigr\vert 
\\
\leq C \frac{2}{p \, \eta ^p} \int _S ^T E_p \phi ' \Bigl( \int _{a_2} ^{b_2} \vert z \vert ^p \, dx \Bigr) \, dt
+ C \eta ^q (p-1) \int _S ^T E_p ^2 \phi ' \, dt .
\end{multline}

{\underline{The fourth term of the right hand side of \eqref{mult-2nd-mult}:}
first we note that since $g$ is globally Lipschitz and $g(0)=0$, there exists $C_g$ such that
$$ \forall s \in \Bbb R, \quad \vert g(s) \vert \leq C_g \vert s \vert, $$
and then we have
$$
\Bigl\vert  z  g(\frac{\rho - \xi}{2}) \Bigl( f'(\rho) + f'(\xi)\Bigr) \Bigr\vert 
\leq C_g  \vert z \vert \, \Bigl \vert \frac{\rho - \xi}{2} \Bigr \vert \Bigl \vert f'(\rho) + f'(\xi)\Bigr \vert .
$$
But, using H\"older estimates (and the fact that $p>2$), we have
\begin{multline*}
\vert \rho - \xi \vert \, \vert f'(\rho) + f'(\xi)\Bigr \vert
\leq (p-1) \Bigl( \vert \rho \vert + \vert \xi \vert \Bigr) \Bigl( \vert \rho \vert ^{p-2} + \vert \xi \vert ^{p-2} \Bigr)
\\
= (p-1) \Bigl( \vert \rho \vert ^{p-1} + \vert \xi \vert ^{p-1} + \vert \rho \vert \vert \xi \vert ^{p-2} + \vert \xi \vert \vert \rho \vert ^{p-2} \Bigr)
\\
\leq C_p \Bigl( \vert \rho \vert ^{p-1} + \vert \xi \vert ^{p-1} \Bigr)
= C_p \Bigl( \vert f(\rho) \vert + \vert f(\xi) \vert \Bigr).
\end{multline*}
Then, we see that
\begin{multline*}
\Bigl\vert \int _0 ^1 \psi _2 (x)   a(x)  \, z  g(\frac{\rho - \xi}{2}) \Bigl( f'(\rho) + f'(\xi)\Bigr) \, dx \Bigr\vert 
\\
\leq \frac{C_g C_p}{2} \int _0 ^1 \psi _2 (x)   a(x)  \, \vert z \vert \, \Bigl( \vert f(\rho) \vert + \vert f(\xi) \vert \Bigr) 
\\
\leq \frac{C_g C_p}{2} \Vert a \Vert _\infty \int _{a_2} ^{b_2} \vert z \vert \, \Bigl( \vert f(\rho) \vert + \vert f(\xi) \vert \Bigr) .
\end{multline*}
Then 
\begin{multline*}
\Bigl\vert - \int _S ^T E_p \phi ' \int _0 ^1 \psi _2 (x)   a(x)  \, z  g(\frac{\rho - \xi}{2}) \Bigl( f'(\rho) + f'(\xi)\Bigr) \Bigr\vert 
\\ 
\leq \frac{C_g C_p}{2} \Vert a \Vert _\infty \int _S ^T E_p \phi ' \int _{a_2} ^{b_2} \vert z \vert \, \Bigl( \vert f(\rho) \vert + \vert f(\xi) \vert \Bigr),
\end{multline*}
and we are in the same situation as in the third term, studied just previously. Therefore, denoting $C''= \frac{C_g C_p}{2} \Vert a \Vert _\infty$, and given $\eta >0$, we have
\begin{multline}
\label{estim-2.4}
\Bigl\vert - \int _S ^T E_p \phi ' \int _0 ^1 \psi _2 (x)   a(x)  \, z  g(\frac{\rho - \xi}{2}) \Bigl( f'(\rho) + f'(\xi)\Bigr) \Bigr\vert 
\\ 
\leq C'' \frac{2}{p \, \eta ^p} \int _S ^T E_p \phi ' \Bigl( \int _{a_2} ^{b_2} \vert z \vert ^p \, dx \Bigr) \, dt
+ C'' \eta ^q (p-1) \int _S ^T E_p ^2 \phi ' \, dt .
\end{multline}

{\underline{Consequence for \eqref{mult-2nd-mult}:} Equations \eqref{mult-2nd-mult} and \eqref{estim-2.1}-\eqref{estim-2.4}
give that there exists some $C_2$ such that, given $\eta >0$ we have
\begin{multline}
\label{mult-2nd-mult-cons1}
\int _S ^T E_p \phi ' \int _0 ^1 \psi _2 (x) \Bigl( \rho f(\rho) + \xi f(\xi) \Bigr)
\\
\leq C_2 E_p (S) ^2 \phi '(S) + C_2 \eta ^q \int _S ^T E_p ^2 \phi ' \, dt 
+ \frac{C_2}{\eta ^p} \int _S ^T E_p \phi ' \Bigl( \int _{a_2} ^{b_2} \vert z \vert ^p \, dx \Bigr) \, dt
\\
+ 2 \int _S ^T E_p \phi ' \int _0 ^1 \psi _2 (x) \, \frac{\rho - \xi}{2} \Bigl( f(\rho) - f(\xi)\Bigr) .
\end{multline}

{\underline{Consequence of \eqref{mult-1st-mult-cons} and \eqref{mult-2nd-mult-cons1}:}
gathering \eqref{mult-1st-mult-cons} with the previous equation yields that
\begin{multline*}
\int _S ^T E_p ^2 \phi ' 
\leq ( \frac{7}{2} + C_1 C_2) E_p (S) ^2  \phi ' (S) + C_1 C_2  \eta ^q \int _S ^T E_p ^2 \phi ' \, dt 
\\
+ \frac{C_1 C_2}{\eta ^p} \int _S ^T E_p \phi ' \Bigl( \int _{a_2} ^{b_2} \vert z \vert ^p \, dx \Bigr) \, dt
\\
+ 2 C_1 \int _S ^T E_p \phi ' \int _0 ^1 \psi _2 (x) \, \frac{\rho - \xi}{2} \Bigl( f(\rho) - f(\xi)\Bigr)
\\
- \int _S ^T E_p \phi ' \int _0 ^1 x a(x) \psi _1 (x) g(\frac{\rho - \xi}{2}) \Bigl( f(\rho) + f(\xi) \Bigr) .
\end{multline*}
Choosing $\eta$ such that $\frac{C_1 C_2}{\eta ^p} =\frac{1}{2}$, we obtain that there exists $C_2 '$ such that
\eqref{mult-2nd-mult-cons} holds. \qed



\subsubsection{Proof of Lemma \ref{lem-3rd-estim} } \hfill

The starting point is \eqref{mult-3rd-mult}, obtained by multiplying the first equation and the second equation of \eqref{eq-invariants} by $E_p \phi ' v $, as in \cite{KMJC2022} (the only difference being here the presence of the multiplier $\phi'$.)

Now, \eqref{mult-3rd-mult} allows us to estimate the localized term in $\vert z \vert ^p$ that appears in 
\eqref{mult-2nd-mult-cons}. This involves the function $v$ defined in \eqref{def-v}. Hence a preliminary step
is to have some estimates on $v$. We recall the following.

\begin{Lemma}(\cite{KMJC2022})
\label{lem-estim-v}
Consider $v$ the solution of \eqref{def-v}. Then there exist positive constants $C, C_p$ such that we have the following estimates:
\begin{equation}
\label{estim-v-Lq}
\forall t >0, \quad \int _0 ^1 \vert v(t,x)\vert  ^q \, dx \leq C C_p E_p (t) ,
\end{equation}
and, for all $\sigma >0$, 
\begin{equation}
\label{estim-v_t-Lq}
\forall t >0, \quad \int _0 ^1 \vert v_t (t,x)\vert  ^q \, dx \leq C_p (p-2)\sigma E_p (t) + \frac{C_p}{\sigma ^{p-2}} \int _0 ^1 \psi _3 (x) \, \vert z_t (t,x)\vert  ^p \, dx .
\end{equation}
\end{Lemma}

\noindent {\it Proof of Lemma \ref{lem-estim-v}.} The proof is based on the following remark:
$\vert f(z) \vert ^q = \vert z \vert ^{(p-1)q} = \vert z \vert ^p$, therefore the right hand side of \eqref{def-v} is in $L^q$.
One can easily solve \eqref{def-v} and find the expression of $v$:
$$ v(t,x)= -x \int _x ^1 (1-y) \psi _3 (y) f(z(t,y)) \, dy + (x-1) \int _0 ^x y \psi _3 (y) f(z(t,y)) \, dy ,$$
and \eqref{estim-v-Lq} derives directly from this integral expression, using H\"older's and Poincar\'e's inequalities.

Moreover,
$$ v_t (t,x)= -x \int _x ^1 (1-y) \psi _3 (y) f'(z(t,y)) z_t (t,y) \, dy + (x-1) \int _0 ^x y \psi _3 (y) f'(z(t,y)) z_t (t,y) \, dy .$$
Then we deduce that
\begin{multline*}
\vert v_t (t,x) \vert ^q 
\leq C \Bigl( \int _0 ^1 \psi _3 (y) \, \vert f'(z(t,y))\vert \, \vert z_t (t,y) \vert \, dy \Bigr) ^q 
\\
\leq C' \Bigl( \int _0 ^1 \psi _3 (y) \, \vert z(t,y)\vert ^{p-2} \, \vert z_t (t,y) \vert \, dy \Bigr) ^q
\\
\leq  C''  \int _0 ^1 \psi _3 (y) ^q \, \vert z(t,y)\vert ^{(p-2)q} \, \vert z_t (t,y) \vert ^q \, dy 
\\
= C''  \int _0 ^1 \Bigl( \sigma \vert z\vert ^p \Bigr) ^{(p-2)/(p-1)} \, \Bigl( \frac{\psi _3 \vert z_t \vert }{\sigma ^{(p-2)/p}} \Bigr) ^{p/(p-1)} \, dy 
\\
\leq C'' \int _0 ^1 \frac{p-2}{p-1} \Bigl( \sigma \vert z\vert ^p \Bigr)  + \frac{1}{p-1}  \Bigl( \frac{\psi _3 \vert z_t \vert }{\sigma ^{(p-2)/p}} \Bigr) ^{p} \, dx 
\\
= C'' \frac{p-2}{p-1}  \sigma \int _0 ^1 \vert z\vert ^p \, dy
+ \frac{C''}{(p-1) \sigma ^{p-2}} \int _0 ^1 \psi _3 ^p \vert z_t \vert ^p \, dy .
\end{multline*}
Using the $L^p$ Poincar\'e's inequality:
$$ \int _0 ^1 \vert z\vert ^p \, dy \leq C_p  \int _0 ^1 \vert z_x \vert ^p \, dy
= C_p  \int _0 ^1 \vert \rho + \xi \vert ^p \, dy
\leq C_p ' \int _0 ^1 \vert \rho \vert ^p + \vert \xi \vert ^p \, dx
= p C_p ' E_p (t) ,$$
which gives
$$ \vert v_t (t,x) \vert ^q 
\leq C'' \frac{p-2}{p-1}  \sigma p C_p ' E_p (t) + \frac{C''}{(p-1) \sigma ^{p-2}} \int _0 ^1 \psi _3 \vert z_t \vert ^p \, dy ,$$
and we conclude integrating w.r.t. $x$. \qed

{\underline{The first and second terms of the right hand side of \eqref{mult-3rd-mult}:}
we have, using \eqref{estim-v-Lq}:
$$ \Bigl \vert \int _0 ^1 v (\rho - \xi) \Bigr\vert
\leq \Bigl( \int _0 ^1 \vert v \vert ^q \Bigr) ^{1/q} \, \Bigl( \int _0 ^1 \vert \rho - \xi \vert ^p \Bigr) ^{1/p} 
\leq C E_p (t) ^{1/q} \, E_p (t) ^{1/p} = C E_p (t) .$$
Therefore
\begin{equation}
\label{estim-3.1}
\Bigl\vert \Bigl[ E_p \phi ' \int _0 ^1 v (\rho - \xi) \Bigr] _S ^T  \Bigr\vert \leq C' E_p (S)^2 \phi '(S) .
\end{equation}
and
\begin{equation}
\label{estim-3.2}
\Bigl\vert - \int _S ^T (E_p \phi ')' \int _0 ^1  v (\rho - \xi) \Bigr\vert \leq C'' E_p (S)^2 \phi '(S) .
\end{equation}

{\underline{The third term of the right hand side of \eqref{mult-3rd-mult}:} fix $\tilde \sigma >0$; then using \eqref{estim-v_t-Lq}, we have
\begin{multline*}
\Bigl \vert \int _0 ^1 v_t (\rho - \xi) \Bigr\vert 
\leq \Bigl( \int _0 ^1 \vert v_t \vert ^q \Bigr) ^{1/q} \, \Bigl( \int _0 ^1 \vert \rho - \xi \vert ^p \Bigr) ^{1/p} 
\\
\leq C \Bigl( C_p (p-2)\sigma E_p (t) + \frac{C_p}{\sigma ^{p-2}} \int _0 ^1 \psi _3 (x) \, \vert z_t (t,x)\vert  ^p \, dx \Bigr) ^{1/q} \, \Bigl( E_p (t) \Bigr) ^{1/p} 
\\
= C \Bigl( C_p (p-2)\frac{\sigma}{\tilde \sigma ^{q/p}} E_p (t) + \frac{C_p}{\sigma ^{p-2} \tilde \sigma ^{q/p}} \int _0 ^1 \psi _3 (x) \, \vert z_t (t,x)\vert  ^p \, dx \Bigr) ^{1/q} \, \Bigl( \tilde \sigma E_p (t) \Bigr) ^{1/p} 
\\
\leq C' \Bigl[ C_p (p-2)\frac{\sigma}{\tilde \sigma ^{q/p}} E_p (t) + \frac{C_p}{\sigma ^{p-2} \tilde \sigma ^{q/p}} \int _0 ^1 \psi _3 (x) \, \vert z_t (t,x)\vert  ^p \, dx + \tilde \sigma E_p (t) \Bigr]
\end{multline*}
Choosing $\sigma$ such that $C_p (p-2)\frac{\sigma}{\tilde \sigma ^{q/p}} = \tilde \sigma$, we obtain that
$$ \Bigl \vert \int _0 ^1 v_t (\rho - \xi) \Bigr\vert 
\leq C' \Bigl[ 2 \tilde \sigma E_p (t)  + \frac{C_p}{\sigma ^{p-2} \tilde \sigma ^{q/p}} \int _0 ^1 \psi _3 (x) \, \vert z_t (t,x)\vert  ^p \, dx  \Bigr] .$$
Therefore 
\begin{multline}
\label{estim-3.3}
\Bigl\vert -  \int _S ^T E_p \phi ' \int _0 ^1  v_t  (\rho - \xi) \Bigr\vert 
\leq 2 C' \tilde \sigma \int _S ^T E_p ^2 \phi ' 
\\
+ C' \frac{C_p}{\sigma ^{p-2} \tilde \sigma ^{q/p}}
\int _S ^T E_p \phi ' \int _0 ^1 \psi _3 (x) \, \vert z_t (t,x)\vert  ^p \, dx .
\end{multline}

{\underline{Consequence for \eqref{mult-3rd-mult}:} Equations \eqref{mult-2nd-mult} and \eqref{estim-3.1}-\eqref{estim-3.3}
give that

\begin{multline}
\label{mult-3rd-mult-cons1}
2 \int _S ^T E_p \phi ' \int _0 ^1 \psi _3 (x) \vert z \vert ^p 
\\
\leq (C'+C'') E_p (S)^2 \phi '(S)
+ 2 C' \tilde \sigma \int _S ^T E_p ^2 \phi ' 
\\
+ C' \frac{C_p}{\sigma ^{p-2} \tilde \sigma ^{q/p}}
\int _S ^T E_p \phi ' \int _0 ^1 \psi _3 (x) \, \vert z_t (t,x)\vert  ^p \, dx
\\ + 2 \int _S ^T E_p \phi ' \int _0 ^1 a (x) \, v g(\frac{\rho - \xi}{2}) .
\end{multline}

{\underline{Consequence of \eqref{mult-3rd-mult-cons1} for \eqref{mult-2nd-mult-cons}:} since $\tilde \psi _1 =1$
on $(a_2, b_2)$, \eqref{mult-2nd-mult-cons} allows us to estimate from above the second term of the right-hand side of \eqref{mult-3rd-mult-cons1} and we obtain that

\begin{multline*}
\int _S ^T E_p ^2 \phi ' 
\leq C_3 (1+ C'+C'') E_p (S) ^2  \phi ' (S) 
+ 2 C_3 C' \tilde \sigma \int _S ^T E_p ^2 \phi ' 
\\
+ C_3 C' \frac{C_p}{\sigma ^{p-2} \tilde \sigma ^{q/p}}
\int _S ^T E_p \phi ' \int _0 ^1 \psi _3 (x) \, \vert z_t (t,x)\vert  ^p \, dx
\\
+ 2C_3 \int _S ^T E_p \phi ' \int _0 ^1 a (x) \, v g(\frac{\rho - \xi}{2}) 
\\
+ 2 C_1 \int _S ^T E_p \phi ' \int _0 ^1 \psi _2 (x) \, \frac{\rho - \xi}{2} \Bigl( f(\rho) - f(\xi)\Bigr)
\\
- 2\int _S ^T E_p \phi ' \int _0 ^1 x a(x) \psi _1 (x) g(\frac{\rho - \xi}{2}) \Bigl( f(\rho) + f(\xi) \Bigr) .
\end{multline*}
Choosing $\tilde \sigma$ such that $2 C_3 C' \tilde \sigma= \frac{1}{2}$, and remembering that $2z_t= \rho - \xi$, we obtain that there exists $C_4$ such that \eqref{mult-3rd-mult-cons} holds. \qed


\section{Proof of Theorem  \ref{thm-stab2}}
\label{sect6}

Here we prove stability results when $1<p<2$, adapting the proof done in the case $p\geq 2$.
Due to the presence of $f'$ (and therefore of the power $p-2$) in the second set of  multipliers in the case $p\geq 2$ (see subsection \ref{sec-sec-mult}), it is not clear if one has the right to use them directly in the case $1<p<2$. However, we are able to adapt our arguments, mainly following and completing \cite[Section 4.2]{KMJC2022}.

\subsection{A modified energy} \hfill

The first thing is to modify the functions $f$ and $F$ just like it was done in \cite[Section 4.2]{KMJC2022}. 
Thereore, we consider, for $p\in (1,2)$, the functions $\tilde f$ and $\tilde F$ defined on $\mathbb{R}$, by
\begin{align}
\tilde f (y)&= (p-1)\int_0^y (|s| + 1)^{p-2} \, ds =\textrm{sgn}(y)\left[(|y|+1)^{p-1}-
1\right], \label{gdef}\\
\tilde F (y)&=\int_0^y \tilde f(s)\,ds =\frac{1}{p}\left[(|y|+1)^p-1\right]-|y|.
\label{Gdef}
\end{align}
This drives us to modify the energy $E_p$ by considering,  for every
$t\in \mathbb{R_+}$ and every solution of \eqref{eq-wave}, the function $\tilde E_p$ defined by
\begin{align}\label{Eepsdef}
\tilde E_p (t)= \int_0^1 \left( \tilde F (\rho) + \tilde F (\xi) \right) \, dx.
\end{align}

First we note that $\tilde E_p$ is nonincreasing: indeed, the function $\tilde f$ is continous on $\Bbb R$, increasing and odd, therefore the function $\tilde F$ is $C^1$ on $\Bbb R$, convex (and even), and then Proposition \ref{prop1} ensures that $\tilde E_p$ is nonincreasing.


\subsection{The key integral inequality on the modified energy} \hfill

We are going to prove the following extension of  Lemma \ref{lem-int-estim}:

\begin{Lemma}
\label{lem-int-estim*}
Let $p \in (1,2)$. Then there exists $C ^{(0)} >0$ (depending on the norm of the initial conditions) such that
\begin{multline}
\label{eq-int-estim*}
\int _S ^T \tilde E_p (t) ^2 \phi ' (t) \, dt
\leq C ^{(0)} \tilde E_p (S) ^2  \phi ' (S) 
\\ + C ^{(0)} \int _S ^T \tilde E_p (t) \phi ' (t) \int _0 ^1 a (x) \, \frac{\rho - \xi}{2} \Bigl( \tilde f(\rho) - \tilde f(\xi)\Bigr) \, dx \, dt .
\end{multline}
\end{Lemma}

Before proving Lemma \ref{lem-int-estim*}, let us prove the decay rate estimates of Theorem \ref{thm-stab2}.


\subsection{Decay rate estimates of $\tilde E_p$ and $E_p$ assuming Lemma \ref{lem-int-estim*}} \hfill

The key integral inequality \eqref{eq-int-estim*} will allow us to prove decay rate estimates of $\tilde E_p$ when $p\in (1,2)$. But $E_p$ is a more natural energy, and we can note that

$$ \Bigl( \int _{\vert \rho \vert \leq 1} \vert \rho \vert ^p \Bigr)^{2/p}
\leq \int _{\vert \rho \vert \leq 1} \vert \rho \vert ^2 \leq \tilde E_p (t) ,$$
and therefore
\begin{multline*}
E_p (t) \leq c \tilde E_p (t) ^{p/2} + c \tilde E_p (t) 
\\
\leq c' (1+ \tilde E_p (0) ^{\frac{2-p}{2}}) \tilde E_p (t) ^{p/2}
\leq c'' (1+  E_p (0) ^{\frac{2-p}{2}}) \tilde E_p (t) ^{p/2} ,
\end{multline*}
therefore a decay rate estimate on $\tilde E_p$ gives immediately a decay rate estimate for $E_p$.

\subsubsection{Theorem \ref{thm-stab2}, (i): $g$ linear} \hfill

If $g(s)=\alpha s$ with some $s>0$, then we choose $\phi (t) = t$, and \eqref{eq-int-estim*} gives us that

\begin{multline*}
\int _S ^T \tilde E_p (t) ^2 \, dt
\leq C ^{(0)} \tilde E_p (S) ^2 
\\+ \frac{C ^{(0)}}{2\alpha} \int _S ^T \tilde E_p (t) \int _0 ^1 a (x) \, g(\rho - \xi) \Bigl( \tilde f(\rho) - \tilde f(\xi)\Bigr) \, dx \, dt 
\\
= (1+\frac{1}{4\alpha}) C ^{(0)} \tilde E_p (S) ^2,
\end{multline*}
and we deduce that $\tilde E_p$ decays exponentially to $0$. Since the problem is linear, the decay is uniform. \qed

\subsubsection{Theorem \ref{thm-stab2}, (ii): $g'(0) >0$} \hfill

The situation is almost the same, the main difference being that the decay rate depends on
$C^{(0)}$, hence on $E_p (0)$. \qed

\subsubsection{Theorem \ref{thm-stab2}, (iii): $g'(0) =0$} \hfill

Here we proceed as in subsection \ref{sec-thm-ii}. the only difference lies in $\Omega _1 ^{(t)}$, where now$\tilde f$ is bounded. Therefore
$$ \tilde f (\rho) - \tilde f (\xi) \vert \leq 2 (p-1) \varepsilon (t) ,$$
and
$$ \int _S ^T \tilde E_p \phi ' \int _{\Omega _1 ^{(t)}} a (x) \, (\rho - \xi) \Bigl( \tilde f(\rho) - \tilde f(\xi)\Bigr) \leq C'' \int _S ^T \tilde E_p (t) \phi ' (t) \varepsilon (t) ^2 \, dt .$$
Choosing $\phi := \phi _m$ defined in \eqref{choice-phi}, and $\varepsilon := \varepsilon _m$ defined in \eqref{choice-eps}, we have
$$ \int _S ^T \tilde E_p (t) ^2 \phi ' (t) \, dt \leq C^{(0)} \tilde E_p (S) ^2 + C^{(0)} \frac{\tilde E_p (S)}{(\phi _m (S) +A)^{2m-1}}.$$
Then Lemma \ref{lem3-COCV} gives that
$$ \tilde E_p (t) \leq \frac{C_m(E_p (0))}{(\phi _m (t) +A) ^{2m}} .$$
Then, as in Theorem \ref{thm-stab}, we have
$$ \tilde E_p (t) \leq \frac{C_m(E_p (0))}{s(t) ^{2m}} \quad \text{ where }
t = \frac{2s(t)}{H(\frac{1}{(2s(t))^m})} .$$
Concerning the natural energy, we have
$$ E_p (t) \leq \frac{C_m(E_p (0))}{s(t) ^{pm}} \quad \text{ where }
t = \frac{2s(t)}{H(\frac{1}{(2s(t))^m})} .$$
If $g_0$ satisfies \eqref{assump-poly-0}, then we deduce from \eqref{growth-phi_n}
that
$$ E_p (t) \leq \frac{C_m(E_p (0))}{t ^{(mp)/[m(q-1)+1]}} .$$
Letting $m\to \infty$, we find a decay of the form
$$ E_p (t) \leq \frac{C_\delta (E_p (0))}{t ^{\frac{p}{q-1} -\delta }}, $$
hence as in the case $p>2$. \qed

\medskip

Therefore it only remains to prove Lemma \ref{lem-int-estim*}. Its proof follows from convex analysis arguments, combined with identities given by the multiplier method.
In the following, we start by explaining the tools that we will use later, in particular some extensions of H\"older's inequality that wil be adapted to our setting.


\subsection{Fenchel's inequality and applications} \hfill
\label{sec-Fenchel}

It is easy to see that the asymptotic behavior of $\tilde F$ as $y\to 0^+$ and $y\to +\infty$ is the following:
$$ \tilde F (y) \sim _{y\to 0} \frac{p-1}{2} y^2, \quad \text{ and } \quad \tilde F (y) \sim _{y\to +\infty} \frac{y^p}{p} .$$
Then one can consider the convex conjugate of $\tilde F$, denoted $\tilde F^*$ and defined
by
\begin{equation}
\label{def-F*}
\forall b \in \Bbb R, \quad  \tilde F ^* (b):= \sup _{y\in \Bbb R} \Bigl( by - \tilde F (y) \Bigr) .
\end{equation}
This function has the following properties, that we will use later:
\begin{itemize}
\item  $\tilde F^*$ takes finite values, thanks to the asymptotic behavior of $\tilde F$,
\item  $\tilde F^*$  is positive on $(0,+\infty)$, equal to $0$ in $0$ (since $\tilde F$ is nonnegative), 
\item $\tilde F^*$ is even, since $\tilde F$ is even, 
\item $\tilde F^*$ is convex, since the derivative of $\tilde F ^*$ is $\tilde f ^{-1}$, that is increasing.
\end{itemize}

We also recall the classical Fenchel's inequality (\cite{JBHU}):
\begin{equation}
\label{Fenchel}
\forall a,b \in \Bbb R, \quad \vert ab \vert \leq \tilde F (\vert a \vert) + \tilde F^* (\vert b \vert) ,
\end{equation}
that is an immediate consequence of the definition of $\tilde F ^*$ (taking $y=a$) but will have nontrivial consequences. (When $\tilde F (y) = \frac{\vert y \vert ^p}{p}$, then $\tilde F ^* (b) = \frac{\vert b \vert ^q}{q}$ where $q$ is the conjugate exponent of $p$, and \eqref{Fenchel}
is the classical H\"older inequality.)

In particular, modifying a little some arguments of \cite[Section 4.2]{KMJC2022}, we have the following

\begin{Lemma}
\label{lem-Fenchel-cons}
There exists $C_p>0$ such that
\begin{equation}
\label{Fenchel-cons<}
\forall a, y \in \Bbb R, \forall \eta \in (0,1), \quad 
\vert a \tilde f(y) \vert \leq \frac{C_p}{\eta ^2} \tilde F (a) + C_p \eta ^2 \tilde F (y) ,
\end{equation}
and 
\begin{equation}
\label{Fenchel-cons>}
\forall a, y \in \Bbb R, \forall \eta \in (0,1), \quad 
\vert a \tilde f(y) \vert \leq C_p \, \eta ^p \tilde F (a) + \frac{C_p}{ \eta ^q} \tilde F (y) .
\end{equation}
\end{Lemma}

\noindent Note that \eqref{Fenchel-cons<} and \eqref{Fenchel-cons>} are very close, 
and allow one to have a small coefficient in front of $\tilde F (y)$ (using \eqref{Fenchel-cons<}) or in front of $\tilde F (a)$ (using \eqref{Fenchel-cons>}).

\medskip

\noindent {\it Proof of Lemma \ref{lem-Fenchel-cons}.} 
First, we prove \eqref{Fenchel-cons<}. Using \eqref{Fenchel}, we have
$$ \vert a \tilde f(y) \vert
= \vert (\frac{a}{\eta}) \, ( \eta \tilde f(y)) \vert
\leq \tilde F (\vert \frac{a}{\eta} \vert) + \tilde F^* (\vert \eta \tilde f(y) \vert).$$
Next, thanks to the asymptotic behavior of $\tilde F$ as $y\to 0^+$ and $y\to +\infty$:
$$ \tilde F (y) \sim _{y\to 0} \frac{p-1}{2} y^2, \quad \tilde F (y) \sim _{y\to +\infty} \frac{y^p}{p},$$
we see that there exists $C_p > C_p ' >0$ such that
\begin{equation}
\label{comp-tildeF}
\forall y \in \Bbb R, \quad C_p ' \min (y^2, y^p) \leq \tilde F (y) \leq C_p \min (y^2, y^p) .
\end{equation}
This yields 
\begin{multline}
\label{comp-tildeF-csq}
\tilde F (\vert \frac{a}{\eta} \vert) \leq C_p \min (\frac{a^2}{\eta ^2}, \frac{a^p}{\eta ^p}) \leq C_p \min (\frac{a^2}{\eta ^2}, \frac{a^p}{\eta ^2})
\\
= \frac{C_p}{\eta ^2} \min (a^2, a^p) 
\leq \frac{C_p }{\eta ^2 C_p '} \tilde F (\vert a \vert ) .
\end{multline}

On the other hand, looking at the definition of $\tilde F^*$, it is clear that, given $b\geq 0$, the function
$$ y \in \Bbb R \mapsto by - \tilde F (y) $$
attains its maximum at $y:= \tilde f ^{-1} (b)$ (remember that $\tilde f$ is increasing), and therefore we always have
$$ \forall b \geq 0, \quad \tilde F ^* (b) = b \, \tilde f ^{-1} (b) - \tilde F (\tilde f ^{-1} (b)) ,$$
or, equivalently,
\begin{equation}
\label{Fenchel-cons2}
\forall y \in \Bbb R, \quad y\, \tilde f(y) = \tilde F(y) + \tilde F ^* ( \tilde f(y)) .
\end{equation}
This allows us to obtain the asymptotic behavior of $\tilde F ^*(b)$ as $b\to 0$ and $b\to +\infty$: first, we have
$$ \tilde F ^* ( \tilde f(y)) \sim _{y\to 0} (p-1)y ^2 - \frac{p-1}{2} y^2 = \frac{p-1}{2} y^2 ,$$
therefore 
$$ \tilde F ^* (b) \sim \frac{b^2}{2(p-1)} \quad \text{ as } b \to 0;$$
next, we have
$$ \tilde F ^* ( \tilde f(y)) \sim _{y\to +\infty} y^p - \frac{y^p}{p} = \frac{p-1}{p} y^p,$$
therefore 
$$ \tilde F ^* (b) \sim \frac{p-1}{p} b ^{\frac{p}{p-1}} \quad \text{ as } b \to +\infty .$$
Since $\frac{p}{p-1} \geq 2$ and since $\tilde F ^*$ is increasing on $\Bbb R_+$ (and even on $\Bbb R$), we see that there exists $C ^* _p > {C_p ^*} ' >0$ such that
\begin{equation}
\label{comp-tildeF*}
\forall b \in \Bbb R, \quad {C_p ^*} ' \max (b^2, b ^{\frac{p}{p-1}} ) \leq \tilde F ^* (b) \leq C_p ^* \max (b^2, b ^{\frac{p}{p-1}} ).
\end{equation}
Now we can conclude:
\begin{multline*}
\tilde F^* (\vert \eta \tilde f(y) \vert) 
\leq C_p ^* \, \max ((\vert \eta \tilde f(y) \vert) ^2, (\vert \eta \tilde f(y) \vert) ^{\frac{p}{p-1}} )
\\
= C_p ^* \, \max ( \eta ^2 \tilde f(\vert y \vert) ^2,  \eta ^{\frac{p}{p-1}} \tilde f(\vert y \vert) ^{\frac{p}{p-1}} )
\leq C_p ^* \, \max ( \eta ^2 \tilde f(\vert y \vert) ^2,  \eta ^{2} \tilde f(\vert y \vert) ^{\frac{p}{p-1}} )
\\
= C_p ^* \eta ^2 \, \max (  \tilde f(\vert y \vert) ^2,  \tilde f(\vert y \vert) ^{\frac{p}{p-1}} ) .
\end{multline*}
But
$$ \max (  \tilde f(\vert y \vert) ^2,  \tilde f(\vert y \vert) ^{\frac{p}{p-1}} ) 
=  \tilde f(\vert y \vert) ^2 \sim (p-1) ^2 y^2 \quad \text{ as } y \to 0 ,$$
and
$$ \max (  \tilde f(\vert y \vert) ^2,  \tilde f(\vert y \vert) ^{\frac{p}{p-1}} ) 
=  \tilde f(\vert y \vert) ^{\frac{p}{p-1}} \sim y^p \quad \text{ as } y \to +\infty ,$$
therefore there exists $\tilde C_p > \tilde C_p ' >0$ such that
\begin{equation}
\label{comp-tildeFf*}
 \forall y \geq 0, \quad 
\tilde C_p ' \min (y^2, y^p)
\leq \max (  \tilde f(\vert y \vert) ^2,  \tilde f(\vert y \vert) ^{\frac{p}{p-1}} )
\leq \tilde C_p \min (y^2, y^p) .
\end{equation}
Then, using \eqref{comp-tildeF}, we obtain that
\begin{equation}
\label{comp-tildeFf*-csq}
\tilde F^* (\vert \eta \tilde f(y) \vert) 
\leq C_p ^* \eta ^2 \, \tilde C_p \min (y^2, y^p) \leq \frac{C_p ^* \tilde C_p}{C_p '} \eta ^2 \tilde F ( \vert y \vert ) .
\end{equation}
The proof of \eqref{Fenchel-cons<} follows from \eqref{comp-tildeF-csq}
and \eqref{comp-tildeFf*-csq}. \qed

\medskip

In the same spirit, we prove \eqref{Fenchel-cons>}: first we can assume that $a,y >0$, and using \eqref{Fenchel}, we can write
$$  a \tilde f(y) = (\eta a) \, \Bigl( \frac{\tilde f(y)}{\eta} \Bigr) 
\leq \tilde F (\eta a) + \tilde F ^* \Bigl( \frac{\tilde f(y)}{\eta} \Bigr) .$$
Using \eqref{comp-tildeF}, we have
$$ \tilde F (\eta a) \leq C_p \min (\eta ^2 a^2, \eta ^p a^p)
\leq C_p \eta  ^p \min ( a^2,  a^p) \leq \frac{C_p}{C_p '} \eta ^p \tilde F (a) .$$
Using successively \eqref{comp-tildeF*}, \eqref{comp-tildeFf*} and \eqref{comp-tildeF}, we have

\begin{multline*}
 \tilde F ^* \Bigl( \frac{\tilde f(y)}{\eta} \Bigr)
\leq C_p ^* \max (\frac{\tilde f(y)^2}{\eta ^2}, \frac{\tilde f(y) ^{\frac{p}{p-1}}}{\eta ^{\frac{p}{p-1}}}  )
\\
\leq \frac{C_p ^*}{\eta ^{\frac{p}{p-1}}} \max (\tilde f(y)^2, \tilde f(y) ^{\frac{p}{p-1}} )
\\
\leq \frac{C_p ^* \tilde C_p}{\eta ^{\frac{p}{p-1}}} \min (y^2, y^p)
\leq \frac{C_p ^* \tilde C_p}{\eta ^{\frac{p}{p-1}} C_p ' } \tilde F (y) .
\end{multline*}
This concludes the proof of \eqref{Fenchel-cons>} and of Lemma \ref{lem-Fenchel-cons}. \qed


As an application, we have the following generalized Poincar\'e's inequality:

\begin{Lemma} \cite[Section 4.2]{K2022}
\label{lem-Hardy-Fenchel}
Let $p\in (1,2)$. Then there exists a positive constant $C_p$ such that, for every absolutely continuous function $z: [0,1] \to \Bbb R$ satisfying $z(0)=0$, one has
\begin{equation}
\label{eq-Hardy-Fenchel}
\int _0 ^1 \tilde F (z(s)) \, ds \leq C_p \int _0 ^1 \tilde F (z'(s)) \, ds .
\end{equation}
\end{Lemma}

\noindent {\it Proof of Lemma \ref{lem-Hardy-Fenchel}.} We have
$$ \tilde F (z(x)) = \tilde F (z(x)) - \tilde F (z(0))
= \int _0 ^x \tilde f (z(s)) \, z'(s) \, ds .$$
Then using \eqref{Fenchel-cons<}, we have
$$ \Bigl \vert \tilde f (z(s)) \, z'(s)  \Bigr \vert 
\leq \frac{C_p}{\eta ^2} \tilde F (z'(s)) + C_p \eta ^2 \tilde F (z(s)) .$$
Therefore
$$ \tilde F (z(x)) \leq \int _0 ^1 \Bigl \vert \tilde f (z(s)) \, z'(s)  \Bigr \vert  \, ds 
\leq \frac{C_p}{\eta ^2} \int _0 ^1 \tilde F (z'(s)) \, ds 
+ C_p \eta ^2 \int _0 ^1 \tilde F (z(s)) \, ds .$$
Integrating with respect to $x$ on $(0,1)$, we get
$$ \int _0 ^1 \tilde F (z(x)) \, dx \leq 
\frac{C_p}{\eta ^2} \int _0 ^1 \tilde F (z'(s)) \, ds 
+ C_p \eta ^2 \int _0 ^1 \tilde F (z(s)) \, ds ,$$
and we obtain \eqref{eq-Hardy-Fenchel} choising $\eta$ small enough. \qed


\subsection{Proof of Lemma \ref{lem-int-estim*}} \hfill

The proof of Lemma \ref{lem-int-estim*} follows adapting the arguments given by \cite[Section 4.2.4]{KMJC2022},
using the following modified multipliers: 
\begin{itemize}
	\item modified first multiplier: we replace $f$ by $\tilde f$, and $E_p$ by $\tilde E_p$;
	\item modified second multiplier: we replace $f'$ by ${\tilde f} '$, and $E_p$ by $\tilde E_p$;
	\item modified third multiplier: we replace $f$ by $\tilde f $ in the elliptic problem \eqref{def-v}, and $v$ and $E_p$ by $\tilde v$ (the solution of the modified elliptic problem, now associated to $\tilde f $) and $\tilde E_p$.
\end{itemize}


\subsubsection{Estimate given by the modified first multiplier} \hfill


When $p\in (1,2)$, and with the modified first multiplier, we have the following estimate, that is a modified version of \eqref{mult-1st-mult-cons}:

\begin{Lemma}
\label{lem-1st-estim*}
There exists $C_1 >0$ independent of the initial conditions such that
\begin{multline}
\label{mult-1st-mult-cons*}
\int _S ^T \tilde E_p (t) ^2 \phi ' (t) \, dt
\leq C_1 \tilde E_p (S) ^2  \phi ' (S)
\\
+ C_1 \int _S ^T \tilde E_p(t)  \phi ' (t) \int _{0} ^1  \psi _2 (x) \Bigl( \rho \tilde f(\rho) + \xi \tilde f (\xi) \Bigr) \, dx \, dt
\\
- \int _S ^T \tilde E_p (t) \phi ' (t) \int _0 ^1 x a(x) \psi _1 (x) g(\frac{\rho - \xi}{2}) \Bigl( \tilde f(\rho) + \tilde f(\xi) \Bigr) \, dx \, dt.
\end{multline}
\end{Lemma}


\subsubsection{Estimate given by the modified second multiplier} \hfill

When $p\in (1,2)$, and with the modified second multiplier, we have the following estimate, which is a modified version of \eqref{mult-2nd-mult-cons}:

\begin{Lemma}
\label{lem-2nd-estim*}
Denote $k^{(0)} := \max (1, E_p (0) ^{\frac{2-p}{p}})$.
Then there exists $C_2 ' >0$ independent of the initial conditions such that
\begin{multline}
\label{mult-2nd-mult-cons*}
\int _S ^T \tilde E_p (t) ^2 \phi ' (t) \, dt
\leq \tilde C_2 ' \, k^{(0)} \tilde E_p (S) ^2  \phi ' (S) 
\\
+ \tilde C_2 ' \, { k ^{(0)} } ^{q} \int _S ^T \tilde E_p (t) \phi '(t)  \Bigl( \int _{a_2} ^{b_2}  \tilde F (z) \, dx \Bigr) \, dt
\\
+ 2 C_1 \int _S ^T \tilde E_p (t) \phi ' (t) \int _0 ^1 \psi _2 (x) \, \frac{\rho - \xi}{2} \Bigl( \tilde f(\rho) - \tilde f(\xi)\Bigr) \, dx \, dt
\\
- 2\int _S ^T \tilde E_p (t) \phi ' (t) \int _0 ^1 x a(x) \psi _1 (x) g(\frac{\rho - \xi}{2}) \Bigl( \tilde f(\rho) + \tilde f(\xi) \Bigr) \, dx \, dt .
\end{multline}
\end{Lemma}
(The only difference with \eqref{mult-2nd-mult-cons} is the fact that the coefficients in \eqref{mult-2nd-mult-cons*} depend on $E_p (0)$, hence on the norm of the initial condition.)


\subsubsection{Estimate given by the modified third multiplier} \hfill

Consider the solution $\tilde v$ of
\begin{equation}
\label{def-vtilde}
\begin{cases}
\tilde v_{xx} = \psi _3 (x) \tilde f(z(t,x)) \quad & x\in (0,1) ,\\
\tilde v(0)=0= \tilde v(1) .
\end{cases}
\end{equation}

When $p\in (1,2)$, and with the modified third multiplier, we have the following estimate, which is a modified version of \eqref{mult-3rd-mult-cons}:

\begin{Lemma}
\label{lem-3rd-estim*}
There exist exists $C_1 ^{(0)}, C_2 ^{(0)}, C_3 ^{(0)}$ that depend on $E_p (0)$ such that

\begin{multline}
\label{mult-3rd-mult-cons*}
\int _S ^T \tilde E_p ^2 \phi ' 
\leq C_1 ^{(0)} \tilde E_p (S) ^2  \phi ' (S) 
+ C_2 ^{(0)} \int _S ^T \tilde E_p \phi ' \int _0 ^1 \psi _3 (x) \tilde F (\rho - \xi ) \, dx
\\
+ C_3 ^{(0)} \int _S ^T \tilde E_p \phi ' \int _0 ^1 a (x) \, \tilde v g(\frac{\rho - \xi}{2}) 
+ 4 C_1 \int _S ^T \tilde E_p \phi ' \int _0 ^1 \psi _2 (x) \, \frac{\rho - \xi}{2} \Bigl( f(\rho) - f(\xi)\Bigr)
\\
- 4\int _S ^T \tilde E_p \phi ' \int _0 ^1 x a(x) \psi _1 (x) g(\frac{\rho - \xi}{2}) \Bigl( \tilde f(\rho) + \tilde f(\xi) \Bigr) .
\end{multline}
\end{Lemma}

(Note that, if it has some utility, we have
$$ C_1 ^{(0)} = \tilde C_2 ' \max (1, E_p (0) ^{\frac{2-p}{p}} ) + \max (1, E_p (0) ^{2\frac{2-p}{p-1}} ),$$
$$ C_2 ^{(0)} = \tilde C_2 ' \max (1, E_p (0) ^{(q ^2+q)\frac{2-p}{p}} ) \max (E_p (0) ^{\frac{2-p}{p}},E_p (0) ^{\frac{2-p}{p-1}}) , $$
$$  C_3 ^{(0)} = \tilde C_3 ' \max (1, E_p (0) ^{\frac{2-p}{p-1}} )).$$

We will need to extend \eqref{estim-v-Lq} and \eqref{estim-v_t-Lq}, and here are 
the natural extensions:
\begin{Lemma}
\label{lem-mai2}
There exists $C$ independent of the initial conditions such that
\begin{equation}
\label{estim-v-Lq*}
\forall t \geq 0, \quad \int _0 ^1 \tilde F ^* (\tilde v(t,x)) \, dx \leq C \max (1, E_p (0) ^{\frac{2-p}{p-1}} ) \, \tilde E_p (t) ,
\end{equation}
and
\begin{equation}
\label{estim-v_t-Lq*}
\forall t \geq 0, \quad \int _0 ^1 \tilde F ( \tilde v_t) \, dx \leq C \max (E_p (0) ^{\frac{2-p}{p}},E_p (0) ^{\frac{2-p}{p-1}}) \int _0 ^1  \psi _3  \tilde F ( z_t ) \, dx  .
\end{equation}
\end{Lemma}


\subsubsection{Proof of Lemma \ref{lem-int-estim*} assuming Lemmas \ref{lem-1st-estim*}-\ref{lem-mai2}} \hfill

We proceed as for the proof of Lemma \ref{lem-int-estim}, the starting point being now \eqref{mult-3rd-mult-cons*}. First, as we did in the previous section, using \eqref{Fenchel}, \eqref{comp-tildeF}, \eqref{comp-tildeF*}, and \eqref{estim-v-Lq*}, there exists $C>0$ such that, for all $\sigma \in (0,1)$, we have
\begin{multline*}
\int _S ^T \tilde E_p \phi ' \int _0 ^1 a (x) \, \tilde v g(\frac{\rho - \xi}{2}) 
\\
\leq C \sigma ^2 \int _S ^T \tilde E_p \phi ' \int _0 ^1 \tilde F ^* (v) 
+ \frac{C}{\sigma ^2} \int _S ^T \tilde E_p \phi ' \int _0 ^1 a (x) \tilde F (\rho - \xi)
\\
\leq C \sigma ^2 \max (1, E_p (0) ^{\frac{2-p}{p-1}} ) \int _S ^T \tilde E_p ^2 \phi '
+ \frac{C}{\sigma ^2} \int _S ^T \tilde E_p \phi ' \int _0 ^1 a (x) \tilde F (\rho - \xi) ,
\end{multline*}
and that allows us to get rid of the third term of the right hand side of \eqref{mult-3rd-mult-cons*} (taking $\sigma >0$ small enough).

Next, as we did in the proof of Lemma \ref{lem-int-estim}, we have
\begin{multline*}
- 4\int _S ^T \tilde E_p \phi ' \int _0 ^1 x a(x) \psi _1 (x) g(\frac{\rho - \xi}{2}) \Bigl( \tilde f(\rho) + \tilde f(\xi) \Bigr)
\\
= 4\int _S ^T \tilde E_p \phi ' \int _0 ^1 x a(x) \psi _1 (x) g(\frac{\rho - \xi}{2}) \Bigl( \tilde f(\rho) - \tilde f(\xi) \Bigr)
\\
+ 4\int _S ^T \tilde E_p \phi ' \int _0 ^1 x a(x) \psi _1 (x) g(\frac{\rho - \xi}{2}) \Bigl( -2 \tilde f(\rho) \Bigr)
\\
\leq 2 \tilde E_p (S) ^2 \phi '(S) 
+ 8 C_g \int _S ^T \tilde E_p \phi ' \int _0 ^1 a(x) \frac{\rho - \xi}{2} \tilde f(\rho) ,
\end{multline*}
and using \eqref{Fenchel-cons<}, we have

\begin{multline*}
\int _S ^T \tilde E_p \phi ' \int _0 ^1 a(x) \frac{\rho - \xi}{2} \tilde f(\rho)
\leq \int _S ^T \tilde E_p \phi ' \int _0 ^1 \frac{C_p}{\eta ^2} \tilde F (a(x) (\rho - \xi) ) + C_p \eta ^2 \tilde F (\rho)
\\
\leq \frac{C_p '}{\eta ^2} \int _S ^T \tilde E_p \phi ' \int _0 ^1 a(x) \tilde F (\rho - \xi ) + C_p \eta ^2 \int _S ^T \tilde E_p ^2 \phi ' .
\end{multline*}
Combining with the previous estimate and taking $\eta >0$ small enough, we can get rid of the fifth term of the righ hand side of \eqref{mult-3rd-mult-cons*}, and we obtain

\begin{multline}
\label{mult-3rd-mult-cons*2}
\int _S ^T \tilde E_p ^2 \phi ' 
\leq C_4 ^{(0)} \tilde E_p (S) ^2  \phi ' (S) 
+ C_5 ^{(0)} \int _S ^T \tilde E_p \phi ' \int _0 ^1 a (x) \tilde F (\rho - \xi ) \, dx
\\
+ 2 C_1 \int _S ^T \tilde E_p \phi ' \int _0 ^1 a (x) \, (\rho - \xi) \Bigl( \tilde f(\rho) - \tilde f(\xi)\Bigr) ,
\end{multline}
with $C_4 ^{(0)}$ and $C_5 ^{(0)}$ for which we can have an explicit estimate in terms of $E_p (0)$. It remains to estimate the second term of the right hand side of \eqref{mult-3rd-mult-cons*2}, and this follows from the following result, extracted from \cite[Lemma 3.2.7]{K2022} but for which we give an alternative proof.

\begin{Lemma}
\label{lem-mai7}
There exists $C_p$ uniform such that, for all $M\geq 1$, we have 

\begin{multline}
\label{eq-mai7}
\forall \rho, \xi \in \Bbb R, \quad
\tilde F (\rho - \xi) \leq \frac{C_p}{M^p} \Bigl( \tilde F (\rho) + \tilde F(\xi) \Bigr) 
\\
+ C_p M^{2-p} (\rho - \xi) (\tilde f(\rho) - \tilde f (\xi) ) .
\end{multline}

\end{Lemma}

\noindent {\it End of the proof of Lemma \ref{lem-int-estim*} assuming Lemma \ref{lem-mai7}.}  It follows from \eqref{mult-3rd-mult-cons*2} and \eqref{eq-mai7} that
\begin{multline}
\label{mult-3rd-mult-cons*3}
\int _S ^T \tilde E_p ^2 \phi ' 
\leq C_4 ^{(0)} \tilde E_p (S) ^2  \phi ' (S) 
+ C_5 ^{(0)} \frac{C_p}{M^p} \int _S ^T \tilde E_p ^2 \phi ' 
\\
+ ( C_5 ^{(0)} C_p M^{2-p} + 2 C_1 ) \int _S ^T \tilde E_p \phi ' \int _0 ^1 a (x) \, (\rho - \xi) \Bigl( \tilde f(\rho) - \tilde f(\xi)\Bigr) .
\end{multline}
Then it is sufficient to take $M$ large enough, and we obtain \eqref{eq-int-estim*}.  \qed

\medskip 
\noindent {\it Proof of Lemma \ref{lem-mai7}.}
Take $M \geq 1$, and divide the first quadrant of the plane $\{(\rho,\xi), \rho, \xi \in \Bbb R \}$ into the following four regions:
$$\mathcal R_1 := \{(\rho,\xi), \rho \geq 0, \xi \geq \frac{M+1}{M} \rho \},
\quad \mathcal R_2 := \{(\rho,\xi), \rho \geq 0, \rho \leq \xi \leq \frac{M+1}{M} \rho \} ,$$
$$\mathcal R_3 := \{(\rho,\xi), \xi \geq 0, \frac{M}{M+1} \rho \leq \xi \leq \rho \},
\quad \mathcal R_4 := \{(\rho,\xi), \xi \geq 0, \xi \leq \frac{M}{M+1} \rho \} .$$
Then, using \eqref{comp-tildeF},
$$ (\rho, \xi) \in \mathcal R_2 \quad \implies
\quad 0 \leq \xi - \rho \leq \frac{\rho}{M} \quad \implies
\tilde F (\xi - \rho) \leq \tilde F (\frac{\rho}{M}) \leq \frac{\tilde F(\rho)}{C_p ' M^p}.$$
In a symmetric way, 
$$ (\rho, \xi) \in \mathcal R_3 \quad \implies
\quad 0 \leq \rho - \xi \leq \frac{\xi}{M} \quad \implies
\tilde F (\rho - \xi) \leq \tilde F (\frac{\xi}{M}) \leq \frac{\tilde F(\xi)}{C_p ' M^p}.$$

Now, 
we note that
$$ \mathcal R_1 = \mathcal R_{1,-} \cup \mathcal R_{1,+},$$
where
$$ \mathcal R_{1,-} = \{(\rho, \xi) \in \mathcal R_1, \xi-\rho \leq 1 \}, \quad
\mathcal R_{1,+} = \{(\rho, \xi) \in \mathcal R_1, \xi-\rho \geq 1 \}.$$
It is clear that $\mathcal R_{1,-}$ is compact:
$$ (\rho, \xi) \in \mathcal R_{1,-} \quad \implies \quad \frac{\rho}{M} \leq \xi - \rho \leq 1 \quad \implies \quad \begin{cases} 0 \leq \rho \leq M, \\ 0 \leq \xi \leq 1+ \rho \leq 1+M . \end{cases} $$
Therefore, if $(\rho, \xi) \in \mathcal R_{1,-}$, we see from \eqref{comp-tildeF} that
$$ \tilde F (\xi - \rho) \leq C_p (\xi - \rho) ^2 ,$$
and, on the other hand, there exists $c_{\rho, \xi} \in (\rho, \xi) \subset [0,M+1]$ such that
$$ (\xi - \rho) (\tilde f(\xi) - \tilde f (\rho) )
= (\xi - \rho)^2 \frac{\tilde f(\xi) - \tilde f (\rho)}{\xi - \rho}
= (\xi - \rho)^2 \tilde f '(c_{\rho, \xi}) \geq (\xi - \rho)^2 \tilde f ' (M+1) .$$
Therefore
$$ (\rho, \xi) \in \mathcal R_{1,-} \quad \implies \quad 
\tilde F (\xi - \rho) \leq \frac{C_p}{\tilde f ' (M+1)} (\xi - \rho) (\tilde f(\xi) - \tilde f (\rho) ) .$$
Since $\tilde f'(x)= (p-1)(1+x)^{p-2}$, there exists $c_p >0$ independent of $M\geq 1$ such that
$$ (\rho, \xi) \in \mathcal R_{1,-} \quad \implies \quad 
\tilde F (\xi - \rho) \leq c_p M^{2-p} (\xi - \rho) (\tilde f(\xi) - \tilde f (\rho) ) .$$

Now, if $(\rho, \xi) \in \mathcal R_{1,+}$, introducing $\tau := \xi - \rho -1$, we see from \eqref{comp-tildeF} that 
$$ \tilde F (\xi - \rho) \leq C_p (\xi - \rho) ^p , \quad \text{ hence } \quad \frac{\tilde F (\xi - \rho)}{\xi - \rho} \leq C_p (\xi - \rho) ^{p-1} = C_p (1+\tau) ^{p-1},$$
and, on the other hand, we have
\begin{multline*}
\tilde f(\xi) - \tilde f (\rho) = (1+ \xi) ^{p-1} - (1+ \rho) ^{p-1}
= \Bigl( [1+ \tau ] + [1+\rho] \Bigr)^{p-1} - (1+ \rho) ^{p-1} 
\\
= (1+\tau) ^{p-1} \Bigl( (1 + \sigma )^{p-1} - \sigma ^{p-1} \Bigr) \quad \text{ where } 
\sigma = \frac{1+ \rho}{1+\tau} .
\end{multline*}
To conclude, it is sufficient to note that, first, the function 
$$  \sigma \in \Bbb R_+ \mapsto (1 + \sigma )^{p-1} - \sigma ^{p-1} $$
is nonincreasing, because $p<2$, and that $\sigma \leq 2 $ if $\rho \leq 1$, and 
$$  \sigma = \frac{1+ \rho}{1+\tau} = \frac{1+ \rho}{\xi-\rho}
\leq \frac{1+ \rho}{\frac{\rho}{M}} = M \frac{1+ \rho}{\rho} \leq 2M \quad \text{ if } \rho \geq 1 ,$$
hence, in any case, $\sigma \leq 2M$, since $M\geq 1$. Therefore, we obtain that,
if $(\rho, \xi) \in \mathcal R_{1,+}$, we have
$$ \tilde f(\xi) - \tilde f (\rho) \geq \Bigl( (1 + 2M )^{p-1} - (2M) ^{p-1} \Bigr) (1+\tau) ^{p-1} ,$$
which gives that
$$ (\rho, \xi) \in \mathcal R_{1,+} \quad \implies \quad 
\tilde F (\xi - \rho) \leq \frac{C_p}{(1 + 2M )^{p-1} - (2M) ^{p-1}} (\xi - \rho) (\tilde f(\xi) - \tilde f (\rho) ) .$$
Since we have
$$ (1 + 2M )^{p-1} - (2M) ^{p-1} \sim \frac{p-1}{(2M)^{2-p}} \quad \text{ as } p\to \infty ,$$
there exists some $c_p >0$ independent of $M \geq 1$ such that
$$ (\rho, \xi) \in \mathcal R_{1} \quad \implies \quad 
\tilde F (\xi - \rho) \leq c_p M^{2-p} (\xi - \rho) (\tilde f(\xi) - \tilde f (\rho) ) .$$
The situation is similar in $\mathcal R_4$.

If $\rho \xi \leq 0$, the proof is much easier:
assume that $\rho \leq 0 \leq \xi$. Then, either $\xi-\rho \leq 1$, or $\xi-\rho \geq 1$.
The subregion where $\xi-\rho \leq 1$ is exactly the triangle bounded by the horizontal axis $(O\rho)$, the vertical axis $(O\xi)$,  and the line $\xi = \rho +1$. Then $\vert \rho \vert, \xi \leq 1$, and we have
$$ \tilde F(\xi - \rho) \leq C_p (\xi-\rho) ^2 = C_p \Bigl( \xi + \vert \rho \vert \Bigr)^2,$$
and on the other hand,
\begin{multline*}
(\xi - \rho) (\tilde f(\xi) - \tilde f (\rho) ) 
= \Bigl( \xi + \vert \rho \vert \Bigr) \Bigl( \tilde f(\xi) + \tilde f (\vert \rho \vert ) \Bigr)
\\
\geq \Bigl( \xi + \vert \rho \vert \Bigr) \Bigl( \tilde f'(1) \xi + \tilde f'(1) \vert \rho \vert  \Bigr) = \tilde f'(1) \Bigl( \xi + \vert \rho \vert \Bigr)^2 \geq \frac{\tilde f'(1)}{C_p} \tilde F(\xi - \rho) .
\end{multline*}
If $\xi-\rho \geq 1$, then 
$$ \frac{\tilde F(\xi - \rho)}{\xi - \rho} \leq C_p (\xi-\rho) ^p ,$$
while
$$ \tilde f(\xi) - \tilde f (\rho) = \tilde f(\xi) + \tilde f (\vert \rho \vert )
\geq c_p \xi ^{p-1} + c_p \vert \rho \vert ^{p-1},$$
and once again there is a uniform $C$ such that
$$ \tilde F(\xi - \rho) \leq C (\xi - \rho) (\tilde f(\xi) - \tilde f (\rho) ) .$$
This concludes the proof of Lemma \ref{lem-mai7}. \qed

\subsubsection{Proof of Lemma \ref{lem-1st-estim*}} \hfill

The modified first multiplier gives immediately an identity completely similar to
\eqref{mult-1st-mult} (with $\tilde f$ and $\tilde E_p$ instead of $f$ and $E_p$), and then we easily obtain \eqref{mult-1st-mult-cons*}. \qed


\subsubsection{Proof of Lemma \ref{lem-2nd-estim*}} \hfill

The strategy is of course to follow the proof of Lemma \ref{lem-2nd-estim},
but there are some modifications that we have to explain.

The starting point is the modified \eqref{mult-2nd-mult}, modified in the sense that $f$ and $E_p$ are replaced by $\tilde f$ and $\tilde E_p$. 
We follow the estimates in the same order as in section \ref{sec-sec-mult}, in the proof
of Lemma \ref{lem-2nd-estim}.

{\underline{Modifications of \eqref{estim-2.1} and \eqref{estim-2.2} when $p\in (1,2)$:}}
First
\begin{multline*}
\Bigl\vert  \int _0 ^1 \psi _2 (x) z \Bigl( \tilde f(\xi) - \tilde f(\rho)\Bigr) \, dx \Bigr\vert
\leq  \int _0 ^1 \Bigl\vert z \tilde f(\xi) \Bigr\vert +  \Bigl\vert z \tilde f(\rho)\Bigr) \Bigr\vert \, dx
\\
\leq C  \Bigl(  \int _0 ^1 \tilde F (z) \, dx + \int _0 ^1 \tilde F (\xi)  + \tilde F (\rho) \, dx \Bigr)
.
\end{multline*}
But $z_x= \frac{\rho}{2} + \frac{\xi}{2}$, therefore, since $\tilde F$ is convex, we have
$$ \tilde F (z_x) = \tilde F \Bigl( \frac{\rho}{2} + \frac{\xi}{2} \Bigr) \leq 
\frac{1}{2} \tilde F (\rho) + \frac{1}{2} \tilde F (\xi) . $$
Therefore, there exists $C>0$ such that
\begin{equation}
\label{2ndmult-mod-0}
\Bigl\vert \int _0 ^1 \psi _2 (x) z \Bigl( \tilde f(\xi) - \tilde f(\rho)\Bigr) \, dx \Bigr\vert
\leq C \tilde E_p (t) .
\end{equation}
This directly implies that
\begin{equation}
\label{2ndmult-mod-1}
\Bigl\vert \Bigl[ \tilde E_p (t) \phi ' (t) \int _0 ^1 \psi _2 (x) z \Bigl( \tilde f(\xi) - \tilde f(\rho)\Bigr) \, dx \Bigr] _S ^T  \Bigr\vert
\leq C \tilde E_p (S) ^2 \phi ' (S),
\end{equation}
and 
\begin{equation}
\label{2ndmult-mod-2}
\Bigl\vert \int _S ^T (\tilde E_p (t) \phi '(t) )' \int _0 ^1 \psi _2 (x) z \Bigl( \tilde f(\rho) - \tilde f(\xi)\Bigr) \, dx \, dt  \Bigr\vert
\leq C \tilde E_p (S) ^2 \phi ' (S),
\end{equation}
hence the modifications of \eqref{estim-2.1} and \eqref{estim-2.2} when $p\in (1,2)$.

{\underline{Modification of \eqref{estim-2.3}:}} using \eqref{Fenchel-cons<}, we have
\begin{multline*}
\Bigl\vert - \int _S ^T \tilde E_p (t) \phi ' (t) \int _0 ^1 \psi _2 ' (x)  \, z \Bigl( \tilde f(\rho) + \tilde f(\xi)\Bigr) \, dx \, dt \Bigr\vert
\\
\leq \int _S ^T \tilde E_p (t) \phi ' (t) \int _{a_2} ^{b_2}   \vert z \tilde f(\rho) \vert + \vert z  \tilde f(\xi) \vert \, dx \, dt
\\
\leq \int _S ^T \tilde E_p (t) \phi ' (t) \int _{a_2} ^{b_2}
\frac{C_p}{\eta ^2} \tilde F(z) + C_p \eta ^2 (\tilde F(\rho) + \tilde F(\xi)) \, dx \, dt .
\end{multline*}
Therefore we obtain the modification of \eqref{estim-2.3}:
\begin{multline}
\label{2ndmult-mod-3}
\Bigl\vert  - \int _S ^T \tilde E_p (t) \phi ' (t) \int _0 ^1 \psi _2 ' (x)  \, z \Bigl( \tilde f(\rho) + \tilde f(\xi)\Bigr) \, dx \, dt \Bigr\vert
\\
\leq \frac{C_p}{\eta ^2} \int _S ^T \tilde E_p (t) \phi ' (t) \Bigl( \int _{a_2} ^{b_2}
 \tilde F(z) \, dx \Bigr) \, dt + C_p \eta ^2 \int _S ^T \tilde E_p (t)^2 \phi ' (t) \, dt .
\end{multline}

{\underline{Modification of \eqref{estim-2.4}:}} since $g$ is globally Lipschitz and $\tilde f'$ is bounded, we have
\begin{multline}
\label{estim-2.4-mai1}
\Bigl\vert - \int _S ^T \tilde E_p (t) \phi ' (t) \int _0 ^1 \psi _2 (x)   a(x)  \, z  g(\frac{\rho - \xi}{2}) \Bigl( \tilde f'(\rho) + \tilde f'(\xi)\Bigr) \, dx \, dt \Bigr \vert
\\
\leq C \int _S ^T \tilde E_p (t) \phi ' (t) \int _{a_2} ^{b_2} \vert z \vert \, \Bigl\vert \frac{\rho - \xi}{2} \Bigr \vert .
\end{multline}
On the other hand, the function $z$ belongs to $L^\infty (\Bbb R_+, W^{1,p} _0 (0,1))$, and
hence to $L^\infty (\Bbb R_+, L^\infty (0,1))$. We can even provide a more precise estimate:
as noted before, since $z(t,0)=0$, we have
\begin{multline*}
\vert z(t,x) \vert ^p  
= \Bigl \vert \int _0 ^x  z_x (t,s) \, ds \Bigr\vert ^p
\leq \Bigl( \frac{1}{2} \int _0 ^1 \vert \rho \vert \, ds + \frac{1}{2} \int _0 ^1 \vert \xi \vert \, ds \Bigr) ^p 
\\
\leq \frac{1}{2} \int _0 ^1 \vert \rho \vert ^p \, ds + \frac{1}{2} \int _0 ^1 \vert \xi \vert ^p \, ds = \frac{p}{2} E_p (t) \leq \frac{p}{2} E_p (0) . 
\end{multline*}
Therefore
$$ \forall t \geq 0, \forall x \in [0,1], \quad \vert z(t,x) \vert \leq C_p ^{(0)} \quad \text{ where } \quad  C_p ^{(0)} := \Bigl( \frac{p}{2} E_p (0) \Bigr) ^{1/p} .$$
Since $\tilde f$ is increasing and $\tilde f '(0) = p-1 >0$, there exists $K^{(0)}$ such that
$$ \forall y \in [- C_p ^{(0)}, C_p ^{(0)}], \quad \vert y \vert  \leq K^{(0)} \tilde f ( \vert y \vert ) .$$
(In fact, since $\tilde f$ is concave, one can choose
$$ K^{(0)} = \frac{C_p ^{(0)}}{\tilde f (C_p ^{(0)})} ,$$
and then there is some $C_p > C_p >0$ independent of $E_p (0)$ such that
$$ C_p ' \max (1, E_p (0) ^{\frac{2-p}{p}} \leq K^{(0)} \leq C_p  \max (1, E_p (0) ^{\frac{2-p}{p}} .)$$
Therefore, 
\begin{equation}
\label{2ndmult-mod-4-inter2}
\forall t \geq 0, \forall x \in [0,1], \quad \vert z (t,x) \vert \leq K^{(0)} \vert \tilde f (z (t,x)) \vert .
\end{equation}
Now we are in position to obtain the extension of \eqref{estim-2.4}: combining \eqref{2ndmult-mod-4-inter2} with \eqref{estim-2.4-mai1}, we see that
$$
\int _{a_2} ^{b_2} \vert z \vert \, \Bigl\vert \frac{\rho - \xi}{2} \Bigr \vert 
\leq K^{(0)} \int _{a_2} ^{b_2} \Bigl( \vert \rho \tilde f (z (t,x)) \vert \Bigr) + \Bigl( \vert \xi \tilde f (z (t,x)) \vert \Bigr) ,
$$
and using twice \eqref{Fenchel-cons>} (with $a=\rho$ and $a=\xi$), we obtain
\begin{multline*}
\int _{a_2} ^{b_2} \vert z \vert \, \Bigl\vert \frac{\rho - \xi}{2} \Bigr \vert 
\leq K^{(0)} \Bigl( \int _{a_2} ^{b_2} \Bigl(  C_p \eta ^p \tilde F (\rho) + \frac{C_p}{ \eta ^q} \tilde F (z) \Bigr) + \Bigl(  C_p \eta ^p \tilde F (\xi) + \frac{C_p}{ \eta ^q} \tilde F (z) \Bigr)
\\
= 2 \frac{K^{(0)} \, C_p}{ \eta ^q} \int _{a_2} ^{b_2} \tilde F (z)
+ K^{(0)}  C_p \eta ^p \int _{a_2} ^{b_2} \tilde F (\rho) + \tilde F (\xi) .
\end{multline*}
Therefore we obtain the extension of \eqref{estim-2.4}:
\begin{multline}
\label{2ndmult-mod-3*}
\Bigl\vert - \int _S ^T \tilde E_p (t) \phi ' (t) \int _0 ^1 \psi _2 (x)   a(x)  \, z  g(\frac{\rho - \xi}{2}) \Bigl( \tilde f'(\rho) + \tilde f'(\xi)\Bigr) \, dx \, dt \Bigr \vert
\\
\leq  \frac{K^{(0)} \, C}{ \eta ^q} \int _S ^T \tilde E_p (t) \phi ' (t) \Bigl( \int _{a_2} ^{b_2} \tilde F (z) \, dx \Bigr) 
+ K^{(0)}  C_p \eta ^p \int _S ^T \tilde E_p (t) ^2 \phi ' (t) \, dt .
\end{multline}

We deduce the modified version of \eqref{mult-2nd-mult-cons1}:
\begin{multline}
\label{mult-2nd-mult-cons1*}
\int _S ^T \tilde E_p (t) \phi ' (t) \int _0 ^1 \psi _2 (x) \Bigl( \rho \tilde f(\rho) + \xi \tilde f(\xi) \Bigr) \, dx \, dt
\\
\leq C \tilde E_p (S) ^2 \phi ' (S) 
+ K^{(0)}  C \eta ^p \int _S ^T \tilde E_p ^2 \phi ' \, dt
 + \frac{K^{(0)} \, C}{ \eta ^q} \int _S ^T \tilde E_p \phi ' \Bigl( \int _{a_2} ^{b_2} \tilde F (z) \, dx \Bigr) \, dt
\\
+ 2 \int _S ^T \tilde E_p (t) \phi ' (t) \int _0 ^1 \psi _2 (x) \, \frac{\rho - \xi}{2} \Bigl( \tilde f(\rho) - \tilde f(\xi)\Bigr) \, dx \, dt ,
\end{multline}
wich immediately gives \eqref{mult-2nd-mult-cons*}. \qed


\subsubsection{Proof of Lemmas \ref{lem-3rd-estim*} and \ref{lem-mai2}} \hfill

\noindent {\it Proof of Lemma \ref{lem-mai2}.}
First, we have an integral formula of $\tilde v$:
$$ \tilde v(t,x)= -x \int _x ^1 (1-y) \psi _3 (y) \tilde f(z(t,y)) \, dy + (x-1) \int _0 ^x y \psi _3 (y) \tilde f(z(t,y)) \, dy ,$$
and, since $\tilde F ^*$ is convex, we derive that
\begin{multline*}
\tilde F ^* (\tilde v(t,x)) \leq x \tilde F ^* \Bigl( - \int _x ^1 (1-y) \psi _3 (y) \tilde f(z(t,y)) \, dy \Bigr) + (1-x) \tilde F ^* \Bigl( - \int _0 ^x y \psi _3 (y) \tilde f(z(t,y)) \, dy \Bigr)
\\
\leq \int _x ^1 \tilde F ^* \Bigl( (1-y) \psi _3 (y) \tilde f(z(t,y))  \Bigr) \, dy
+ \int _0 ^x \tilde F ^* \Bigl( y \psi _3 (y) \tilde f(z(t,y)) \Bigr) \, dy
\\
\leq 2 \int _0 ^1 \tilde F ^* \Bigl( \tilde f(z(t,y))  \Bigr)\, dy .
\end{multline*}
Since $\tilde f$ is concave, we have
$$ \vert \tilde f(z(t,y)) \vert \leq (p-1) \vert z(t,y) \vert ,$$
and therefore by monotonicity of $\tilde F^*$:
$$ \tilde F ^* \Bigl( \tilde f(z(t,y))  \Bigr)
\leq \tilde F ^* \Bigl( (p-1) z(t,y)  \Bigr),$$
and using \eqref{comp-tildeF*}
\begin{multline*}
\tilde F ^* \Bigl( (p-1) z(t,y)  \Bigr) 
\leq C_p ^* \max ((p-1)^2 z(t,y) ^2 , (p-1)^q \vert z(t,y) \vert ^q) 
\\
\leq C' \max (z(t,y) ^2 , \vert z(t,y) \vert ^q) 
\leq C' \max (z(t,y) ^2 , {C_p ^{(0)}} ^{q-2}  z(t,y)^2) 
\\
= C' z(t,y) ^2 \max (1, {C_p ^{(0)}} ^{q-2}) .
\end{multline*}
Using \eqref{comp-tildeF}, one easily see that
$$ \forall y \in [- C_p ^{(0)}, C_p ^{(0)}], \quad 
y^2 \leq \frac{\max(1, {C_p ^{(0)}}^{2-p})}{C_p '} \tilde F(y) ,$$
and therefore 
$$ \tilde F ^* \Bigl( (p-1) z(t,y)  \Bigr)
\leq C' \max (1, {C_p ^{(0)}} ^{q-2}) \frac{\max(1, {C_p ^{(0)}}^{2-p})}{C_p '} \tilde F(z(t,x)) .$$
Therefore
$$ \tilde F ^* (\tilde v(t,x)) \leq 2 C' \max (1, {C_p ^{(0)}} ^{q-2}) \frac{\max(1, {C_p ^{(0)}}^{2-p})}{C_p '} \int _0 ^1 \tilde F(z ) \, dx .$$
We conclude using the generalized Poincar\'e's inequality \eqref{eq-Hardy-Fenchel}
that there is $C''$ independent of the initial conditions such that
$$ \tilde F ^* (\tilde v(t,x)) \leq C'' \max (1, {C_p ^{(0)}} ^{q-2}) \, \max(1, {C_p ^{(0)}}^{2-p}) \, \int _0 ^1 \tilde F(\rho ) + \tilde F (\xi ) \, dx .$$
Using the expression of $C_p ^{(0)}$, we see that
$$ \max (1, {C_p ^{(0)}} ^{q-2}) \, \max(1, {C_p ^{(0)}}^{2-p}) \leq C \max (1, E_p (0) ^{\frac{2-p}{p-1}} ),$$
which gives \eqref{estim-v-Lq*} (integrating w.r.t. $x$).  

To prove \eqref{estim-v_t-Lq*}, we start from the integral formula of $\tilde v_t$:
$$ \tilde v_t (t,x)= -x \int _x ^1 (1-y) \psi _3 (y) \tilde f'(z(t,y)) z_t (t,y) \, dy + (x-1) \int _0 ^x y \psi _3 (y) \tilde f'(z(t,y)) z_t (t,y) \, dy .$$
Then
\begin{multline*}
 \tilde F ^* (\tilde v_t) \leq \tilde F ^* \Bigl( \int _0 ^1 \psi _3 (y) \vert z_t (t,y)\vert \, dy  \Bigr) = \tilde F ^* \Bigl( \frac{1}{2} \int _{\vert z_t \vert \leq 1} 2 \psi _3 \vert z_t \vert + \frac{1}{2} \int _{\vert z_t \vert \geq 1} 2 \psi _3 \vert z_t \vert \Bigr) .
\\
\leq \frac{1}{2} \tilde F ^* \Bigl( \int _{\vert z_t \vert \leq 1} 2 \psi _3 \vert z_t \vert \Bigr) + \frac{1}{2} \tilde F ^* \Bigl( \int _{\vert z_t \vert \geq 1} 2 \psi _3 \vert z_t \vert \Bigr) .
\end{multline*}
Since there exists $c_0 >0$ such that
$$ \forall y \in [0, 1], \quad \tilde f (y) \geq c_0 y ,$$
we deduce that
$$ \int _{\vert z_t \vert \leq 1} 2 \psi _3 \vert z_t \vert
\leq \frac{2}{c_0} \int _{\vert z_t \vert \leq 1} \psi _3 \tilde f(\vert z_t \vert ),
$$
and therefore, 
$$ \tilde F ^* \Bigl( \int _{\vert z_t \vert \leq 1} 2 \psi _3 \vert z_t \vert \Bigr)
\leq C \int _0 ^1 \tilde F^* (\psi _3 \tilde f(\vert z_t \vert))
\leq C' \int _0 ^1 \psi _3  \tilde F ( \tilde f(\vert z_t \vert) ) 
\leq C'' \int _0 ^1 \psi _3  \tilde F (  z_t )  $$ 
as we have seen previously, and that is what we want.

On the other hand, 
$$ \int _{\vert z_t \vert \geq 1} 2 \psi _3 \vert z_t \vert 
\leq \int _0 ^1 \vert \rho \vert + \vert \xi \vert 
\leq  \Bigl( \int _0 ^1 \vert \rho \vert ^p \Bigr)^{1/p} + \Bigl( \int _0 ^1 \vert \xi \vert ^p \Bigr)^{1/p} \leq C_p E_p (0) ^{1/p}  ,$$
and then it follows from \eqref{comp-tildeF*} that
\begin{multline*}
\forall b \in [0, C_p E_p (0) ^{1/p}], \quad
\tilde F ^* (b) \leq C_p ^* \max (b^2, b ^{\frac{p}{p-1}} )
\\
\leq C_p ^* \max (b^2, (C_p E_p (0) ^{1/p}) ^{q-2} b ^{2} )
= \Bigl( C_p ^* \max (1, (C_p E_p (0) ^{1/p}) ^{q-2} ) \Bigr) \, b^2 .
\end{multline*}
Therefore,
$$ \tilde F ^* \Bigl( \int _{\vert z_t \vert \geq 1} 2 \psi _3 \vert z_t \vert \Bigr)
\leq C_p ^* \max (1, (C_p E_p (0) ^{1/p}) ^{q-2} ) 
\Bigl( \int _{\vert z_t \vert \geq 1} 2 \psi _3 \vert z_t \vert \Bigr) ^2 .$$
Now it is sufficient to note that there is some $c_1 >0$ such that
$$ \forall y \geq 1, \quad \tilde f (y) \geq c_1 y^{p-1} .$$
This yields that
$$ \vert z_t \vert \geq 1 \quad \implies \quad \vert z_t \vert \leq \Bigl( \frac{\tilde f (\vert z_t \vert)}{c_1} \Bigr) ^{1/p-1} .$$
Then
\begin{multline*}
\int _{\vert z_t \vert \geq 1} 2 \psi _3 \vert z_t \vert 
\leq \int _{\vert z_t \vert \geq 1} 2 \psi _3  \Bigl( \frac{\tilde f (\vert z_t \vert)}{c_1} \Bigr) ^{1/p-1}
\leq \int _0 ^1 2 \psi _3  \Bigl( \frac{\tilde f (\vert z_t \vert)}{c_1} \Bigr) ^{1/p-1}
\\
\leq C \Bigl( \int _0 ^1  \psi _3  \tilde f (\vert z_t \vert) ^{p/p-1} \Bigr) ^{1/p}
\leq C' \Bigl( \int _0 ^1  \psi _3  \tilde F (\vert z_t \vert) \Bigr) ^{1/p} .
\end{multline*}
We can conclude that
\begin{multline*}
\tilde F ^* \Bigl( \int _{\vert z_t \vert \geq 1} 2 \psi _3 \vert z_t \vert \Bigr)
\leq C_p ^* \max (1, (C_p E_p (0) ^{1/p}) ^{q-2} ) \Bigl( \int _0 ^1  \psi _3  \tilde F ( z_t ) \Bigr) ^{2/p}
\\
\leq C_p ^* \max (1, (C_p E_p (0) ^{1/p}) ^{q-2} ) \Bigl( \tilde E_p (0) \Bigr) ^{\frac{2}{p}-1} \, \int _0 ^1  \psi _3  \tilde F ( z_t ), 
\end{multline*}
and finally, we get
\begin{multline*}
\tilde F ^* (\tilde v_t) \leq C_p ^* \max (1, (C_p E_p (0) ^{1/p}) ^{q-2} ) \Bigl( \tilde E_p (0) \Bigr) ^{\frac{2}{p}-1} \, \int _0 ^1  \psi _3  \tilde F ( z_t )
\\
\leq C \max (E_p (0) ^{\frac{2-p}{p}},E_p (0) ^{\frac{2-p}{p-1}}) \int _0 ^1  \psi _3  \tilde F ( z_t ) . \end{multline*}
This concludes the proof of Lemma \ref{lem-mai2}. \qed

\medskip

\noindent {\it Proof of Lemma \ref{lem-3rd-estim*}.} 
First, as we proved \eqref{mult-3rd-mult}, we have now:
\begin{multline}
\label{mult-3rd-mult*}
2 \int _S ^T \tilde E_p \phi ' \int _0 ^1 \psi _3 (x) z \tilde f (z) 
= \Bigl[ \tilde E_p \phi ' \int _0 ^1 \tilde v (\rho - \xi) \Bigr] _S ^T 
\\
- \int _S ^T (\tilde E_p \phi ')' \int _0 ^1  \tilde v (\rho - \xi)
-  \int _S ^T \tilde E_p \phi ' \int _0 ^1  \tilde v_t  (\rho - \xi)
\\
+ 2 \int _S ^T \tilde E_p \phi ' \int _0 ^1 a (x) \, \tilde v g(\frac{\rho - \xi}{2}) .
\end{multline}
Note that the function $y\mapsto y \tilde f(y)$ is even, and satisfies
$$ y \tilde f(y) \sim _{y\to 0} (p-1) x^2, \quad \text{ and } \quad y \tilde f(y) \sim _{y\to +\infty} y^p,$$
there exists some $c_p > c_p ' >0$ such that
$$ \forall y \in \Bbb R, \quad c_p ' \tilde F (y)  \leq y \tilde f(y) \leq c_p \tilde F (y).$$
Therefore, using Lemma \ref{lem-mai2} and \eqref{mult-2nd-mult-cons*}, 
we see that the second term that appears in \eqref{mult-2nd-mult-cons*}
can be estimated as follows:
$$ \int _S ^T \tilde E_p (t) \phi '(t)  \Bigl( \int _{a_2} ^{b_2}  \tilde F (z) \, dx \Bigr) \, dt \leq \frac{1}{c_p '} \int _S ^T \tilde E_p (t) \phi '(t) \Bigl( \int _{0} ^{1} \psi _3 (x)  z \tilde f (z) \, dx \Bigr) \, dt ,$$
that is estimated using the identity \eqref{mult-3rd-mult*} and the estimates of Lemma \ref{lem-mai2}. First, using \eqref{Fenchel} and \eqref{estim-v-Lq*}, we have
$$ \Bigl\vert \int _0 ^1 \tilde v (\rho - \xi)\Bigr\vert
\leq \int _0 ^1 \tilde F^* (\tilde v) + \tilde F (\rho - \xi) 
\leq C \max (1, E_p (0) ^{\frac{2-p}{p-1}} ) \tilde E_p (t) ,$$ 
and therefore
$$ \Bigl\vert \Bigl[ \tilde E_p \phi ' \int _0 ^1 \tilde v (\rho - \xi) \Bigr] _S ^T 
- \int _S ^T (\tilde E_p \phi ')' \int _0 ^1  \tilde v (\rho - \xi) \Bigr\vert \leq C \max (1, E_p (0) ^{\frac{2-p}{p-1}} ) \tilde E_p (S) ^2 \phi '(S) .$$
Next, once again using \eqref{Fenchel}, \eqref{comp-tildeF}, \eqref{comp-tildeF*}, and \eqref{estim-v_t-Lq*}, we have: given $\sigma \in (0,1)$, 
\begin{multline*}
\Bigl\vert \int _0 ^1  \tilde v_t  (\rho - \xi) \Bigr\vert 
= \Bigl\vert \int _0 ^1  \frac{\tilde v_t}{\sigma} \, \Bigl( \sigma  (\rho - \xi) \Bigr)\Bigr\vert 
\\
\leq \int _0 ^1 \tilde F ^* \Bigl( \frac{\tilde v_t}{\sigma} \Bigr) + \tilde F \Bigl( \sigma  (\rho - \xi) \Bigr)
\leq \frac{C}{\sigma ^q} \int _0 ^1 \tilde F ^* (\tilde v_t)
+ C \sigma ^p  \int _0 ^1 \tilde F (\rho) + \tilde F (\xi) 
\\
\leq C \sigma ^p \tilde E _p (t) + \frac{C}{\sigma ^q} \max (E_p (0) ^{\frac{2-p}{p}},E_p (0) ^{\frac{2-p}{p-1}}) \int _0 ^1  \psi _3  \tilde F ( z_t ) \, dx . 
\end{multline*} 
The last estimates give \eqref{mult-3rd-mult-cons*}. \qed



\end{document}